\documentclass[11pt]{article}
\usepackage[a4paper,margin=2.4cm]{geometry}
\usepackage{amsmath,amssymb,amsthm,booktabs,enumitem,array,float}
\usepackage{mathrsfs}
\usepackage{xcolor}
\usepackage[colorlinks,linkcolor=black,citecolor=blue!45!black,urlcolor=blue!50!black]{hyperref}
\usepackage{graphicx}
\usepackage{multirow,adjustbox}
\usepackage{algorithm}
\usepackage{algpseudocode}
\newtheorem{definition}{Definition}
\usepackage{tikz}\usetikzlibrary{positioning,arrows.meta,calc}
\usepackage[affil-it]{authblk}
\usepackage{fancyhdr}
\newcommand{\bd}{\partial\Omega}\newcommand{\dn}{\partial_n}\newcommand{\GD}{\Gamma_D}\newcommand{\GN}{\Gamma_N}

\newcommand{\Om}{\Omega}
\newcommand{\R}{\mathbb R}\newcommand{\qint}{\widetilde{\int}}\newcommand{\pre}{\mathrm{pre}}\newcommand{\un}{\mathrm{un}}
\theoremstyle{plain}\newtheorem{theorem}{Theorem}[section]\newtheorem{proposition}[theorem]{Proposition}
\newtheorem{lemma}[theorem]{Lemma}\newtheorem{corollary}[theorem]{Corollary}
\theoremstyle{definition}
\theoremstyle{plain}\newtheorem{remark}[theorem]{Remark}
\title{\vspace{-1.5cm}\textsc{FokasNets: Completing missing boundary data in Physics-Informed Neural Networks and Operators via the Global Relation}}
\author[1]{Kyriakos C. Georgiou}
\author[1]{Khemraj Shukla}
\author[2,3]{Athanassios S. Fokas}
\author[1]{George Em Karniadakis\thanks{Corresponding author: \texttt{george\_karniadakis@brown.edu}}}
\affil[1]{Division of Applied Mathematics, Brown University, Providence, RI, USA}
\affil[2]{Department of Applied Mathematics and Theoretical Physics, University of Cambridge, Cambridge, United Kingdom}
\affil[3]{Mathematics Research Center, Academy of Athens, Athens, Greece}
\date{}

\begin{document}
\maketitle
\vspace{-0.8cm}

\begin{abstract}\noindent
We introduce FokasNets, a framework that combines Physics-Informed Neural Networks (PINNs) and operators (PINOs) with the Global Relation (GR) formulation of the Fokas method (also referred to as the Unified Transform), for mixed boundary problems. The general framework considers lifted representations of the neural network models, where the Dirichlet boundary conditions are imposed by construction and the neural network learns the interior correction. The Global Relation is used as an additional solve during training that provides the missing boundary conditions to the lifted neural network, making training more efficient and accurate in the interior and on the boundary. 
We develop the construction of FokasNets for elliptic and evolution PDEs (Fokas-PINNs), and the corresponding PDE solution operator problems (Fokas Operator Networks - FONets), using the Global Relation formulations for these PDE classes. Furthermore, we develop the necessary theoretical framework for FokasNets, including error and stability analysis and show that FokasNets are universal approximators. We provide comprehensive numerical experiments, using examples with manufactured solutions and literature benchmarks. Compared to standard PINNs and Physics-Informed DeepONets and their lifted versions, FokasNets achieve significant improvements in solution accuracy in the interior and on the boundaries, both for single PDEs and particularly for solution operators, where the improvement reaches multiple orders of magnitude and fixed threshold errors are reached faster during training. 
\end{abstract}


\section{Introduction}\label{sec:introduction}

In recent years, the use of neural network models for approximating the solution of Partial Differential Equations (PDEs) has given rise to the field of Scientific Machine Learning (SciML), producing a vast amount of scientific literature and progress in fields ranging from physics to medicine and many more. The key breakthrough that ushered in this era was Physics-Informed Neural Networks (PINNs) \cite{raissi2019physics}, where the solution is modelled by a neural network and training is then performed by minimizing a loss function that encapsulates the PDE, boundary and initial condition residuals. Hence, PINNs are able to provide solutions only using the governing equations without the need for data. This approach has been incredibly successful, having been used in problems in physics \cite{raissi2020hidden,mao2020physics,cai2021physics}, biomedical engineering \cite{kissas2020machine,arzani2021uncovering}, medicine \cite{sahli2020physics,herrero2022ep,kharazmi2021identifiability} and many more. Extending this approach, the next important milestone was the development of the Deep Operator Network (DeepONet) \cite{lu2021learning}, which is able to learn the solution operator, i.e., a family of PDE solutions. Although first introduced using data-driven training, training can also incorporate the equation formulations, when known, giving rise to Physics-Informed Neural Operators (PINOs) \cite{wang2021learning,goswami2023physics}, and the corresponding Physics-Informed Deep Operator Networks (PI-DeepONets). Applications of operator learning and PINOs have been equally broad and fruitful \cite{goswami2022physics,wang2023long,mao2021deepm,oommen2022learning,lu2022multifidelity}. 

An important issue that often arises in PINNs and PINOs, and is one of the motivations of the present work, is the handling and effect of boundary and initial conditions in the Physics-Informed training. Specifically, as these conditions enter the loss function as additional penalty terms, their imposition is ``soft'' and the resulting accuracy of the neural network approximation on the domain depends on optimization. Furthermore, practical details such as the scales of the loss terms and sampling during training can result in issues such as higher errors in the boundary conditions than in the interior, which in turn hinder the accuracy of the interior solution \cite{wang2021understanding}. An established remedy is to encode the boundary information into the approximation function itself: a prescribed component satisfies the boundary conditions by construction and a trainable correction, multiplied by a function vanishing on the boundary, preserves them \cite{lagaris1998artificial}. This idea has been systematized through the Theory of Functional Connections (TFC) \cite{mortari2017theory,mortari2019multivariate, leake2022theory}, which parametrizes the entire set of functions satisfying a given family of linear constraints, and has been combined with neural networks \cite{leake2020deep} to yield hard-constrained PINNs. Constructions based on approximate distance functions and transfinite interpolation \cite{sukumar2022exact,dong2026novel} extend exact imposition to Dirichlet, Neumann and Robin conditions on general geometries. However, as the constrained expression must be derived anew for each combination of constraint type, PDE class and geometry, the resulting formulas quickly become involved and require additional derivations. In practice, therefore, hard-constrained PINNs most often retain the simplest form, in which only the Dirichlet data are encoded, through a lift that reproduces the trace values and extends smoothly into the interior, while the remaining conditions are left to the loss. As will be described, we deliberately adopt this simplest form in the present work, with the lift realized by standard extension operators (harmonic extension, Coons transfinite interpolation \cite{coons1967surfaces,gordon1973construction} or Hermite interpolation, depending on the geometry), and address its main limitation by a different mechanism. Finally, the enforcement of boundary conditions in operator learning is significantly less developed. In particular, exact enforcement has been limited mainly to Dirichlet and periodic boundary conditions, which have been built into DeepONets \cite{lu2022comprehensive,brecht2023improving}, into the kernel of a neural operator by a boundary correction, which also imposes Neumann data up to a finite-difference approximation \cite{saad2022guiding}, and, more recently, into every layer of a kernel-integral operator by restricting its output to the span of Dirichlet Laplacian eigenfunctions \cite{stuart2026enforcing}. 

\par The natural question is how to handle mixed boundary conditions, i.e., when the Dirichlet data are missing or incomplete. The missing values are then part of the solution to be approximated. Leaving these values undetermined means that they become additional trainable terms, leaving their determination to the training process, with the issues described above. The Fokas method and the corresponding Global Relation, on the other hand, provide a deterministic, closed-form relation that couples the two types of conditions on the entire boundary, which allows us to determine the missing trace needed for training.

The Global Relation (GR) is a central ingredient of the Fokas method, or
unified transform, introduced for linear and certain integrable nonlinear
PDEs \cite{fokas1997unified,fokas2012synthesis}. For linear
elliptic equations, the GR is obtained by using suitable solutions of the adjoint equation to yield an identity that connects the Dirichlet and Neumann traces. The full Fokas method then uses the GR together with spectral representations to characterize unknown
boundary values and reconstruct the solution. Numerical implementations of the method recover missing traces through boundary expansions and spectral collocation
\cite{fokas2009semi,fornberg2011numerical} and concrete strategies for the configurations of these numerical methods have been developed for elliptic problems on polygonal domains \cite{hashemzadeh2015numerical,colbrook2019hybrid} and on curvilinear and circular domains \cite{colbrook2020extending,crowdy2015transform}. Furthermore, analytical work has also established stability properties and convergent Galerkin formulations for the underlying spectral Dirichlet-to-Neumann map \cite{ashton2013spectral,ashton2016numerical}. The method has been developed for linear elliptic equations, including the Laplace, modified Helmholtz and Helmholtz equations \cite{colbrook2019unified,fulton2004analytical,davis2014spectrally} and the biharmonic equation \cite{crowdy2004explicit,dimakos2015poisson}, for linear evolution equations on the half-line and the finite interval, such as the heat, linear Schr\"odinger and linearized KdV equations \cite{fokas2002new,fokas2001two,deconinck2014method,kesici2018numerical}, and for integrable nonlinear equations such as the nonlinear Schr\"odinger, KdV and sine--Gordon equations on the half-line \cite{fokas2002integrable,de2019hybrid,fokas2012unified,pelloni2015advances}. In the integrable case the GR is obtained from the associated Lax pair and remains closed in the spectral functions. However, for nonlinear equations without this structure, the GR depends on the solution of the PDE in the interior and cannot be closed by analytical means alone, and for this reason the Fokas method has not been extended to such problems.

In this work, we combine the complementary advantages of the PINN/PINO approach to learn nonlinear dynamics based on the PDE formulations with the power of the Fokas method to address the issues explained above, by developing FokasNets. FokasNets are a PINN and PINO framework, whereby training is performed by combining the GR with the lifted neural network representation that enforces the Dirichlet trace. In general, during training, the GR provides the Dirichlet values on the entire boundary, which are then used to construct the lifted neural network. Hence, the GR becomes an additional step in the optimization algorithm, used to complete and constrain the missing boundary trace. In this way, the two components have distinct roles: the GR enforces nonlocal compatibility
conditions between the learned and prescribed boundary data, and the Physics-Informed training learns the neural network approximation that satisfies the PDE pointwise. 
We call the resulting family of models trained in this way \emph{FokasNets}: when combined with PINNs for solving single PDEs we use the term \emph{Fokas-PINN} and when combined with PINOs for approximating the solution operator of a PDE family, \emph{Fokas-Operator Network (FONet)}. 

Generally, hybrid schemes that embed numerical solvers in neural network architectures are an established approach in SciML. One such research direction considers hidden layer outputs as the solution of an equation
and differentiates it through that equation, as in e.g., differentiable
optimization layers \cite{amos2017optnet,agrawal2019differentiable}, deep
equilibrium models \cite{bai2019deep}, adjoint-based neural ODEs
\cite{chen2018neural}, or Fredholm Neural Networks \cite{georgiou2025fredholm,georgiou2025fredholmelliptic, georgiou2026learning}, which instead unroll layers as consecutive iterations in a fixed point scheme. A general framework for the implicit differentiation and training used in such methods is given in \cite{blondel2022efficient}. Differentiable PDE solvers have been placed inside training loops to correct or accelerate simulations \cite{um2020solver,kochkov2021machine,belbute2020combining,patsatzis2025gorinns}, and differentiable PDE-constrained layers have been used to enforce physical constraints in
neural representations \cite{negiar2022learning}. Other methods rely on underlying mathematical theory that allows improved training efficiency, such as e.g., variable projection
\cite{golub1973differentiation,newman2021train,georgiou2025heatnets,fabiani2025random,fabiani2025randonets}.

The present work adds to this line of research with the following key contributions:
\begin{enumerate}
    \item We develop FokasNets by targeting a fundamental difficulty of Physics-Informed training (namely, the handling of boundary and initial conditions) and embedding the main component of the highly accurate and efficient Fokas method to tackle precisely this difficulty. In this way, FokasNets provide a principled mathematical method to tackle the difficulties that may arise in PINN and PINO training.

    \item We detail the FokasNet constructions for elliptic and evolution PDEs, using the corresponding GR representations and provide the necessary discretization schemes and matrix formulations of the resulting finite Global Relation, as well as the mechanisms for the lifted neural network representation and for the propagation of gradients through the GR during Physics-Informed training.

    \item We prove important theoretical results further motivating the efficacy of FokasNets. In particular, we prove that FokasNets are universal approximators and further establish error analysis and stability results. 

    \item We provide in-depth numerical results across a suite of PDE and solution operator examples, showing that FokasNets can achieve significantly higher accuracies compared to standard neural network approximations and their lifted counterparts, across both the interior and boundaries (the improvement is of multiple orders of magnitude in some cases). This is particularly evident in the solution operator approximations with FONets, where imposing exact boundary conditions is considerably harder than the single PDE cases. Finally, we provide insight into the computational cost of FokasNets and show that this architecture can be used to dramatically improve training time compared to standard PINNs, through second-order training. 
\end{enumerate}

The remainder of this paper is organized as follows. In  Section~\ref{sec:prelim} we provide the necessary mathematical preliminaries, including PINNs, PINOs and the general Fokas Method and Global Relation formulation. In Section~\ref{sec:training-method} we introduce FokasNets, detailing the corresponding Global Relation residual and Physics-Informed training scheme. In Section~\ref{sec:consistency-error} we develop the theoretical framework of FokasNets. Section~\ref{sec:pde-constructions} constructs the FokasNets explicitly for the elliptic and evolution PDEs we consider in this work. Finally, in Section~\ref{sec:numerics} we present numerical results using manufactured and literature benchmarks, before concluding.

\section{Problem setting and preliminaries}
\label{sec:prelim}

In this section we review the necessary mathematical preliminaries and concepts that will be used for the development of the proposed methods, including the PDE classes considered, the Physics-Informed training approach and the Global Relation formulation of the Fokas method.

Throughout, let $d\in\mathcal D$ denote the data for a given PDE, i.e., the forcing term, physical parameters, prescribed boundary and initial functions, where the admissible set $\mathcal D$ contains data for which the chosen problem and solution class are well defined. Furthermore, let $\Om\subset\R^{n_x}$ be a bounded connected Lipschitz domain with
piecewise smooth boundary. We will be considering elliptic and evolution PDEs, and their corresponding solution operator problems. The former are defined on the spatial domain $\Om$, while the latter are defined on the space-time domain $Q_{\mathcal T}=\Om\times(0,\mathcal T)$, for some $\mathcal T>0$.

We will be considering the class of second-order, semilinear elliptic PDEs, with mixed boundary conditions, given by:
\begin{equation}\label{eq:bvp}
 -\Delta u(x)+R(u(x);\lambda)=F, \,\, x\in \Om,
 \qquad \text{with } \quad u|_{\GD}=g_D,
 \qquad \dn u|_{\GN}=w_N.
\end{equation}
where $x\in\R^{n_x}$ and $\GD$ and $\GN$ are disjoint relatively open boundary portions, with
$\bd=\overline\GD\cup\overline\GN$. The parameter $\lambda \in \mathbb{R}$
enters the reaction $R$ and defines an entire family of problems, as will be considered in the context of operator learning. In the mixed problem, $g$ is unknown on $\GN$ and $w$ is unknown on $\GD$. We denote the complete value and
normal-derivative traces by $g=u|_{\bd}$ and $w=\dn u|_{\bd}$, and set:
\begin{equation}\label{eq:q}
 q[u;d]:=R(u;\lambda)-F,\qquad \Delta u=q[u;d].
\end{equation}
In this setting, the datum is the vector $d = (F,g_D,w_N,\lambda)$.

Next, we consider first-order evolution equations, where $\mathcal L_x$ is a
linear spatial differential operator of order $r$ with sufficiently smooth
real coefficients and $\mathcal N$ is a nonlinear operator. Since
$\mathcal L_x$ is of order $r$, each boundary point carries $r$ traces, which
we collect in the \emph{trace collection}:
\begin{equation}\label{eq:trace-collection}
 \boldsymbol\gamma u:=\bigl(u,\dn u,\ldots,\partial_n^{r-1}u\bigr)\big|_{\bd}.
\end{equation}
A boundary condition prescribes, on each portion of $\bd$, a subset of these
traces. We write \emph{$\pre$} for the index set of the prescribed traces and \emph{$\un$} for
the complementary, unknown ones, and denote the corresponding components by
$\boldsymbol\gamma_\pre u$ and $\boldsymbol\gamma_\un u$. This notation applies to both classes of problems (see Remark \ref{rem:trace-example} below), and will be used throughout the paper. With $(x,t)\in Q_{\mathcal T}$,
the initial-boundary-value problem is then:
\begin{equation}\label{eq:evo-bvp}
 \partial_tu+\mathcal L_xu+\mathcal N[u;\lambda]=F \text{ on } Q_{\mathcal T},
 \qquad u(x,0)=u_0(x) \text{ on } \Om,
 \qquad \boldsymbol\gamma_\pre u=\boldsymbol d_\partial \text{ on } \bd\times(0,\mathcal T),
\end{equation}
where $\boldsymbol d_\partial:\bd\times(0,\mathcal T)\to\R^{|\pre|}$ is the
prescribed boundary datum, i.e., $\boldsymbol d_\partial(x,t)$ collects the
values of the $|\pre|$ prescribed traces at the boundary point $x$ and time
$t$, in the order of $\boldsymbol\gamma_\pre u$. The number of prescribed
traces may differ between boundary portions, and the prescribed set is
chosen so that the problem is well posed. The datum
here is, therefore, $d=(F,u_0,\boldsymbol d_\partial,\lambda)$, and admissible
data satisfy the initial--boundary compatibility conditions, e.g.,
$\boldsymbol\gamma_\pre u_0=\boldsymbol d_\partial(\cdot,0)$ whenever these traces
are defined. The unknown traces $\boldsymbol\gamma_\un u$ are part of the
solution, and the Global Relation of the Fokas method is precisely the identity that couples them to the prescribed ones. 

\begin{remark}[Example of the trace collection]\label{rem:trace-example}
As an illustration, consider the mixed second-order elliptic problem \eqref{eq:bvp}. Here, the trace collection is
$\boldsymbol\gamma u=(u,\dn u)|_{\bd}$, and the mixed boundary condition of
\eqref{eq:bvp} prescribes $u$ on $\GD$ and $\dn u$ on $\GN$. Hence, we have
$\boldsymbol\gamma_\pre u=(u|_{\GD},\dn u|_{\GN})$ with
$\boldsymbol d_\partial=(g_D,w_N)$, while
$\boldsymbol\gamma_\un u=(\dn u|_{\GD},u|_{\GN})$ collects the unknown Neumann
and Dirichlet traces. The same holds for the evolution problem, with a second-order spatial differential operator ($r=2$), with $\boldsymbol{d}_{\partial}$ now depending on both space and time. 
\end{remark}

\subsection{Physics-Informed Neural Networks and Neural Operators}
\label{sec:prelim-neural-approximations}

We now recall the Physics-Informed Neural Network (PINN) and Physics-Informed Neural Operator (PINO) formulations for the PDEs described above \cite{raissi2019physics,lu2021learning,wang2021learning}. 
Let $\mathcal V$ be a space of sufficiently smooth real-valued functions on
$\overline\Om$ (elliptic case) or on $\overline Q_{\mathcal T}$ (evolution
case), on which the differential operator and the prescribed conditions of
\eqref{eq:bvp} or \eqref{eq:evo-bvp} are defined. Both problems can be
written in the common form:
\begin{equation}\label{eq:prelim-common-form}
 \mathcal R_{\mathrm{PDE}}[u;d]=0\ \text{in the interior},\quad
 \boldsymbol\gamma_\pre u=\boldsymbol d_\partial\ \text{on the boundary},\quad
 u(\cdot,0)=u_0\ \text{(evolution only)},
\end{equation}
where the differential residual is defined as:
\begin{equation}\label{eq:prelim-pde-residual}
 \mathcal R_{\mathrm{PDE}}[u;d]:=
 \begin{cases}
 -\Delta u+R(u;\lambda)-F,&\text{elliptic},\\
 \partial_tu+\mathcal L_xu+\mathcal N[u;\lambda]-F,&\text{evolution}.
 \end{cases}
\end{equation}
Furthermore, for a finite sample set $\mathscr S=\{X_a\}$, where $X_a = x_a$ and $X_a = (x_a,t_a)$, for the elliptic and evolution problems respectively, and positive weights $\omega_a$, define the weighted mean square over the sample by:
\begin{equation}\label{eq:fpinn-discrete-norm}
 \|f\|_{\mathscr S}^2:=\frac{\sum_a\omega_a|f(X_a)|^2}{\sum_a\omega_a}.
\end{equation}

Consider a neural network $u_\theta:\overline\Om\to\R$ (or
$\overline Q_{\mathcal T}\to\R$) with trainable weights
$\theta\in\Theta\subset\R^{N_\theta}$, smooth enough that
$u_\theta\in\mathcal V$. Given interior, boundary and initial sample sets $\mathscr S_{\mathrm{PDE}}$,
$\mathscr S_{\rm BC}$, $\mathscr S_0$ and positive weights $\lambda_{\mathrm{PDE}},\lambda_{\rm BC},\lambda_0 >0$, the Physics-Informed loss function is given by:
\begin{equation}\label{eq:prelim-physics-objective}
 \mathcal J(\theta; d)
   :=\lambda_{\mathrm{PDE}}\|\mathcal R_{\mathrm{PDE}}[u_{\theta};d]\|_{\mathscr S_{\mathrm{PDE}}}^2
     +\lambda_{\rm BC}\|\boldsymbol\gamma_\pre u_{\theta}-\boldsymbol d_\partial\|_{\mathscr S_{\rm BC}}^2
     +\lambda_0\|u_{\theta}(\cdot,0)-u_0\|_{\mathscr S_0}^2,
\end{equation}
and a Physics-Informed Neural Network (PINN) is a neural network trained by minimizing \eqref{eq:prelim-physics-objective} \cite{raissi2019physics}. (In the above, $\lambda_0 = 0$ for elliptic problems.) Note that for the mixed second-order problems the boundary term splits into a Dirichlet and a Neumann part, $\mathcal J_D$ and $\mathcal J_N$, calculated over the corresponding sample sets $\mathscr S_D\subset\GD$, $\mathscr S_N\subset\GN$ with weights $\lambda_D,\lambda_N$. The result of training is the solution approximation $u_{\theta^*}$, where $\theta^*$ denotes the trained weights. Notice that, since the boundary and initial conditions
enter \eqref{eq:prelim-physics-objective} only through penalty terms, they are satisfied only approximately, to an accuracy determined by the optimization. This is precisely the cause of the high boundary errors that can arise in PINN training, and the motivation for the FokasNet framework developed in this work.

For operator learning, let $\mathcal S:\mathcal D\to\mathcal V$,
$d\mapsto u(\cdot;d)$, be the solution map of the problem over the
admissible data. A \emph{neural operator} is a neural network
$u_\theta:\overline\Om\times\mathcal D\to\R$ (or
$\overline Q_{\mathcal T}\times\mathcal D\to\R$), with weights shared
across the whole family, which defines the approximation:
\begin{equation}\label{eq:prelim-neural-operator}
 \mathcal S_\theta(d):=u_\theta(\cdot;d)\approx\mathcal S(d).
\end{equation}
A Physics-Informed Operator (PINO) is a neural operator trained by averaging the Physics-Informed loss over training data $d_1,\ldots,d_{N_d}$, yielding the loss function: 
\begin{equation}\label{eq:prelim-operator-objective}
 \mathcal J_{\mathrm{op}}(\theta):=\frac1{N_d}\sum_{i=1}^{N_d}
   \mathcal J(\theta;d_i),
\end{equation}
so that no solution data are needed. In this work we will be focusing on  the DeepONet
branch-trunk architecture~\cite{lu2021learning} for $u_\theta(\cdot;d)$, trained through \eqref{eq:prelim-operator-objective}~\cite{wang2021learning}. However, we note that FokasNets can be used with any neural operator architecture, including Fourier Neural Operators (FNOs) \cite{li2020fourier} and their Physics-Informed variants \cite{goswami2023physics}.

\subsection{The Fokas method and the Global Relation}
\label{sec:prelim-global-relations}

The Fokas method, or unified transform, was introduced in \cite{fokas1997unified} for linear and integrable nonlinear PDEs
and is presented in \cite{fokas2012synthesis}. Further explorations have been given in \cite{fokas2001two} for evolution equations on an interval and \cite{ashton2008fundamental} for its formulation in terms of the fundamental differential form. The method starts from the \emph{Global Relation} (GR), an identity obtained by testing the PDE against solutions of its formal adjoint, which couples all boundary traces of the solution, prescribed and unknown, to the data. As mentioned, in the full method, the GR is combined with a spectral integral representation of the solution. Importantly, the numerical treatment of the GR for elliptic problems has been the topic of important work, such as
\cite{fornberg2011numerical,colbrook2019hybrid,hashemzadeh2015numerical} and further analysed in
\cite{ashton2016numerical}. FokasNets use the key GR component, in order to determine the unknown traces $\boldsymbol\gamma_\un u$
from the prescribed ones. We now present the necessary notation and concepts pertaining to the GR, as will be used throughout.

Let $\mathcal L_x^*$ denote the formal adjoint of $\mathcal L_x$ and
define the boundary concomitant $\mathscr B_{\mathcal L}$ through the identity:
\begin{equation}\label{eq:prelim-evo-lagrange}
 \int_{\bd}\mathscr B_{\mathcal L}(u,\varphi)\,\mathrm ds
  :=\int_\Om\bigl(u\mathcal L_x^*\varphi-\varphi\mathcal L_xu\bigr)\,\mathrm dx,
\end{equation}
which, by integration by parts, is a bilinear expression in the traces
$\boldsymbol\gamma u$ and $\boldsymbol\gamma\varphi$. For example, $\mathcal L_x=-\Delta$ is formally self-adjoint,
$\mathcal L_x^*=-\Delta$, and Green's second identity
$\int_\Om(\varphi\Delta u-u\Delta\varphi)\,\mathrm dx
 =\int_{\bd}(\varphi\dn u-u\dn\varphi)\,\mathrm ds$
gives the concomitant
$\mathscr B_{\mathcal L}(u,\varphi)=\varphi\dn u-u\dn\varphi$.
The GR is obtained by choosing the test functions $\varphi$ in the spaces
of solutions of the adjoint equation:
\begin{align}
 \mathcal H_E&:=\{\varphi\in C^2(\overline\Om;\mathbb C):\Delta\varphi=0\},
   \label{eq:prelim-harmonic-tests}\\
 \mathcal H_T&:=\{\varphi\in C^{1,r}(\overline Q_{\mathcal T};\mathbb C):
   -\partial_t\varphi+\mathcal L_x^*\varphi=0\},
   \label{eq:prelim-evo-adjoint}
\end{align}
for the elliptic and evolution problems, respectively. In the Fokas
method these are the one-parameter families of exponential solutions,
e.g., $\varphi(x,t)=e^{-\mathrm ikx-\omega(k)t}$ with $k\in\mathbb C$ and
$\omega(k)$ determined by $\mathcal L_x$, which turn the GR into a
relation between Fourier-type transforms of the traces. However, any other
adjoint solutions (test functions), e.g., harmonic polynomials in the elliptic case, are
equally admissible, and we use both in the sequel. Note that, in the case of complex-valued tests, their real and imaginary parts provide two real
identities via the GR. The following proposition provides the GR for both problem classes.

\begin{proposition}[Global Relations
\cite{fokas1997unified,fokas2012synthesis}]
\label{prop:prelim-global-relations}
Let $u$ be a sufficiently smooth solution of a problem belonging to the class of elliptic PDEs
\eqref{eq:bvp}. Then, for every $\varphi\in\mathcal H_E$, the Global Relation is given by:
\begin{equation}\label{eq:GR}
 \int_{\bd}(\varphi\dn u-u\dn\varphi)\,\mathrm ds=\int_\Om\varphi q[u;d]\,\mathrm dx.
\end{equation}
Analogously, let $u$ be a sufficiently smooth solution of a problem belonging to the class of evolution PDEs
\eqref{eq:evo-bvp}. Then, for every $\varphi\in\mathcal H_T$ and every
$0<\tau\le\mathcal T$, the Global Relation is given by:
\begin{equation}\label{eq:prelim-evo-GR} 
\begin{aligned}
 \int_0^\tau\!\int_{\bd}\mathscr B_{\mathcal L}(u,\varphi)\,\mathrm ds\,\mathrm dt
 ={}&\int_0^\tau\!\int_\Om\varphi q[u;d]\,\mathrm dx\,\mathrm dt\\
 &+\int_\Om\varphi(x,\tau)u(x,\tau)\,\mathrm dx
  -\int_\Om\varphi(x,0)u_0(x)\,\mathrm dx,
\end{aligned}
\end{equation}
where $q[u;d]:=\mathcal N[u;\lambda]-F$.
\end{proposition}

In the above, \eqref{eq:GR} is Green's second identity with a harmonic test function, whereas \eqref{eq:prelim-evo-GR} follows by integrating $\partial_tu+\mathcal L_xu+q[u;d]=0$ against $\varphi$ over $\Om\times(0,\tau)$ and using \eqref{eq:prelim-evo-lagrange} together with the adjoint equation (see \cite{fokas2012synthesis}
for details). In both identities the left-hand side contains the completed trace collection $\boldsymbol\gamma u$, i.e., both $\boldsymbol\gamma_\pre u$
and $\boldsymbol\gamma_\un u$, while the right-hand side contains the source. However, note that the evolution relation additionally contains the initial moment,
which is known (prescribed) data, as well as the terminal moment $\int_\Om\varphi(x,\tau)u(x,\tau)\,\mathrm dx$, which is not. Hence, this term depends on the underlying approximation of the solution, just as the term depending on the source (which appears in both identities).

\begin{remark}[Example of the Global Relation]\label{rem:gr-example}
For the mixed elliptic problem \eqref{eq:bvp}, recall from
Remark~\ref{rem:trace-example} that
$\boldsymbol\gamma_\pre u=(u|_{\GD},\dn u|_{\GN})=(g_D,w_N)$ and
$\boldsymbol\gamma_\un u=(\dn u|_{\GD},u|_{\GN})$. Splitting the boundary
integral in \eqref{eq:GR} over $\GD$ and $\GN$ and moving the prescribed
traces to the right-hand side gives, for every $\varphi\in\mathcal H_E$:
\begin{equation}\label{eq:GR-example}
 \underbrace{\int_{\GD}\varphi\,\dn u\,\mathrm ds
   -\int_{\GN}u\,\dn\varphi\,\mathrm ds}_{\text{unknown traces }\boldsymbol\gamma_\un u}
 =\underbrace{\int_{\GD}g_D\,\dn\varphi\,\mathrm ds
   -\int_{\GN}\varphi\,w_N\,\mathrm ds}_{\text{prescribed traces }\boldsymbol\gamma_\pre u}
 +\underbrace{\int_\Om\varphi\,q[u;d]\,\mathrm dx}_{\text{source}}.
\end{equation}
For a fixed test $\varphi$, the left-hand side is a linear functional of
the unknown traces, with coefficients given by the traces of $\varphi$,
while the right-hand side contains the prescribed data and the source. Each test function therefore supplies one equation for
$\boldsymbol\gamma_\un u$, whose left-hand side is linear in the unknown
traces. Notice how, for $R\neq0$ the
source term depends on the interior values of $u$, which couples the determination of the traces to the solution itself. The latter is the case in which the Fokas method is difficult to apply directly, and the
case that FokasNets address by using the neural network approximation $u_\theta$.
\end{remark}

\section{Fokas Neural Networks}
\label{sec:training-method}

With the above frameworks at hand, we can now develop the Fokas Neural
Networks (FokasNets). The idea can be succinctly described as follows: instead of imposing the prescribed conditions through penalty terms, we create a neural network approximation whose traces are, by construction, a certain category of boundary values identically (e.g., the Dirichlet values). Hence, some of these traces are prescribed through $\boldsymbol d_\partial$, whereas others are a subset of the unknown traces $\boldsymbol\gamma_\un u$, which must be approximated. Then, during training, these unknown traces are determined by solving the GR with the current neural network approximation, and are then fed back into the neural network approximation satisfying the necessary boundary conditions. Hence, at each training iteration, the unknown traces are updated based on the current neural network approximation, which is then updated through the Physics-Informed training step. We first provide the necessary ingredients to define FokasNets at the continuous level, and subsequently provide details regarding the discretization (Section~\ref{sec:fpinn-finite-gr}) and the training procedures  (Section~\ref{sec:training-objective}), which are required for the implementation of FokasNets.

\subsection{Continuous formulation}\label{sec:fpinn-continuous}

To achieve the desired formulation, we first consider the lifting that identically satisfies the prescribed boundary conditions. As mentioned, implementing the lift that enforces the Dirichlet data values is the most standard in such implementations (e.g., \cite{lagaris1998artificial}), due to its relative simplicity. This is the approach we adopt for FokasNets, but note that other lifting functions can also be used that impose flux values and more \cite{dong2026novel}. With the present approach, the unknown traces $\boldsymbol\gamma_\un u$ can be split into those that the lift encodes (Dirichlet), and the remaining ones, which the lift does not depend on and which enter only through the Global Relation. We make this rigorous in the definition below, which will allow us to construct the lifted neural network approximation.

\begin{definition}[Lift and gate]\label{def:lift}
Consider a problem of the form \eqref{eq:bvp} or \eqref{eq:evo-bvp}, with trace collection $\boldsymbol\gamma u$ and prescribed traces $\boldsymbol\gamma_\pre u=\boldsymbol d_\partial$, and write $\Sigma$ for the boundary given by $\Sigma=\bd$ for the elliptic problem and $\Sigma=\bd\times(0,\mathcal T)$ for the evolution problem. Consider a partition of the components of the trace collection $\boldsymbol\gamma u$ into those that will be \emph{encoded} by the lift and the remaining ones, indicated by the subscripts ${\rm e}$ and ${\rm a}$. The encoded components comprise the prescribed traces that the lift reproduces and the unknown traces that it carries; for a value-only lift, which is the standard choice adopted in this work, these are the Dirichlet values on $\GD$ and on $\GN$, respectively, while the Neumann data on $\GN$ are not encoded. Let $|\un_{\rm e}|$ and $|\un_{\rm a}|$ denote the number of encoded and non-encoded unknown components at a given boundary point, so that $|\un_{\rm e}|+|\un_{\rm a}|=|\un|$. A candidate for the unknown traces is a pair $(\boldsymbol\zeta,\boldsymbol\zeta')$ of functions:
\begin{equation}\label{eq:fpinn-candidate-traces}
 \boldsymbol\zeta:\Sigma\to\R^{|\un_{\rm e}|},\qquad
 \boldsymbol\zeta':\Sigma\to\R^{|\un_{\rm a}|},
\end{equation}
with the regularity required by the problem. Then:
\begin{enumerate}
\item A \emph{lift} is a map $L(\cdot;d):\mathcal Z\to\mathcal V$, affine
in $\boldsymbol\zeta$, which reproduces the encoded prescribed traces and
the encoded candidate traces, i.e., satisfies the identities:
\begin{equation}\label{eq:fpinn-lift-traces}
 \bigl(\boldsymbol\gamma_\pre L(\boldsymbol\zeta;d)\bigr)_{\rm e}=(\boldsymbol d_\partial)_{\rm e},\qquad
 \bigl(\boldsymbol\gamma_\un L(\boldsymbol\zeta;d)\bigr)_{\rm e}=\boldsymbol\zeta,
 \qquad\text{for every }\boldsymbol\zeta\in\mathcal Z.
\end{equation}
For evolution problems, additionally $L(\boldsymbol\zeta;d)(\cdot,0)=u_0$.
Prescribed traces that are not encoded, such as the Neumann data for a
value-only lift, are not reproduced by the lift and are imposed through
the corresponding term of the loss.
\item A \emph{gate} is a fixed function $\beta:\overline\Om\to\R$ (or
$\overline Q_{\mathcal T}\to\R$) such that
$\bigl(\boldsymbol\gamma(\beta v)\bigr)_{\rm e}=0$ for every smooth
$v$, i.e., the correction does not alter any encoded trace, and
$\beta(\cdot,0)=0$ in the evolution case. Specifically for a lift that encodes
Dirichlet values only, $\beta|_{\bd}=0$.
\end{enumerate}
\end{definition}

As an illustrative example, in the mixed elliptic case, $\boldsymbol\zeta:\GN\to\R$
is the candidate for $u|_{\GN}$ and $\boldsymbol\zeta':\GD\to\R$ for
$\dn u|_{\GD}$ and the full Dirichlet trace required by the lift is $(g_D,\boldsymbol\zeta)$ (for the evolution problem, the same are defined on $\GN\times(0,\mathcal T)$
and $\GD\times(0,\mathcal T)$, respectively). 

Finally, given a neural network
$v_\theta$, the lifted neural network approximation is defined as:
\begin{equation}\label{eq:fpinn-trial-field}
 u_\theta(X;\boldsymbol\zeta,d):=L(\boldsymbol\zeta;d)(X)+\beta(X)\,v_\theta(X;d),
\end{equation}
with $X = x$ and $X = (x,t)$ for the elliptic and evolution problems, respectively. We will refer to this simply as the neural approximation or solution approximation, for brevity, when there is no risk of confusion. Therefore, by construction, the encoded prescribed traces of $u_\theta(\cdot;\boldsymbol\zeta,d)$ equal $(\boldsymbol d_\partial)_{\rm e}$ and its encoded unknown traces equal $\boldsymbol\zeta$, for every $\theta$ and $\boldsymbol\zeta$; the non-encoded prescribed traces (the Neumann data, for a value-only lift) are left to the loss. The neural correction is activated in the interior, where the lift extends the boundary conditions smoothly in order to facilitate Physics-Informed training \cite{sukumar2022exact,dong2026novel}.

We now turn to the Global Relation, which will be used to determine the unknown traces $\boldsymbol\zeta$ and $\boldsymbol\zeta'$ during training. 
In analogy with the differential residual $\mathcal R_{\mathrm{PDE}}$, we
define the same quantity for the Global Relation below, which is a key component and contribution of our approach.

\begin{definition}[Global Relation residual and completion]\label{def:gr-residual}
The Global Relation residual of $u$ is the functional $r[u;d]$ defined by:
\begin{equation}\label{eq:fpinn-gr-residual}
\begin{aligned}
 r[u;d](\varphi)&:=\int_\Om\varphi\,q[u;d]\,\mathrm dx
   -\int_{\bd}\mathscr B_{\mathcal L}(u,\varphi)\,\mathrm ds, \,\,\,\, \varphi\in\mathcal H_E,\\
 r[u;d](\varphi)&:=\int_0^\tau\!\!\int_\Om\varphi\,q[u;d]\,\mathrm dx\,\mathrm dt
  +\int_\Om\varphi(\cdot,\tau)u(\cdot,\tau)\,\mathrm dx
  -\int_\Om\varphi(\cdot,0)u_0\,\mathrm dx\\
 &\qquad-\int_0^\tau\!\!\int_{\bd}\mathscr B_{\mathcal L}(u,\varphi)\,\mathrm ds\,\mathrm dt,
 \,\,\,\,\, \varphi\in\mathcal H_T,\ 0<\tau\le\mathcal T,
\end{aligned} 
\end{equation}
for the elliptic and the evolution problem, respectively, and where $u\in\mathcal V$. Furthermore, for the lifted approximation $u_\theta(\cdot;\boldsymbol\zeta,d)$ of
\eqref{eq:fpinn-trial-field}, whose encoded trace coordinates are the prescribed data $\boldsymbol d_\partial$ and $\boldsymbol\zeta$, split the
concomitant by bilinearity into the contributions of the prescribed and of
the unknown traces by
$\mathscr B_{\mathcal L}(u_\theta,\varphi)
 =\mathscr B_{\mathcal L}^{\pre}(\boldsymbol d_\partial,\varphi)
 +\mathscr B_{\mathcal L}^{\un}(\boldsymbol\zeta,\boldsymbol\zeta';\varphi)$,
and define:
\begin{equation}\label{eq:fpinn-gr-residual-parts}
\begin{aligned}
 a(\boldsymbol\zeta,\boldsymbol\zeta';\varphi)
   &:=-\int_{\bd}\mathscr B_{\mathcal L}^{\un}(\boldsymbol\zeta,\boldsymbol\zeta';\varphi)\,\mathrm ds,
 &\ell(\varphi;d)
   &:=\int_{\bd}\mathscr B_{\mathcal L}^{\pre}(\boldsymbol d_\partial,\varphi)\,\mathrm ds,
 &\mathcal M[u;d](\varphi)
   &:=\int_\Om\varphi\,q[u;d]\,\mathrm dx,
\end{aligned}
\end{equation}
for the elliptic problem and, for the evolution problem:
\begin{equation}\label{eq:fpinn-gr-residual-parts-evo}
\begin{aligned}
 a(\boldsymbol\zeta,\boldsymbol\zeta';\varphi)
   &:=-\int_0^\tau\!\!\int_{\bd}\mathscr B_{\mathcal L}^{\un}(\boldsymbol\zeta,\boldsymbol\zeta';\varphi)\,\mathrm ds\,\mathrm dt,\\
 \ell(\varphi;d)
   &:=\int_0^\tau\!\!\int_{\bd}\mathscr B_{\mathcal L}^{\pre}(\boldsymbol d_\partial,\varphi)\,\mathrm ds\,\mathrm dt
      +\int_\Om\varphi(\cdot,0)u_0\,\mathrm dx,\\
 \mathcal M[u;d](\varphi)
   &:=\int_0^\tau\!\!\int_\Om\varphi\,q[u;d]\,\mathrm dx\,\mathrm dt
      +\int_\Om\varphi(\cdot,\tau)u(\cdot,\tau)\,\mathrm dx.
\end{aligned}
\end{equation}
Then the Global Relation residual of the (lifted) neural approximation can be written as: 
\begin{equation}\label{eq:fpinn-gr-residual-structure}
 r[u_\theta(\cdot;\boldsymbol\zeta, \boldsymbol\zeta', d);d](\varphi)
 =a(\boldsymbol\zeta,\boldsymbol\zeta';\varphi)-\ell(\varphi;d)
  +\mathcal M[u_\theta(\cdot;\boldsymbol\zeta,d);d](\varphi).
\end{equation}
Finally, fix a family $\Phi$ of test functions in $\mathcal H_E$, or in
$\mathcal H_T$ together with upper limits $\tau$. For given weights
$\theta \in \Theta$, a \emph{completion} of the lifted approximation is a pair of
unknown traces
$(\boldsymbol\zeta^*,\boldsymbol\zeta'^*)=(\boldsymbol\zeta^*(\theta,d),\boldsymbol\zeta'^*(\theta,d))\in\mathcal Z\times\mathcal Z'$
for which the Global Relation residual
\eqref{eq:fpinn-gr-residual-structure} vanishes on $\Phi$:
\begin{equation}\label{eq:fpinn-completion-continuous}
 r[u_\theta(\cdot;\boldsymbol\zeta^*, \boldsymbol\zeta'^*, d);d](\varphi)
 =a(\boldsymbol\zeta^*,\boldsymbol\zeta'^*;\varphi)-\ell(\varphi;d)
 +\mathcal M[u_\theta(\cdot;\boldsymbol\zeta^*,d);d](\varphi)=0,
 \,\,\,\,\,\text{for all }\varphi\in\Phi.
\end{equation}
\end{definition}

Consequently, \eqref{eq:fpinn-completion-continuous} is the classical Global
Relation of the Fokas method, but which now depends on the neural approximation in the interior and therefore on $\theta$ and, through the lift, on $\boldsymbol\zeta$. This creates a system of equations for the unknown traces. In other words, in FokasNets, the unknown traces are not trained but computed: they are the traces for which the current lifted
approximation is compatible with the Global Relation, and they change as the neural network weights are updated. Inserting the completion
into the lifted approximation yields the FokasNet solution approximation. 

We can now exactly define the proposed neural approximation models, FokasNets. As mentioned, we refer to FokasNets as the general framework, which is divided into the Fokas-PINN, for the single PDE problem and FONet, for the solution operator, both defined below. 

\begin{definition}[Fokas-Physics-Informed Neural Network (Fokas-PINN)]\label{def:fokas-pinn}
Let $u_\theta(\cdot;\boldsymbol\zeta,d)$ be the solution approximation
\eqref{eq:fpinn-trial-field} and $(\boldsymbol\zeta^*(\theta,d), \boldsymbol\zeta'^*(\theta,d)) $ a
completion approximation computed to minimize the Global Relation residual of Definition \ref{def:gr-residual}. Then, the Fokas-PINN is the neural network approximation given by:
\begin{equation}\label{eq:fpinn-completed-field}
 u_\theta(X;d):=u_\theta(X;\boldsymbol\zeta^*(\theta,d),d)
   =L(\boldsymbol\zeta^*(\theta,d);d)(X)+\beta(X)\,v_\theta(X;d),
\end{equation}
with neural weights $\theta \in \Theta$ trained by minimizing the Physics-Informed loss function
\eqref{eq:prelim-physics-objective}:
\begin{equation}\label{eq:loss-general}
 \mathcal J(\theta;d)= \mathcal J(\theta;\boldsymbol{\zeta}^*,d).
\end{equation}
\end{definition}

Extending the above, we formulate the definition for the operator learning setting.

\begin{definition}[Fokas Operator Network (FONet)]\label{def:fokas-pinn-operator}
Let $v_\theta(\cdot;d)$ be a neural operator correction as in
\eqref{eq:prelim-neural-operator}, with weights shared across the data
family, and let the lift, gate and test family be fixed across the
family. Consider the datumwise completion approximation $\boldsymbol\zeta^*(\theta,d)$, which minimizes the Global Relation residual of Definition \ref{def:gr-residual} for each $d\in\mathcal D$. Then, the
Fokas Operator Network (FONet) is the neural operator approximation given by:
\begin{equation}\label{eq:opfield}
 \mathcal S_\theta(d):=u_\theta(\cdot;d)
  =L(\boldsymbol\zeta^*(\theta,d);d)+\beta\,v_\theta(\cdot;d),
\end{equation}
with shared neural network weights $\theta \in \Theta$ for the entire family of PDEs, trained on data $d_1,\ldots,d_{N_d}$ by minimizing the Physics-Informed Operator loss function \eqref{eq:prelim-operator-objective}:
\begin{equation}\label{eq:operator-loss}
 \mathcal J_{\mathrm{op}}(\theta):=\frac1{N_d}\sum_{i=1}^{N_d}\mathcal J(\theta; \boldsymbol{\zeta}^*, d_i).
\end{equation}
\end{definition}

We now describe how we numerically approximate and construct the necessary quantities for FokasNets, particularly the Global Relation residual. The following implementation details are given in their general form, with the precise, PDE class-specific construction outlined in Section \ref{sec:pde-constructions}.

\subsection{Numerical implementation of the Global Relation residual}
\label{sec:fpinn-finite-gr}\label{sec:fpinn-lift}

To move from the continuous formulation to the finite version which will be implemented, we describe the necessary steps to turn the Global Relation residual \eqref{eq:fpinn-gr-residual} and the corresponding completion condition
\eqref{eq:fpinn-completion-continuous}, into a finite system of equations. 

We begin by considering the trace representation and lift functions. For the trace, expand the candidate unknown traces in finite bases, by:
\begin{equation}\label{eq:fpinn-full-trace-representation}
 \boldsymbol\zeta=\sum_{j=1}^m z_j\boldsymbol\zeta_j,\qquad
 \boldsymbol\zeta'=\sum_{l=1}^{m'}z'_l\boldsymbol\zeta'_l,
\end{equation}
so that $z\in\R^m$ and $z'\in\R^{m'}$ are the coefficient vectors of $\boldsymbol\zeta$ and $\boldsymbol\zeta'$, respectively. In the evolution construction, the basis
functions also depend on time, i.e.,  $\boldsymbol\zeta_j=\boldsymbol\zeta_j(x,t)$. The complete trace collection of the approximation is then
$\boldsymbol t(z,z';d)=\boldsymbol t^{\mathrm{dat}}(d)+\sum_j\boldsymbol\zeta_jz_j+\sum_l\boldsymbol\zeta'_lz'_l$,
where $\boldsymbol t^{\mathrm{dat}}(d)$ carries the prescribed components.
 We write $u_\theta(\cdot;z,d)$ for the approximation
\eqref{eq:fpinn-trial-field} with $\boldsymbol\zeta$ given by
\eqref{eq:fpinn-full-trace-representation}. 

We now consider the test functions for the GR and the integral quadrature. 
Choose finitely many tests $\varphi_k$ from $\mathcal H_E$ or
$\mathcal H_T$. After splitting complex identities into real and
imaginary parts, we obtain a system of $p$ real equations, indexed by
$k=1,\ldots,p$; in
evolution problems the index also identifies the upper time limit
$\tau_k$. Integrals over a set $S$ are approximated by quadrature, which we denote simply by:
\begin{equation}\label{eq:fpinn-quadrature-functional}
 \qint_S[f]:=\sum_{a=1}^{N_S}\omega_{S,a}f(X_{S,a})\approx\int_Sf\,\mathrm d\mu_S,
\end{equation}
with the nodes $\{X_{S,a}\}_{a=1}^{N_S}$, the weights $\{\omega_{S,a}\}_{a=1}^{N_S}$. Evaluating the three
groups comprising the Global Relation residual \eqref{eq:fpinn-gr-residual-structure} by quadrature for the tests $\varphi_k$ defines the boundary functional $\mathfrak b[\boldsymbol t]=(\mathfrak b_k[\boldsymbol t])_{k=1}^p\in\R^p$, with $\mathfrak b_k[\boldsymbol t]\in\R$ linear in the trace collection $\boldsymbol t$, the moments $\mathbf m(z;\theta,d)\in\R^p$ of the lifted approximation $u_\theta(\cdot;z,d)$ and the known initial contributions $\mathbf c^0(d)\in\R^p$, with entries given by:
\begin{equation}\label{eq:fpinn-signed-boundary-functional}
 \mathfrak b_k[\boldsymbol t]=
 \begin{cases}
 -\qint_{\bd}[\mathscr B_{\mathcal L}[\boldsymbol t,\varphi_k]],&\text{elliptic},\\
 -\qint_{\bd\times(0,\tau_k)}[\mathscr B_{\mathcal L}[\boldsymbol t,\varphi_k]],&\text{evolution},
 \end{cases}
\end{equation}
\begin{equation}\label{eq:fpinn-unreduced-moments}
 \mathbf m_k(z;\theta,d)=
 \begin{cases}
 \qint_\Om[\varphi_k(R(u_\theta(\cdot;z,d);\lambda)-F)],&\text{elliptic},\\
 \qint_{\Om\times(0,\tau_k)}[\varphi_k(\mathcal N[u_\theta(\cdot;z,d);\lambda]-F)]
 +\qint_\Om[\varphi_k(\cdot,\tau_k)u_\theta(\cdot,\tau_k;z,d)],&\text{evolution},
 \end{cases}
\end{equation}
\begin{equation}\label{eq:fpinn-known-initial-moments}
 (\mathbf c^0(d))_k=
 \begin{cases}0,&\text{elliptic},\\
 \qint_\Om[\varphi_k(\cdot,0)u_0],&\text{evolution},
 \end{cases}
\end{equation}
for $k=1,\ldots,p$; the moments depend on $z$ through the lift and on $\theta$ through the neural correction. (Note that, for a general function $u\in\mathcal V$ we write $\mathbf m[u;d]$ for the same expressions with $u$ in place of $u_\theta(\cdot;z,d)$. This is consistent with the above use, for which we have $\mathbf m(z;\theta,d)=\mathbf m[u_\theta(\cdot;z,d);d]$. This functional form will be used in the analysis of Section~\ref{sec:consistency-error}.) Here, $\mathscr B_{\mathcal L}[\boldsymbol t,\varphi]$ denotes the
concomitant evaluated on the trace collection $\boldsymbol t$. (For the elliptic problem, with $\boldsymbol t=(g,w)$ the value and flux components, $-\mathscr B_{\mathcal L}[\boldsymbol t,\varphi_k]=g\,\dn\varphi_k-\varphi_kw$,
so that $\mathfrak b_k[\boldsymbol t]=\qint_{\bd}[g\,\dn\varphi_k-\varphi_kw]$.) Hence, the $k$-th GR equation $r[u;d](\varphi_k)=0$ becomes:
\begin{equation}\label{eq:fpinn-component-gr}
 \mathfrak b_k[\boldsymbol t(z,z';d)]+\mathbf m_k(z;\theta,d)=\mathbf c^0_k(d),
 \qquad k=1,\ldots,p.
\end{equation}

To gather all the GR equations into a single system, recall that  $\boldsymbol t(z,z';d)=\boldsymbol t^{\mathrm{dat}}(d)+\sum_j\boldsymbol\zeta_jz_j+\sum_l\boldsymbol\zeta'_lz'_l$. Hence, $\mathfrak b_k$ is also linear in the trace coefficients and we can write:
\begin{equation}\label{eq:fpinn-b-linear}
 \mathfrak b_k[\boldsymbol t(z,z';d)]
 =\mathfrak b_k[\boldsymbol t^{\mathrm{dat}}(d)]
  +\sum_{j=1}^m\mathfrak b_k[\boldsymbol\zeta_j]\,z_j
  +\sum_{l=1}^{m'}\mathfrak b_k[\boldsymbol\zeta'_l]\,z'_l,
 \qquad k=1,\ldots,p.
\end{equation}
Finally, by defining $\mathsf B\in\R^{p\times m}$, $\mathsf B'\in\R^{p\times m'}$ and $\mathbf f\in\R^p$, with entries:
\begin{equation}\label{eq:fpinn-assembly-entries}
 (\mathsf B)_{kj}:=\mathfrak b_k[\boldsymbol\zeta_j],\qquad
 (\mathsf B')_{kl}:=\mathfrak b_k[\boldsymbol\zeta'_l],\qquad
 (\mathbf f(d))_k:=(\mathbf c^0(d))_k-\mathfrak b_k[\boldsymbol t^{\mathrm{dat}}(d)],
\end{equation}
we gather all $p$ GR equations \eqref{eq:fpinn-component-gr} to produce the system:
\begin{equation}\label{eq:fpinn-joint-gr}
 \mathsf Bz+\mathsf B'z'+\mathbf m(z;\theta,d)=\mathbf f(d).
\end{equation}

As mentioned, the auxiliary coefficients $z'$ do not
enter the approximation as they are not incorporated in the lift. Therefore, the corresponding unknown coefficients $z'$ can be eliminated to produce a lower dimensional system only for the needed trace coefficients $z$. This is done by a simple elimination in the system \eqref{eq:fpinn-joint-gr}. Formally, choose a full-row-rank matrix $\mathsf T\in\R^{n\times p}$ with:
\begin{equation}\label{eq:fpinn-annihilator}
 \mathsf T\mathsf B'=0,\qquad \ker\mathsf T=\operatorname{range}(\mathsf B').
\end{equation}
In other words, the rows of $\mathsf T$ select the combinations of the tested equations that are free of $z'$. The kernel condition makes the reduced equations equivalent to the existence of some $z'$ satisfying \eqref{eq:fpinn-joint-gr}. With $\nu=\operatorname{rank}\mathsf B'$, a full elimination produces a system with $n=p-\nu$
equations. (If there are no auxiliary traces, $\mathsf T=I$ and $n=p$.) We can therefore define the reduced quantities:
\begin{equation}\label{eq:fpinn-reduced-matrices}
\begin{aligned}
 \mathsf K&:=\mathsf T\mathsf B\in\R^{n\times m},&\qquad
 \mathbf b(d)&:=\mathsf T\mathbf f(d)\in\R^n,\\
 \mathcal M(z;\theta,d)&:=\mathsf T\mathbf m(z;\theta,d)\in\R^n,&\qquad
 \mathcal M[u;d]&:=\mathsf T\mathbf m[u;d],
\end{aligned}
\end{equation}
so that the reduced system reads:
\begin{equation}\label{eq:reduced}
  \mathsf Kz+\mathcal M(z;\theta,d)=\mathbf b(d), 
\end{equation}
a system of $n$ equations for the $m$ unknown trace coefficients $z$. This construction gives rise to the following definition. 

\begin{definition}[Finite Global Relation residual]\label{def:fpinn-finite-gr}
The finite Global Relation residual of the neural approximation $u_{\theta}(\cdot; z,d)$ is given by:
\begin{equation}\label{eq:fpinn-moment-residual}
 r(z;\theta,d):=\mathsf Kz+\mathcal M(z;\theta,d)-\mathbf b(d)\in\R^n.
\end{equation}
\end{definition}

\begin{remark}[A pedagogical example]\label{rem:interval-example}
Consider a 1D elliptic PDE $-u''+R(u;\lambda)=F$ on $\Om=(-1,1)$
with $u(-1)=g_D$ and $\dn u(1)=u'(1)=w_N$, so that $\GD=\{-1\}$,
$\GN=\{1\}$, and $\dn=\pm \partial_x$ at $x=\pm 1$. The trace collection
consists of four numbers, and the unknown ones are
$\boldsymbol\zeta=u(1)$, needed by the lift, and
$\boldsymbol\zeta'=\dn u(-1)$, which enters only the global relation. Thus, we have
$m=m'=1$ and \eqref{eq:fpinn-full-trace-representation} reads:
\[
 \boldsymbol t(z,z';d)=
 \begin{pmatrix}u(-1)\\ u(1)\\ \dn u(-1)\\ \dn u(1)\end{pmatrix}
 =\underbrace{\begin{pmatrix}g_D\\0\\0\\w_N\end{pmatrix}}_{\boldsymbol t^{\mathrm{dat}}(d)}
 +\underbrace{\begin{pmatrix}0\\1\\0\\0\end{pmatrix}}_{\boldsymbol\zeta_1}z
 +\underbrace{\begin{pmatrix}0\\0\\1\\0\end{pmatrix}}_{\boldsymbol\zeta'_1}z'.
\]
Since $\mathcal H_E=\operatorname{span}\{1,x\}$ on an interval, there are exactly two tests, $\varphi_1=1$ and $\varphi_2=x$
($p=2$), and \eqref{eq:fpinn-assembly-entries} gives:
\[
 \mathsf B=\begin{pmatrix}0\\1\end{pmatrix},\qquad
 \mathsf B'=\begin{pmatrix}-1\\1\end{pmatrix},\qquad
 \mathbf f(d)=\begin{pmatrix}w_N\\w_N+g_D\end{pmatrix},\qquad
 \mathbf m[u;d]=\begin{pmatrix}\qint_\Om[q[u;d]]\\ \qint_\Om[x\,q[u;d]]\end{pmatrix}.
\]
Adding the two rows eliminates $z'$, i.e., $\mathsf T=(1\;\;1)$, and
\eqref{eq:fpinn-reduced-matrices} reduces to the scalar relation
$z-(2w_N+g_D)+\qint_\Om[(1+x)\,q[u;d]]=0$ for the missing value
$u(1)$. For the corresponding evolution problem,
$u_t-u_{xx}+\mathcal N[u;\lambda]=F$ with the same boundary partition and
initial datum $u_0$, the trace representation has the same four-slot
structure with time-dependent entries,
$\boldsymbol t^{\mathrm{dat}}(d)=(g_D(t),\,u_0(1),\,0,\,w_N(t))^\top$ and
basis functions $\psi_j(t)$ with $\psi_j(0)=0$ in the unknown slots. The full constructions are given in Section~\ref{sec:pde-constructions}.
\end{remark}

\subsection{Physics-Informed training with the finite Global Relation completion}
\label{sec:training-objective}\label{sec:operator}
The finite completion condition is $r(z;\theta,d)=0$, a system of $n$
equations for the $m$ unknown coefficients $z$. In general, it is a nonlinear system in
$z$, since $\mathcal M(z;\theta,d)$ depends on $z$ through the
approximation $u_\theta(\cdot;z,d)$, and it need not be square, since the number $n$ of retained tested equations is a discretization choice independent of $m$. Hence, the completion
is defined as a zero of a square defining map
$\mathcal C(\cdot;\theta,d):\R^m\to\R^m$ derived from $r$, which may be solved
iteratively for fixed $\theta$. We use three such rules, labelled (R),
(LS) and (P) throughout. With $J_r:=D_zr\in\R^{n\times m}$, the defining map is given by:
\begin{equation}\label{eq:fpinn-defining-map}
 \mathcal C(z;\theta,d):=
 \begin{cases}
 r(z;\theta,d),&\text{(R) root of the square GR system},\\[1mm]
 J_r(z;\theta,d)^\top r(z;\theta,d),&\text{(LS) nonlinear least squares},\\[1mm]
 \Pi\,r(z;\theta,d),&\text{(P) fixed compression}.
 \end{cases}
\end{equation}
Rule (R) applies when $n=m$ and solves the finite global relations
exactly, whereas rule (LS) applies to a rectangular system and is the stationarity condition of the least-squares objective:
\begin{equation}\label{eq:fpinn-ls-objective}
 \Psi(z;\theta,d):=\tfrac12\,\|r(z;\theta,d)\|_2^2,
\end{equation}
i.e., $\mathcal C=\nabla_z\Psi$, so the selected completion is an interior
local minimizer of $\Psi$. Finally, rule (P) compresses the
$n$ equations to $m$ by a fixed matrix $\Pi\in\R^{m\times n}$, independent of $(z,\theta)$, chosen so that $\Pi r=0$ is a square system whose Jacobian is close to the identity. We use it for the space--time evolution construction, where $\Pi$ is the pseudo-inverse of the linear part of the residual; the construction and the reasons for this choice are given in Section~\ref{sec:evo-completion-P}. Rules (R) and (LS) are used for the stationary and time-stepping constructions of Section~\ref{sec:pde-constructions}.

\begin{figure}[htbp]
\centering
\IfFileExists{Figures/FokasNets.pdf}{\includegraphics[width=\linewidth]{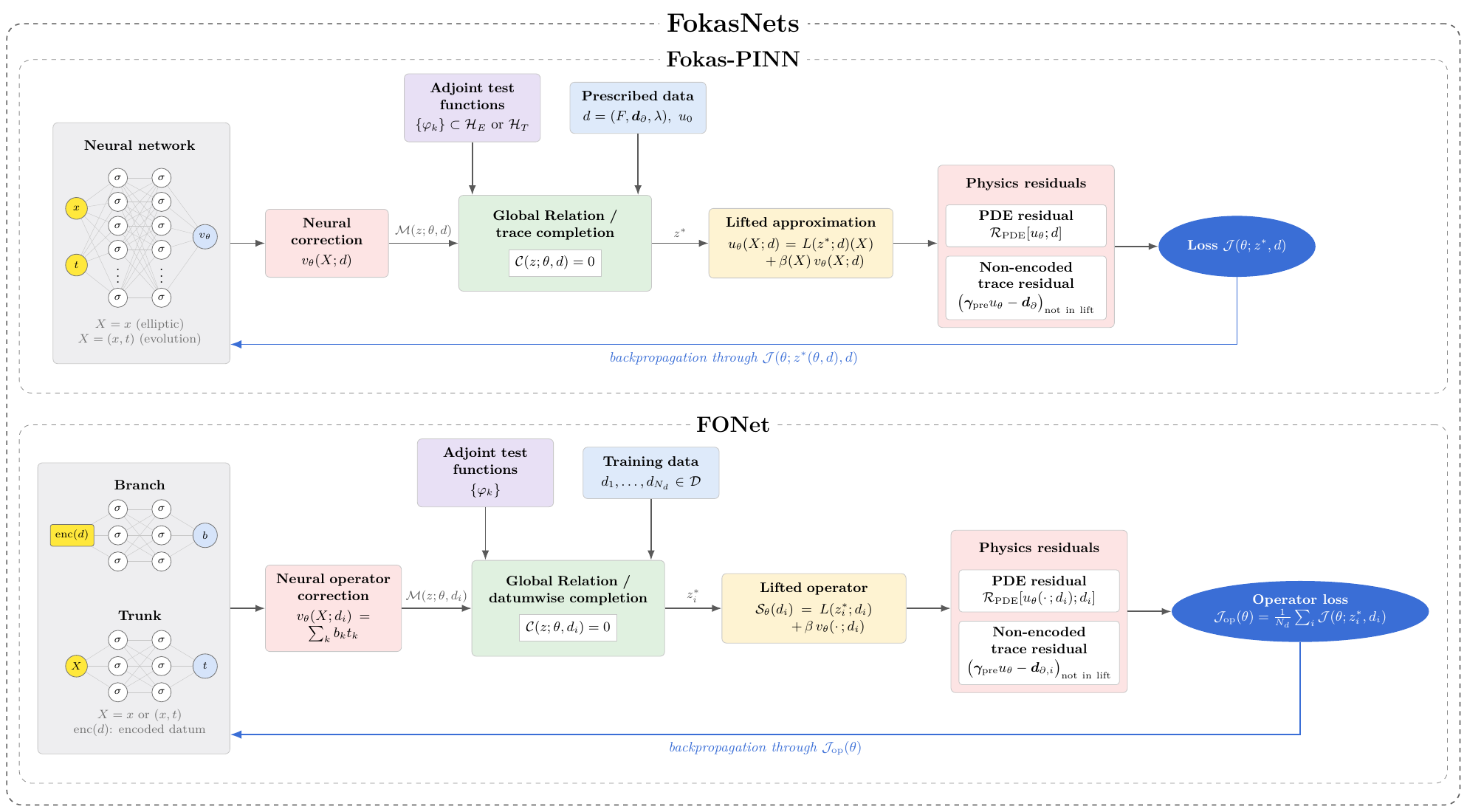}}{%
\begin{tikzpicture}[node distance=8mm and 13mm,>=Latex,
 every node/.style={font=\small},
 block/.style={draw,rounded corners,align=center,text width=43mm,minimum height=13mm}]
 \node[block] (field) {Approximation\\$L(z;d)+\beta v_\theta$};
 \node[block,right=of field] (moments) {Source and terminal\\moments};
 \node[block,below=of moments] (completion) {Trace completion\\$\mathcal C(z;\theta,d)=0$};
 \node[block,below=of field] (loss) {Completed approximation\\and physics loss};
 \draw[->] (field)--(moments);
 \draw[->] (moments)--(completion);
 \draw[->] (completion)--node[below]{$z^*$}(loss);
 \draw[->] (completion.east)--++(7mm,0)|-(field.north);
 \draw[->,dashed] (loss.west)--++(-8mm,0)|-
 node[pos=.75,above]{implicit gradient to $\theta$}(field.north west);
\end{tikzpicture}}
\caption{The FokasNets architectures. }
\label{fig:fokas-pinn}
\end{figure}

In all cases, coefficients $z^*(\theta,d)$ provide the discrete approximation $\sum_j z^*_j(\theta,d)\,\boldsymbol\zeta_j$ of the continuous completion $\boldsymbol\zeta^*(\theta,d)$ of Definition \ref{def:gr-residual}.  The Fokas-PINN and FONet are then
implemented with the completed approximation
$u_\theta(X;d)=u_\theta(X;z^*(\theta,d),d)$ of \eqref{eq:fpinn-completed-field} and \eqref{eq:opfield}.
Figure~\ref{fig:fokas-pinn} summarizes the architecture and the described training process for FokasNets.

We can now consider how the GR affects the Physics-Informed training process. In the finite formulation, the loss functions \eqref{eq:loss-general} and
\eqref{eq:operator-loss} are evaluated with the completed coefficients
$z^*(\theta,d)$ in place of $\boldsymbol\zeta^*(\theta,d)$, i.e.,
$\mathcal J(\theta;d)=\mathcal J(\theta;z^*(\theta,d),d)$, where
$\mathcal J(\theta;z,d)$ denotes the loss \eqref{eq:prelim-physics-objective}
evaluated at the lifted approximation $u_\theta(\cdot;z,d)$ for $z$
independent of the neural weights. Training minimizes $\mathcal J(\theta;d)$, or $\mathcal J_{\mathrm{op}}(\theta)$,
with respect to $\theta$, which requires
the gradient $\nabla_\theta\mathcal J(\theta;d)$. By the chain rule, this produces:
\begin{equation}\label{eq:fpinn-chain-rule}
 \nabla_\theta\mathcal J(\theta;d)
   =\Bigl[\partial_\theta\mathcal J(\theta;z,d)
    +(D_\theta z^*)^\top\nabla_z\mathcal J(\theta;z,d)\Bigr]_{z=z^*(\theta,d)},
\end{equation}
where the partial derivatives $\partial_{\theta}$ hold $z$ fixed. The first term is the usual
PINN gradient, with the completed coefficients treated as constants and the
second accounts for the fact that the completion itself depends on the updating of the neural network weights. The additional complexity lies in the sensitivity $D_\theta z^*$, since the completion is not given by an explicit formula but as a zero of the
defining map, $\mathcal C(z^*(\theta,d);\theta,d)=0$, computed by an
iterative solver. Differentiating this identity with respect to $\theta$
gives:
\begin{equation}\label{eq:fpinn-implicit-identity}
 J_C\,D_\theta z^*+\partial_\theta\mathcal C=0,
 \qquad J_C:=D_z\mathcal C\in\R^{m\times m},
\end{equation}
so that, whenever the Jacobian $J_C$ of the defining map is nonsingular,
$D_\theta z^*$ follows from a single linear solve, without
differentiating through the iterations of the inner solver. Unlike the
generally rectangular $J_r\in\R^{n\times m}$, $J_C$ is square, and its
expression depends on the completion rule:
\begin{equation}\label{eq:fpinn-ls-hessian}
 J_C=
 \begin{cases}
 J_r,&\text{(R)},\\[1mm]
 J_r^\top J_r
    +\displaystyle\sum_{i=1}^n r_iD_z^2r_i,
       &\text{(LS)},\\[1mm]
 \Pi J_r,&\text{(P)}.
 \end{cases}
\end{equation}
Finally, inserting $D_\theta z^*$ from \eqref{eq:fpinn-implicit-identity}
into \eqref{eq:fpinn-chain-rule} and regrouping yields an adjoint form of
the gradient, in which a single transpose solve with $J_C$ for an
adjoint vector $\mathbf p\in\R^m$ replaces the $m\times N_\theta$
sensitivity matrix $D_\theta z^*$. The following
proposition makes this additional training step rigorous.

\begin{proposition}[Implicit differentiation of the Global Relation completion]
\label{prop:fpinn-implicit-differentiation}
Fix $d$ and a selected completion $z_0=z^*(\theta_0,d)$.
Assume $\mathcal C$ is continuously differentiable in $(z,\theta)$
near $(z_0,\theta_0)$, $\mathcal C(z_0;\theta_0,d)=0$, and
$J_C(z_0;\theta_0,d)$ is nonsingular. Furthermore, assume the lift and gate are
independent of $\theta$, the gated neural correction is continuously
differentiable into $\mathcal V$, and $\mathcal J(\theta;z,d)$ is continuously
differentiable in $(z,\theta)$. Then, a locally unique differentiable completion branch exists, with:
\begin{equation}\label{eq:fpinn-implicit-sensitivity}
 D_\theta z^*=-(J_C^\star)^{-1}\partial_\theta\mathcal C,
\end{equation}
where $J^{\star}_C$ means evaluation at the selected completion, and:
\begin{equation}\label{eq:fpinn-field-sensitivity}
 D_\theta u_\theta
    =\beta D_\theta v_\theta+D_zL\,D_\theta z^*,
\end{equation}
where $D_zL$ is independent of $z$, since the lift is affine in $z$. Furthermore, the loss gradient is given by:
\begin{equation}\label{eq:fpinn-adjoint-gradient}
 \nabla_\theta\mathcal J(\theta;d)
    =\Bigl[\partial_\theta\mathcal J(\theta;z,d)
       -(\partial_\theta\mathcal C)^\top\mathbf p\Bigr]_{z=z^*(\theta,d)},
\end{equation}
where the adjoint vector $\mathbf p\in\R^m$ is the solution of the linear
system $(J_C^\star)^\top\mathbf p=\nabla_z\mathcal J(\theta;z,d)\big|_{z=z^*(\theta,d)}$.
For least squares, if $z_0$ is an interior local minimizer, this branch
consists of strict local minimizers for nearby $\theta$.
\end{proposition}

\begin{proof}
Since $\mathcal C$ is continuously differentiable near $(z_0,\theta_0)$,
$\mathcal C(z_0;\theta_0,d)=0$ and $D_z\mathcal C(z_0;\theta_0,d)=J_C(z_0;\theta_0,d)$
is nonsingular, the Implicit Function Theorem provides the completion branch, i.e., neighbourhoods
$U$ and $V$ and a unique continuously
differentiable map $z^*(\cdot,d):U\to V$ with $z^*(\theta_0,d)=z_0$ and
$\mathcal C(z^*(\theta,d);\theta,d)=0$ for all $\theta\in U$. Differentiating this identity with respect to
$\theta$ yields
$J_C^\star D_\theta z^*+\partial_\theta\mathcal C=0$,
which proves \eqref{eq:fpinn-implicit-sensitivity}. Differentiating
$u_\theta=L(z^*;d)+\beta v_\theta$ by the chain rule proves
\eqref{eq:fpinn-field-sensitivity}, with $D_zL$ constant since $L$ is
affine in $z$.
The chain rule gives:
\[
 \nabla_\theta\mathcal J(\theta;d)
   =\partial_\theta\mathcal J(\theta;z,d)+(D_\theta z^*)^\top\nabla_z\mathcal J(\theta;z,d).
\]
Substituting \eqref{eq:fpinn-implicit-sensitivity} into the second term, we obtain:
\[
 (D_\theta z^*)^\top\nabla_z\mathcal J
   =-\bigl((J_C^\star)^{-1}\partial_\theta\mathcal C\bigr)^\top\nabla_z\mathcal J
   =-(\partial_\theta\mathcal C)^\top\,(J_C^\star)^{-\top}\nabla_z\mathcal J
   =-(\partial_\theta\mathcal C)^\top\mathbf p,
\]
where $\mathbf p:=(J_C^\star)^{-\top}\nabla_z\mathcal J$, i.e., the solution of
$(J_C^\star)^\top\mathbf p=\nabla_z\mathcal J$. Hence, we obtain
\eqref{eq:fpinn-adjoint-gradient}, which allows us to calculate the Jacobian by obtaining
$\mathbf p\in\R^m$ from one linear solve, whereas
$D_\theta z^*\in\R^{m\times N_\theta}$ would require $N_\theta$ of them. 
\par For the last part of the statement, notice that in the (LS) formulation the defining
Jacobian is the Hessian of $\Psi$. At an interior local minimum it
is positive semidefinite, and nonsingularity makes it positive definite. Furthermore, continuity preserves positive definiteness near the selected point, so nearby stationary completions are strict local minima.
\end{proof}

In practice the gradient of
Proposition~\ref{prop:fpinn-implicit-differentiation} is computed by backpropagation, so we only
need the completion step to carry the derivative
\eqref{eq:fpinn-implicit-sensitivity}. The inner solver runs with
$\theta$ fixed and outside the differentiation graph, returning
$\widetilde z$. Writing $\bar z$ and $\bar J_C$ for $\widetilde z$ and
$D_z\mathcal C(\widetilde z;\theta,d)$, both held constant during
differentiation, the lifted approximation is evaluated at
\begin{equation}\label{eq:fpinn-ift-output}
 z_{\mathrm{out}}=\bar z-\bar J_C^{-1}\,\mathcal C(\bar z;\theta,d).
\end{equation}
At an exact completion the correction vanishes, so
$z_{\mathrm{out}}=z^*$, and since only $\mathcal C$ depends on $\theta$,
its derivative is
$-\bar J_C^{-1}\partial_\theta\mathcal C=D_\theta z^*$ by
\eqref{eq:fpinn-implicit-sensitivity}. Backpropagating the loss through
$u_\theta=L(z_{\mathrm{out}};d)+\beta v_\theta$ therefore returns
\eqref{eq:fpinn-adjoint-gradient}; the backward pass through the solve
with $\bar J_C$ is precisely the transpose solve for $\mathbf p$. For a
completion computed to finite tolerance, $z_{\mathrm{out}}$ is one
Newton step from $\widetilde z$, and its value and derivative converge
to the exact ones as $\widetilde z\to z^*$
(Appendix~\ref{app:implicit}).

Finally, we note that, for the FONet, each datum has its own completion, defining Jacobian and
adjoint vector, and the gradient of $\mathcal J_{\mathrm{op}}$ is the average of the datumwise gradients \eqref{eq:fpinn-adjoint-gradient}.
Shared trunk evaluations and their derivatives can be reused across a batch, as described in Section~\ref{sec:numerics}. The complete evaluation of the loss and its gradient is summarized in
Algorithm~\ref{alg:fokas-evaluation}.

\begin{algorithm}[!h]
\caption{Training of FokasNets (Fokas--PINN: $N_d=1$; FONet: $N_d>1$)}
\label{alg:fokas-evaluation}
\begin{algorithmic}[1]
\Require Data $d_1,\ldots,d_{N_d}$; trace bases, lift $L$ and gate $\beta$;
 tests $\varphi_k$ and quadrature; completion rule (R), (LS) or (P);
 network $v_\theta$; tolerance $\mathrm{tol}$; validation set
\Ensure Trained weights $\theta^\star$ and completed approximations
 $u_{\theta^\star}(\cdot;d)=L(z^*(\theta^\star,d);d)+\beta v_{\theta^\star}(\cdot;d)$
\Statex \textit{Offline assembly (independent of $\theta$)}
\State Assemble $\mathsf B$, $\mathsf B'$ by \eqref{eq:fpinn-assembly-entries}, the annihilator
 $\mathsf T$ of \eqref{eq:fpinn-annihilator} and $\mathsf K=\mathsf T\mathsf B$; for rule (P), the matrix $\Pi$
\State For $i=1,\ldots,N_d$: assemble $\mathbf b(d_i)=\mathsf T\mathbf f(d_i)$ and set the initial guess $z_i$
\State Tabulate the lift and gate, and their derivatives, at the collocation and quadrature nodes
\Statex \textit{Training}
\State Initialize $\theta$
\For{$s=1,\ldots,S$}
 \For{$i=1,\ldots,N_d$} \Comment{datumwise, in parallel}
  \State $z\gets z_i$ \Comment{warm start}
  \Repeat \Comment{completion with $\theta$ fixed, outside the AD graph}
   \State Evaluate $\mathcal M(z;\theta,d_i)$ from $u_\theta(\cdot;z,d_i)$ and form
    $r=\mathsf Kz-\mathbf b(d_i)+\mathcal M(z;\theta,d_i)$
   \State Form $\mathcal C(z;\theta,d_i)$ by \eqref{eq:fpinn-defining-map} and $J_C=D_z\mathcal C$ by \eqref{eq:fpinn-ls-hessian}
   \State $z\gets z-J_C^{-1}\mathcal C(z;\theta,d_i)$ \Comment{backtracking on $\|\mathcal C\|_2$}
  \Until{$\|\mathcal C(z;\theta,d_i)\|_2\le\mathrm{tol}$}
  \State \textbf{if} $J_C$ is singular, or not positive definite for (LS), \textbf{then} reject this evaluation
  \State $z_i\gets z$;\quad $\bar z\gets z$,\ $\bar J_C\gets J_C$ \Comment{held constant in AD}
  \State $z_{\mathrm{out}}\gets\bar z-\bar J_C^{-1}\mathcal C(\bar z;\theta,d_i)$ \Comment{stop-gradient output \eqref{eq:fpinn-ift-output}}
  \State $u_i\gets L(z_{\mathrm{out}};d_i)+\beta\,v_\theta(\cdot;d_i)$ and
   $\mathcal J_i\gets\mathcal J(\theta;z_{\mathrm{out}},d_i)$ \Comment{non-encoded loss terms only}
 \EndFor
 \State $\mathcal J\gets\frac1{N_d}\sum_i\mathcal J_i$ and $\nabla_\theta\mathcal J$ by backpropagation
  \Comment{equals \eqref{eq:fpinn-adjoint-gradient}}
 \State Update $\theta$ with the optimizer; each line-search trial repeats lines 6--17
 \State Every $N_{\mathrm{val}}$ iterations, evaluate the validation loss and keep the best $\theta^\star$
\EndFor
\State \Return $\theta^\star$
\end{algorithmic}
\end{algorithm}

\begin{remark}
Regarding training, we note the terms of the Physics-Informed loss function that are encoded by the lift vanish identically, and can therefore be omitted from the loss function. What remains is the PDE residual, as well as the loss terms for traces that the lift does not encode. For FokasNets specifically, the vanishing terms are the completed Dirichlet values and initial conditions. Naturally, if we formulate the lifted neural network approximation such that other boundary conditions are satisfied identically, this will also change which of the loss terms vanish.
\end{remark}

\section{Error analysis and Universal Approximation for FokasNets}
\label{sec:consistency-error}

In this section, we provide a rigorous analysis of the approximation errors and stability of the FokasNet training algorithm via the Global Relation completion. These results then allow us to formulate the Universal Approximation theorem for FokasNets. We note that, although the universal approximation capability of standard neural networks is well established, its preservation under the additional structural constraints imposed by the Global Relation completion is not automatic. For this reason, the results of this section establish that, under suitable consistency and uniform completion-stability conditions, the FokasNet architecture retains the ability to approximate the exact PDE solution arbitrarily accurately. 

In the results that follow, we will use $\|\cdot\|$ for the norm used in the estimate, on functions belonging to $\mathcal V$. Coefficient and residual vectors use Euclidean norms, and matrices use their induced
operator norms. The uniform norm is used when the moments depend continuously on pointwise values, and derivative-dependent moments require a norm controlling their derivative arguments. These will be made specific where necessary in the analysis below.

\subsection{Consistency and completion stability}
\label{sec:consistency-completion}

We first provide some necessary preliminaries. Fix a datum $d$, an exact solution $u^\dagger=u^\dagger(\cdot;d)$
and a discretization. Recall the approximation and finite global relation residual:
\[
 u_\theta(\cdot;z,d)=L(z;d)+\beta v_\theta(\cdot;d),
 \qquad
 r(z;\theta,d)=\mathsf Kz-\mathbf b(d)
                  +\mathcal M[u_\theta(\cdot;z,d);d],
\]
and $J_r=D_zr$ and $J_C=D_z\mathcal C$.
We choose represented exact trace coefficients
$z^\dagger=z^\dagger(d)$ and measure the Global Relation discrepancy by:
\begin{equation}
\label{eq:ca-insertion-defect}
 \delta
   :=\bigl\|\mathsf Kz^\dagger-\mathbf b(d) +\mathcal M[u^\dagger;d]\bigr\|_2.
\end{equation}

To convert this error into an approximation-error bound, we
compare the solution approximation using $z^\dagger$ with the approximation obtained after completion. For a fixed network, set:
\begin{equation}
\label{eq:ca-field-error}
 e:=\|u_\theta(\cdot;z^\dagger,d)-u^\dagger\|,
 \qquad
 C_E:=\|\mathcal E(d)\|_{\R^m\to(\mathcal V,\|\cdot\|)},
 \qquad \mathcal E(d):=D_zL(\cdot;d),
\end{equation}
where, since the lift is affine in $z$, $\mathcal E(d)$ is a constant
linear map and $L(z;d)-L(z';d)=\mathcal E(d)(z-z')$ for all $z,z'$.
Assume there exist constants $C_M\ge0$ and $\varepsilon_0>0$
such that:
\begin{equation}
\label{eq:ca-moment-continuity}
 \|\mathcal M[u;d]-\mathcal M[u^\dagger;d]\|_2
 \le C_M\|u-u^\dagger\|,
\end{equation}
for every $u\in\mathcal V$ satisfying
$\|u-u^\dagger\|<\varepsilon_0$. Also, for fixed $(\theta,d)$, define:
\begin{equation}\label{eq:chi}
 \chi:=
 \begin{cases}
  1, & \text{(R)},\\[1mm]
  \|J_r(z^\dagger;\theta,d)\|_2,
     & \text{(LS)},\\[1mm]
  \|\Pi\|_2,&\text{(P)}.
 \end{cases}
\end{equation}

We begin with the following Proposition that supplies both existence of a completion and an estimate of its effect on the solution approximation. 

\begin{proposition}[Global Relation completion stability and approximation error]
\label{prop:ca-completion-error}
Fix $(\theta,d)$, suppose \eqref{eq:ca-moment-continuity} holds with $e<\varepsilon_0$ and assume one
of the following conditions is satisfied:

\begin{enumerate}
\item $\mathcal C(\cdot;\theta,d)$ is continuous and strongly monotone on $\R^m$, that is,
\[
 \langle\mathcal C(z;\theta,d)-\mathcal C(w;\theta,d),z-w\rangle
       \ge\mu\|z-w\|_2^2,
 \qquad \mu>0,
\]
and set also $C_S=\mu^{-1}$.

\item $\mathcal C(\cdot;\theta,d)$ is continuously differentiable on a neighborhood of
$B:=\overline B_{\rho}(z^\dagger) = 
\{z\in\mathbb R^m:\|z-z^\dagger\|_2\le\rho\}
$ and 
there is a nonsingular matrix $J_{C,0}$ such that:
\begin{equation}
\label{eq:ca-contraction-assumption}
 \sup_{z\in B}
 \|I-J_{C,0}^{-1}D_z\mathcal C(z;\theta,d)\|_2
       \le q<1.
\end{equation}
With $C_S:=\|J_{C,0}^{-1}\|_2/(1-q)$, assume:
\begin{equation}
\label{eq:ca-basin-condition}
 C_S\chi(\delta+C_Me)<\rho.
\end{equation}
Only for the least-squares completion (LS), also assume
$D_z\mathcal C(z;\theta,d)\succ0$ on $B$.
\end{enumerate}

Then, a Global Relation completion $z^*$ exists uniquely on $\R^m$ in the first case and uniquely in $B$ in the second, and:
\begin{equation}
\label{eq:ca-trace-error}
 \|z^*-z^\dagger\|_2
       \le C_S\chi(\delta+C_Me).
\end{equation}
In the least squares completion (LS), it is respectively the unique global minimizer
or a strict local minimizer of the objective $\Psi$ of \eqref{eq:fpinn-ls-objective}.
Furthermore, if a numerical completion $\widetilde z$ satisfies:
\[
 \|\mathcal C(\widetilde z;\theta,d)\|_2\le\mathrm{tol},
\]
and belongs to $B$ in the local case, then, with
$\widetilde u_\theta=u_\theta(\cdot;\widetilde z,d)$, it holds that:
\begin{equation}
\label{eq:ca-total-error}
 \|\widetilde u_\theta-u^\dagger\|
 \le e+C_EC_S
       \bigl[\chi(\delta+C_Me)+\mathrm{tol}\bigr].
\end{equation}
\end{proposition}

\begin{proof}
Plugging $z^\dagger$ into the global relation residual gives:
\[
\begin{aligned}
 r(z^\dagger;\theta,d)
   =\mathsf Kz^\dagger-\mathbf b(d) +\mathcal M[u^\dagger;d]+\mathcal M[u_\theta(\cdot;z^\dagger,d);d] -\mathcal M[u^\dagger;d].
\end{aligned}
\]
Consequently, with $\delta, e$ and $\chi$ defined as above, we obtain:
\begin{equation}
 \|r(z^\dagger;\theta,d)\|_2
       \le\delta+C_Me,
\end{equation}
and hence, depending on the completion rule, we can write: 
\begin{equation}
 \|\mathcal C(z^\dagger;\theta,d)\|_2
       \le\chi(\delta+C_Me),
\end{equation}
with $\chi$ as in \eqref{eq:chi}.

Under assumption~(1), strong monotonicity and the
Cauchy-Schwarz inequality give:
\[
\begin{aligned}
 \mu\|z-w\|_2^2
 &\le
 \langle \mathcal C(z;\theta,d)-\mathcal C(w;\theta,d),z-w\rangle\\
 &\le
 \|\mathcal C(z;\theta,d)-\mathcal C(w;\theta,d)\|_2
 \,\|z-w\|_2.
\end{aligned}
\]
For $z\ne w$, dividing by $\mu\|z-w\|_2$ and using
$C_S=\mu^{-1}$ yields:
\begin{equation}
\label{eq:ca-inverse-bound}
 \|z-w\|_2
 \le C_S\|\mathcal C(z;\theta,d)-\mathcal C(w;\theta,d)\|_2.
\end{equation}
The inequality also holds when $z=w$.

To establish existence of a zero, apply strong monotonicity
with $w=0$. Adding $\langle\mathcal C(0;\theta,d),z\rangle$
to both sides gives the first inequality below, and
Cauchy--Schwarz gives the second:
\[
\begin{aligned}
 \langle\mathcal C(z;\theta,d),z\rangle
 \ge \mu\|z\|_2^2+\langle\mathcal C(0;\theta,d),z\rangle
 \ge \mu\|z\|_2^2-\|\mathcal C(0;\theta,d)\|_2\|z\|_2.
\end{aligned}
\]
On the sphere $\|z\|_2=R$, the last expression is
$R\bigl(\mu R-\|\mathcal C(0;\theta,d)\|_2\bigr)$,
which is strictly positive for
$R>\|\mathcal C(0;\theta,d)\|_2/\mu$.
For a contradiction, suppose that $\mathcal C(\cdot;\theta,d)$ has no zero in the
closed ball $\overline B_R(0)$. Then, the map:
\[
 T(z)=-R\frac{\mathcal C(z;\theta,d)}
                 {\|\mathcal C(z;\theta,d)\|_2}
\]
is continuous and maps $\overline B_R(0)$ into itself, since
$\|T(z)\|_2=R$ for every $z$ in the ball. By Brouwer's fixed-point
theorem, there is a point $\bar z=T(\bar z)$ with $\|\bar z\|_2=R$.
The fixed-point equation gives:
\[
 \langle\mathcal C(\bar z;\theta,d),\bar z\rangle
 =-R\|\mathcal C(\bar z;\theta,d)\|_2<0,
\]
contradicting the strict positivity on the sphere. Thus, a zero exists.
Applying \eqref{eq:ca-inverse-bound} to any two zeros $z_1, z_2$ proves uniqueness.

Under assumption~(2), define
$G(z)=z-J_{C,0}^{-1}\mathcal C(z;\theta,d)$.
For $z,w\in B$, convexity of the ball and the fundamental
theorem of calculus give:
\[
 G(z)-G(w)
   =\int_0^1
       \bigl[I-J_{C,0}^{-1}
           D_z\mathcal C(w+s(z-w);\theta,d)\bigr](z-w)\,ds.
\]
Thus, $\|G(z)-G(w)\|_2\le q\|z-w\|_2$. For $z\in B$, we obtain:
\[
\begin{aligned}
 \|G(z)-z^\dagger\|_2
   \le q\rho
       +\|J_{C,0}^{-1}\|_2\|\mathcal C(z^\dagger;\theta,d)\|_2\le q\rho
       +\|J_{C,0}^{-1}\|_2\,\chi(\delta+C_Me)
    <\rho.
\end{aligned}
\]
By Banach's fixed-point theorem, there exists a unique zero in $B$. Hence, to obtain the coefficient stability estimate under assumption~(2), the definition of $G$ gives:
\[
 z-w=G(z)-G(w)
 +J_{C,0}^{-1}
  \bigl(\mathcal C(z;\theta,d)-\mathcal C(w;\theta,d)\bigr).
\]
For $z,w\in B$, the contraction estimate therefore yields:
\[
 (1-q)\|z-w\|_2
 \le \|J_{C,0}^{-1}\|_2\,
 \|\mathcal C(z;\theta,d)-\mathcal C(w;\theta,d)\|_2.
\]
Thus, \eqref{eq:ca-inverse-bound} also holds on $B$,
with $C_S=\|J_{C,0}^{-1}\|_2/(1-q)$.

Applying \eqref{eq:ca-inverse-bound} to $z^\dagger,z^*$
proves \eqref{eq:ca-trace-error}, and applying it to
$\widetilde z,z^*$ gives
$\|\widetilde z-z^*\|_2\le C_S\mathrm{tol}$.
For rule (LS),
$\mathcal C=\nabla_z\Psi$, so the completion satisfies
$\nabla_z\Psi(z^*;\theta,d)=0$. Under assumption~(1), integrating the strong monotonicity
inequality along the segment from $z^*$ to $z$ gives:
\[
 \Psi(z;\theta,d)-\Psi(z^*;\theta,d)
 \ge \frac{\mu}{2}\|z-z^*\|_2^2,
\]
for $z\in\mathbb R^m$. Hence, $z^*$ is the unique global minimizer. Under assumption~(2), the completion lies in the interior of $B$ and:
\[
 D_z^2\Psi(z^*;\theta,d)
 =D_z\mathcal C(z^*;\theta,d)\succ0.
\]
Together with stationarity, the second-order sufficient
condition implies that $z^*$ is a strict local minimizer.

Finally, the affine lift gives:
\[
 \widetilde u_\theta-u^\dagger
   =u_\theta(\cdot;z^\dagger,d)-u^\dagger
        +\mathcal E(d)(\widetilde z-z^\dagger).
\]
The triangle inequality and the preceding coefficient bounds prove \eqref{eq:ca-total-error}.
\end{proof}

\subsection{Universal approximation}
\label{sec:ca-universal-approximation}

Proposition~\ref{prop:ca-completion-error} transfers an accurate solution approximation through the GR completion. Therefore, to obtain an approximation theorem, we must also show that a gated neural correction can produce such an approximation. 

To this end, let $\mathcal K$ be the compact closure of the spatial or
space-time domain and let $\mathscr N\subset C(\mathcal K)$ be the class
of scalar feedforward networks of arbitrary finite width with a linear
output layer and a fixed activation function $\sigma$. We assume that
$\mathscr N$ is dense in $C(\mathcal K)$, which holds for every
continuous, bounded and nonconstant $\sigma$
\cite{hornik1991approximation}. When, in addition, $\sigma\in C^s$, the
class is dense in $C^s(\mathcal K)$, which is used in
Section~\ref{sec:ca-derivatives-time}.

\begin{lemma}[Gated approximation]
\label{lem:ca-gated-density}
Let $\beta\in C(\mathcal K)$ satisfy $\beta\ge0$, with nonempty zero
set $\mathcal Z_\beta=\beta^{-1}(0)$. Then, for $f,L\in C(\mathcal K)$, the following holds:
\begin{equation}
\label{eq:ca-gated-best}
 \inf_{v\in\mathscr N}\|L+\beta v-f\|_\infty
       =\|L-f\|_{C(\mathcal Z_\beta)}.
\end{equation}
\end{lemma}

\begin{proof}
Notice that the right-hand side is a lower bound because the gated correction vanishes on $\mathcal Z_\beta$.
Set $w=f-L$, $t=\|w\|_{C(\mathcal Z_\beta)}$ and, for $\varepsilon>0$,
$q_\varepsilon=w/(\beta+\varepsilon)$.
Then:
\[
 L+\beta q_\varepsilon-f
       =-\frac{\varepsilon w}{\beta+\varepsilon}.
\]
For any $\epsilon>0$, the compact set
$F_\epsilon=\{X\in \mathcal K:|w(X)|\ge t+\epsilon\}$ is disjoint from $\mathcal Z_\beta$.
If it is nonempty, $\beta$ has a positive minimum $c_\epsilon$
on $F_\epsilon$, and hence:
\[
 \left\|\frac{\varepsilon w}{\beta+\varepsilon}\right\|_\infty
       \le\max\left\{t+\epsilon,\frac{\varepsilon\|w\|_\infty}{c_\epsilon}\right\}.
\]
If $F_\epsilon$ is empty, the norm is at most $t+\epsilon$.
The lower bound on $\mathcal Z_\beta$ therefore implies convergence of
this norm to $t$ as $\varepsilon\downarrow0$.

Then, given $\epsilon>0$, choose $\varepsilon>0$ such that
$ \|L+\beta q_\varepsilon-f\|_\infty<t+\epsilon/2.$
Since $q_\varepsilon\in C(\mathcal K)$, the standard Universal Approximation Theorem for neural networks provides $v\in\mathscr N$ satisfying:
\[
 \|v-q_\varepsilon\|_\infty
 <\frac{\epsilon}{2(1+\|\beta\|_\infty)}.
\]
The triangle inequality then gives
$\|L+\beta v-f\|_\infty<t+\epsilon$ and, since $\epsilon$ is arbitrary and $t$ is the lower bound,
\eqref{eq:ca-gated-best} follows.
\end{proof}

The result above pertains to the approximation with $z=z^\dagger$ fixed. With this and Proposition~\ref{prop:ca-completion-error}, we can now formulate the below Universal Approximation Theorems for the FokasNets, starting with the Fokas-PINN.

\begin{theorem}[Universal approximation for Fokas--PINNs]
\label{thm:ca-fpinn-uat}
Fix the problem data $d$ and a discretization. Let $\mathcal K$ be compact,
$u^\dagger\in\mathcal V\cap C(\mathcal K)$,
$L(z^\dagger;d)\in C(\mathcal K)$, and
$\beta\in C(\mathcal K)$ with $\beta\ge0$.
Assume that the neural correction class is dense in
$C(\mathcal K)$ in the uniform norm and that the corresponding solution
approximations belong to $\mathcal V$ and set:
\begin{equation}
\label{eq:ca-trace-mismatch}
 \mathcal Z_\beta:=\{X\in \mathcal K:\beta(X)=0\},
 \qquad
 \varepsilon_{\rm tr}
 :=\|L(z^\dagger;d)-u^\dagger\|_{C(\mathcal Z_\beta)},
\end{equation}
with $\varepsilon_{\rm tr}=0$ if $\mathcal Z_\beta=\varnothing$. Furthermore, for $\varepsilon_0>0$ and every network satisfying:
\[
 e:=\|u_\theta(\cdot;z^\dagger,d)-u^\dagger\|_\infty
 <\varepsilon_0,
\]
assume the moment estimate and the same completion
alternative of
Proposition~\ref{prop:ca-completion-error}, and that these hold with common bounds
$C_E,C_M,C_S,\chi$, independent of the network.
In the local alternative, also assume that $\rho, q$ and, for rule (LS), the positive-definiteness of $D_z\mathcal C$ on $B$ in Proposition~\ref{prop:ca-completion-error}(2) hold with constants independent of the network. Finally, suppose:
\[
 \varepsilon_{\rm tr}<\varepsilon_0,
 \qquad
 C_S\chi
   (\delta+C_M\varepsilon_{\rm tr})<\rho
 \quad\text{in the local alternative},
\]
where $a:=1+C_EC_S\chi C_M$, and $c:=C_EC_S\chi$.

Then, for every $\epsilon>0$ satisfying:
\[
 \varepsilon_{\rm tr}+\epsilon<\varepsilon_0,
 \qquad
 C_S\chi
   [\delta+C_M(\varepsilon_{\rm tr}+\epsilon)]
 <\rho
 \quad\text{in the local alternative},
\]
there exist a finite network and a completion
$z^*(\theta,d)$ supplied by
Proposition~\ref{prop:ca-completion-error} such that:
\begin{equation}
\label{eq:ca-fpinn-uat-bound}
 \|u_\theta-u^\dagger\|_\infty
 \le a(\varepsilon_{\rm tr}+\epsilon)+c\delta.
\end{equation}
If, along a sequence of refinements satisfying
these assumptions, we have:
\begin{equation}
\label{eq:ca-fpinn-refinement}
 a\varepsilon_{\rm tr}+c\delta
 \rightarrow0,
\end{equation}
then, for every $\varepsilon>0$, there exist a
discretization and a finite network with a completion
such that:
\[
 \|u_\theta-u^\dagger\|_\infty<\varepsilon.
\]
\end{theorem}

\begin{proof}
Apply Lemma~\ref{lem:ca-gated-density} with
$L=L(z^\dagger;d)$, $f=u^\dagger$, and
$\beta=\beta$. For every admissible $\epsilon>0$,
the lemma gives a finite network satisfying $e<\varepsilon_{\rm tr}+\epsilon<\varepsilon_0$. In the local alternative, this also gives $C_S\chi(\delta+C_Me)<\rho$.

Thus, Proposition~\ref{prop:ca-completion-error}
provides a completion for this network and 
its error estimate with $\mathrm{tol}=0$ yields
\[
\begin{aligned}
 \|u_\theta-u^\dagger\|_\infty
 \le e+C_EC_S\chi
       (\delta+C_Me)
 =ae+c\delta
 \le a(\varepsilon_{\rm tr}+\epsilon)
       +c\delta,
\end{aligned}
\]
proving \eqref{eq:ca-fpinn-uat-bound}.

Now let $\varepsilon>0$. By
\eqref{eq:ca-fpinn-refinement}, choose a discretization
such that:
\[
 a\varepsilon_{\rm tr}+c\delta
 <\frac{\varepsilon}{2}.
\]
The strict neighborhood and basin inequalities allow
a sufficiently small $\epsilon>0$ satisfying the statement's
conditions and $a\epsilon<\frac{\varepsilon}{2}$. Applying \eqref{eq:ca-fpinn-uat-bound} gives
$\|u_\theta-u^\dagger\|_\infty<\varepsilon$.
\end{proof}

We next extend the approximation result to a shared
branch-trunk \cite{lu2021learning} correction over a data family. This requires
uniform completion estimates and an encoder that represents
the correction required by the family.
Let $\mathcal D_c\subset\mathcal D$
be a compact metric data family. Suppose the selected solution
map is continuous into $C(\mathcal K)$, $z^\dagger(d)$ is continuous,
and the lift is jointly continuous in $(z,d)$ into $C(\mathcal K)$.
Keep $\mathcal K$, the gate and its zero set fixed over the family
at each resolution. Take $\delta$ and
$\varepsilon_{\rm tr}$ as suprema over $\mathcal D_c$,
and use constants $a=1+C_EC_S\chi C_M$ and $c=C_EC_S\chi$ independent of $d$ and of the network.

Let $\mathrm{enc}:\mathcal D_c\to\mathcal D_{\mathrm{enc}}\subset\R^{N_{\mathrm{enc}}}$
be the continuous encoder used by the branch network,
with $\mathcal D_{\mathrm{enc}}=\mathrm{enc}(\mathcal D_c)$.
Suppose there is $\hat w\in C(\mathcal D_{\mathrm{enc}}\times \mathcal K)$
satisfying:
\begin{equation}
\label{eq:ca-encoding-error}
 \sup_{d\in\mathcal D_c,\,X\in \mathcal K}
 \left|
  u^\dagger(X;d)-L(z^\dagger(d);d)(X)
       -\hat w(\mathrm{enc}(d),X)
 \right|
       \le\varepsilon_{\rm enc}.
\end{equation}
Here $w^\dagger(d,X):=u^\dagger(X;d)-L(z^\dagger(d);d)(X)$ is the
correction that the gated network must represent for the datum $d$, and
\eqref{eq:ca-encoding-error} asks that it be reproducible, up to
$\varepsilon_{\rm enc}$, by a continuous function of the encoding alone. (The exact form of the lifted
operator, with the datumwise completion inside the lift, is written out explicitly in \eqref{eq:numerics-op-fields}--\eqref{eq:numerics-op-trace}
of Section~\ref{sec:numerics}.)

\begin{theorem}[Universal approximation for FONets]
\label{thm:ca-operator-uat}
Assume the compact-family and encoding conditions above.
For every branch--trunk correction in the approximation neighborhood
of Theorem~\ref{thm:ca-fpinn-uat}, suppose its moment and
completion conditions hold uniformly over $\mathcal D_c$.
Assume also that $\mathcal C(z;\theta,d)$ is jointly
continuous in $(z,d)$ for each fixed correction network.
Require:
\[
 \varepsilon_{\rm tr}+\varepsilon_{\rm enc}<\varepsilon_0,
\]
and, in the local case:
\[
 C_S\chi
  [\delta+C_M
       (\varepsilon_{\rm tr}+\varepsilon_{\rm enc})]
       <\rho.
\]
For every sufficiently small $\epsilon>0$, one finite
branch--trunk network, with datum-wise completion, satisfies:
\begin{equation}
\label{eq:ca-operator-uat-bound}
 \sup_{d\in\mathcal D_c}
 \|\mathcal S_{\theta}(d)-u^\dagger(\cdot;d)\|_\infty
 \le a
    (\varepsilon_{\rm tr}+\varepsilon_{\rm enc}+\epsilon)
      +c\delta.
\end{equation}
The completed operator is continuous on $\mathcal D_c$.
If the assumptions hold along refinements and moreover $a(\varepsilon_{\rm tr}+\varepsilon_{\rm enc}) +c\delta\rightarrow 0$, then, the completed operators approximate the selected solution map arbitrarily accurately, uniformly over $\mathcal D_c$. The branch--trunk rank and network widths may increase with the requested accuracy.
\end{theorem}

\begin{proof}
Finite sums $\sum_{j=1}^{p_{\mathrm{op}}} a_j(p)b_j(X)$ of continuous
factors form an algebra on the compact set $\mathcal D_{\mathrm{enc}}\times \mathcal K$.
They contain constants and separate points, so
Stone--Weierstrass gives their density in
$C(\mathcal D_{\mathrm{enc}}\times \mathcal K)$.
Approximate each factor by a scalar network.
For one product,
\[
 \|a_jb_j-B_jT_j\|_\infty
 \le\|a_j\|_\infty\|b_j-T_j\|_\infty
      +\|a_j-B_j\|_\infty\|T_j\|_\infty.
\]
Since the sum is finite, arbitrarily small factor errors
give an arbitrarily small total error.
Combining the factors into multioutput branch and trunk
networks proves density of the branch-trunk class used
in the operator construction
\cite{chen1995universal,lu2021learning}.

Put
$w(d,X)=u^\dagger(X;d)-L(z^\dagger(d);d)(X)$.
For $\varepsilon>0$, define
$q_\varepsilon(p,X)=\hat w(p,X)/(\beta(X)+\varepsilon)$.
Then, we have:
\[
\begin{aligned}
 \beta(X)q_\varepsilon(\mathrm{enc}(d),X)-w(d,X)
   = -\frac{\varepsilon w(d,X)}
                  {\beta(X)+\varepsilon}+\frac{\beta(X)}
                  {\beta(X)+\varepsilon}
         \bigl(\hat w(\mathrm{enc}(d),X)-w(d,X)\bigr).
\end{aligned}
\]
The first term has uniform norm tending to
$\varepsilon_{\rm tr}$, by the compact-set argument
of Lemma~\ref{lem:ca-gated-density} applied on
$\mathcal D_c\times \mathcal K$.
The second is bounded by $\varepsilon_{\rm enc}$.
Choose $\varepsilon$ and then one branch--trunk approximation
to $q_\varepsilon$ so that:
\[
 \sup_{d\in\mathcal D_c}
 \|u_\theta(\cdot;z^\dagger(d),d)-u^\dagger(\cdot;d)\|_\infty
       <\varepsilon_{\rm tr}+\varepsilon_{\rm enc}+\epsilon.
\]
Choose $\epsilon$ within the uniform approximation and basin neighborhoods.
Proposition~\ref{prop:ca-completion-error} then gives a
completion for every datum and proves
\eqref{eq:ca-operator-uat-bound}.

For continuity, fix this network and let $d_n\to d$.
The coefficient error estimate and continuity of
$z^\dagger(d_n)$ imply boundedness of the completions.
Any convergent subsequence has a limit $\bar z$ satisfying
$\mathcal C(\bar z;\theta,d)=0$ by joint continuity.
In the local case it also belongs to
$\overline B_{\rho}(z^\dagger(d))$.
Uniqueness in that ball, or global uniqueness in the first
alternative, identifies the limit with $z^*(\theta,d)$.
Thus the full coefficient sequence converges.
Continuity of the lift and correction proves continuity
of the completed operator.
Finally, choose the resolution and then $\epsilon$ as in
Theorem~\ref{thm:ca-fpinn-uat} to obtain arbitrary accuracy.
\end{proof}

\subsection{Derivative errors}
\label{sec:ca-derivatives-time}

Theorem~\ref{thm:ca-fpinn-uat} controls the completed approximation in
the uniform norm. The quantities minimized in training, however, are the
PDE residual and the residuals of the prescribed traces that the lift does
not encode, which involve derivatives of $u_\theta$. The following
corollary shows that the same argument, carried out in a norm that
controls these derivatives, yields convergence of the residuals for both
problem classes.

\begin{corollary}[PDE residual convergence]
\label{cor:ca-stronger-norm}
Let $\|\cdot\|_\star$ denote the norm of $C^2(\overline\Om)$ for the
elliptic problem \eqref{eq:bvp} and of $C^{1,r}(\overline Q_{\mathcal T})$
for the evolution problem \eqref{eq:evo-bvp}. Suppose the completion
assumptions of Theorem~\ref{thm:ca-fpinn-uat} hold with every norm
replaced by $\|\cdot\|_\star$, and that there is a correction
$v^\dagger$ with:
\begin{equation}
\label{eq:ca-strong-reference}
 \|L(z^\dagger;d)+\beta v^\dagger-u^\dagger\|_\star\le\varepsilon_{\rm ref},
 \qquad
 \inf_{\theta\in\Theta}\|\beta(v_\theta-v^\dagger)\|_\star=0,
\end{equation}
where $\varepsilon_{\rm ref}$ plays the role of $\varepsilon_{\rm tr}$ in
$\|\cdot\|_\star$ (Lemma~\ref{lem:ca-gated-density} is specific to the
uniform norm and no longer provides it). Assume
$\varepsilon_{\rm ref}<\varepsilon_0$ and, in the local case,
$C_S\chi(\delta+C_M\varepsilon_{\rm ref})<\rho$. Then, for every
$\epsilon>0$ satisfying the same inequalities with
$\varepsilon_{\rm ref}+\epsilon$ in place of $\varepsilon_{\rm ref}$, there
is a finite network with:
\begin{equation}
\label{eq:ca-strong-bound}
 \|u_\theta-u^\dagger\|_\star\le a(\varepsilon_{\rm ref}+\epsilon)+c\,\delta,
\end{equation}
with the constants $a,c$ measured in $\|\cdot\|_\star$. Consequently, if
$R(\cdot;\lambda)$, respectively $\mathcal N[\cdot;\lambda]$, is locally
Lipschitz from $\|\cdot\|_\star$ into the uniform norm, then:
\begin{equation}
\label{eq:ca-elliptic-residual-bound}
 \|\mathcal R_{\mathrm{PDE}}[u_\theta;d]\|_\infty
   \le C_{\mathcal R}\,\|u_\theta-u^\dagger\|_\star,
 \qquad
 \|\boldsymbol\gamma_\pre u_\theta-\boldsymbol d_\partial\|_\infty
   \le\|u_\theta-u^\dagger\|_\star,
\end{equation}
where $C_{\mathcal R}$ depends only on the differential operator and on
the local Lipschitz constant of the nonlinearity (for the elliptic
problem, $C_{\mathcal R}=n_x+C_R$ with $C_R$ the Lipschitz constant of
$R(\cdot;\lambda)$ near the range of $u^\dagger$). In particular, the PDE
residual and the residuals of the prescribed traces not encoded by the
lift tend to zero as $\|u_\theta-u^\dagger\|_\star\to0$.
\end{corollary}

\begin{proof}
Choose a network with $\|\beta(v_\theta-v^\dagger)\|_\star<\epsilon$, so
that $e<\varepsilon_{\rm ref}+\epsilon$ by the triangle inequality. The
proof of Theorem~\ref{thm:ca-fpinn-uat}, with
Proposition~\ref{prop:ca-completion-error} applied in $\|\cdot\|_\star$,
gives \eqref{eq:ca-strong-bound}. 

Since $\mathcal R_{\mathrm{PDE}}[u^\dagger;d]=0$ and
$\boldsymbol\gamma_\pre u^\dagger=\boldsymbol d_\partial$, the first bound
in \eqref{eq:ca-elliptic-residual-bound} follows from the local Lipschitz
continuity of $\mathcal R_{\mathrm{PDE}}$ with respect to $\|\cdot\|_\star$;
for the elliptic problem, with $w=u_\theta-u^\dagger$,
$\|-\Delta w\|_\infty\le n_x\|w\|_{C^2}$ and
$\|R(u_\theta;\lambda)-R(u^\dagger;\lambda)\|_\infty\le C_R\|w\|_\infty$.
Finally, the second bound holds because the prescribed traces not encoded by the
lift are derivatives of $u_\theta-u^\dagger$ of order at most $r-1$,
which $\|\cdot\|_\star$ controls.
\end{proof}

\begin{remark}
Note that hypothesis \eqref{eq:ca-strong-reference} holds, for instance, when $\beta\in C^s$ and $v^\dagger$ admits a $C^s$ extension to a neighbourhood
of the closed domain, with $s=2$ for the elliptic and $s=\max(1,r)$ for
the evolution problem, by the $C^s$ density of $\mathscr N$ \cite{hornik1991approximation} and the product rule.  
\end{remark}

\section{Elliptic and evolution FokasNet constructions}
\label{sec:pde-constructions}

In this Section, we make explicit the FokasNet constructions for the
elliptic and evolution PDE classes, i.e., the trace bases, lifts, gates,
tests and quadratures that enter the general assembly of
Section~\ref{sec:fpinn-finite-gr}, together with the resulting matrices
and Jacobians. These realizations are used in the numerical experiments
of Section~\ref{sec:numerics}. Throughout, $d$ and the discretization are
fixed and their dependence is suppressed where unambiguous.

\subsection{Elliptic problems}
\label{sec:elliptic-linear}\label{prop:elliptic-realisation}

For the elliptic problem \eqref{eq:bvp}, we fix boundary interpolation nodes $x_i^\partial\in\bd$, $i=1,\ldots,M$ and partition their indices into the sets $\pre$ and $\un$ of prescribed and missing boundary values, with $m=|\un|$. We choose functions $\ell_0,\ell_1,\ldots,\ell_M\in\mathcal V$
and a smooth gate $\beta$ with:
\[
 \ell_0(x_i^\partial)=0,\qquad \ell_j(x_i^\partial)=\delta_{ij},
 \qquad i,j=1,\ldots,M,\qquad
 \beta|_{\bd}=0,\quad \beta>0\ \text{in }\Om .
\]
Writing $\hat g(z)\in\R^M$ for the nodal values,
$\hat g_\pre(z)=\hat g_\pre:=(g_D(x_i^\partial))_{i\in\pre}$ and
$\hat g_\un(z)=z$, the affine lift can be written as: 
\begin{equation}\label{eq:affine-lift-method}
 L(z;d)(x)=\ell_0(x)+\sum_{j=1}^M\ell_j(x)\,\hat g_j(z).
\end{equation}
We refer to \eqref{eq:affine-lift-method} as a \emph{nodal lift}: the
trace is represented by its values at the boundary nodes, which the
cardinal functions $\ell_j$ interpolate and extend into the domain.
By the interpolation properties, $L(z;d)(x_i^\partial)=\hat g_i(z)$,
and since $\beta|_{\bd}=0$ the lifted approximation
$u_\theta(x;z,d)=L(z;d)(x)+\beta(x)v_\theta(x)$ satisfies
$u_\theta(x_i^\partial;z,d)=g_D(x_i^\partial)$ for $i\in\pre$ and
$u_\theta(x_i^\partial;z,d)=z_i$ for $i\in\un$, for every $z$ and
$\theta$. 

We now consider the test function and quadrature for the GR, in accordance with Section \ref{sec:fpinn-finite-gr}. We use real, harmonic tests $\{\varphi_k\}_{k=1}^p\subset\mathcal H_E$, a
boundary quadrature and interior quadrature nodes $\{x_a\}_{a=1}^{n_A}$
with weights $\omega_a>0$, and define:
\begin{equation}\label{eq:elliptic-matrices}
\begin{aligned}
 A_{aj}&=\ell_j(x_a),\quad A\in\R^{n_A\times M},
 &\qquad H_{ak}&=\varphi_k(x_a),\quad Q_0=H^\top\operatorname{diag}(\omega_a)\in\R^{p\times n_A},\\
 K_{kj}&=\qint_{\bd}[\ell_j\,\dn\varphi_k],\quad K\in\R^{p\times M},
 &\qquad (C_\pre)_{kl}&=\qint_{\bd}[\varphi_k\zeta'_l],\quad C_\pre\in\R^{p\times m'},
\end{aligned}
\end{equation}
where $\zeta'_l$, $l=1,\ldots,m'$, are the basis functions of the unknown  flux on $\GD$. In the above, for the matrix $X$ whose columns are indexed by the boundary nodes, $X_{:,\pre}$ and $X_{:,\un}$ denote the column blocks indexed by $\pre$ and $\un$, in the ordering of the nodes. The boundary value and flux of the approximation are
$g=\ell_0|_{\bd}+\sum_j\ell_j|_{\bd}\,\hat g_j(z)$ and
$w=w_0+\sum_l\zeta'_lz'_l$, where $w_0$ is the prescribed flux (extended
by zero on $\GD$). Inserting them into the discretized form of
\eqref{eq:GR} and collecting the coefficients of $z$, $z'$ and the data gives the system:
\begin{equation}\label{eq:general-joint-system}
 K_{:,\un}\,z-C_\pre\,z'+Q_0\,\mathbf q=\mathbf f,
 \qquad
 \mathbf f_k=\qint_{\bd}[\varphi_kw_0-\ell_0\,\dn\varphi_k]-(K_{:,\pre}\hat g_\pre)_k,
\end{equation}
where $\mathbf q \in \mathbb{R}^{n_A}$ with entries $q(x_a)$ is the sampled source and
$(Q_0\mathbf q)_k=\sum_a\omega_a\varphi_k(x_a)q_a$. This is precisely \eqref{eq:fpinn-joint-gr} with
$\mathsf B=K_{:,\un}$ and $\mathsf B'=-C_\pre$. With the annihilator $\mathsf T$ of $C_\pre$ as in \eqref{eq:fpinn-annihilator}, the reduced quantities \eqref{eq:fpinn-reduced-matrices} are
$\mathsf K=\mathsf TK_{:,\un}\in\R^{n\times m},
\mathbf b=\mathsf T\mathbf f\in\R^n,
Q=\mathsf TQ_0\in\R^{n\times n_A}$, with $n=p-\operatorname{rank}C_\pre$, and multiplication of \eqref{eq:general-joint-system} by $\mathsf T$ gives the reduced system:
\begin{equation}\label{eq:general-reduced-system}
 \mathsf Kz+Q\,\mathbf q=\mathbf b,
\end{equation}
which is the elliptic form of \eqref{eq:reduced}.

In our setting, the source $\mathbf q$ depends on the solution approximation. Evaluating the lift at the quadrature nodes gives the vector $\mathbf u_A(z) \in \mathbb{R}^{n_A}$, with entries:
\[
\begin{aligned}
(\mathbf u_A(z))_a =  u_\theta(x_a;z,d)
  =A\hat g(z)+\ell_0(x_a)+\beta(x_a)v_\theta(x_a)
  =A_{:,\un}z+ a_\theta(x_a),
\end{aligned}
\]
with $a_\theta:=A_{:,\pre}\hat g_\pre+\ell_0(x_a)+\beta(x_a)v_\theta(x_a)$, for $a=1, \dots, n_A$.
The sampled source is then given by:
\begin{equation}
  \mathbf q=R(\mathbf u_A(z);\lambda)-\mathbf F 
\end{equation}
where $\mathbf F \in \mathbb{R}^{n_A}$, with entries $\mathbf F_a = F(x_a)$ and $R$ acting componentwise. Inserting it into \eqref{eq:general-reduced-system}, the finite GR
residual \eqref{eq:fpinn-moment-residual} becomes:
\begin{equation}\label{eq:residual}
 r(z;\theta,d)=\mathsf Kz-\mathbf b+Q\bigl(R(\mathbf u_A(z);\lambda)-\mathbf F\bigr),
\end{equation}
i.e., $\mathcal M(z;\theta,d)=Q(R(\mathbf u_A(z);\lambda)-\mathbf F)$.
For $R(\cdot;\lambda)\in C^1$, the chain rule gives
$D_z\mathbf q=D_RA_{:,\un}$, and the Jacobian becomes:
\begin{equation}\label{eq:jacobian}
 J_r(z;\theta,d)=\mathsf K+Q\,D_R\,A_{:,\un},
 \qquad
 D_R=\operatorname{diag}\bigl(R'((\mathbf u_A(z))_a;\lambda)\bigr)_{a=1}^{n_A}.
\end{equation}

For the mixed elliptic problem, the loss \eqref{eq:prelim-physics-objective}
reduces to
\begin{equation}\label{eq:elliptic-training-loss}
 \mathcal J(\theta;z,d)
   =\lambda_{\mathrm{PDE}}\|-\Delta u_\theta+R(u_\theta;\lambda)-F\|_{\mathscr S_{\mathrm{PDE}}}^2
   +\lambda_N\|\dn u_\theta-w_N\|_{\mathscr S_N}^2,
\end{equation}
since the Dirichlet values are exact at the nodes $x_i^\partial$,
$i\in\pre$ and the derivatives entering it are calculated using the lifted representation and are given by:
\begin{equation}\label{eq:spatial-derivatives-method}
 \Delta u_\theta=\Delta L(z;d)+\beta\Delta v_\theta+2\nabla\beta\cdot\nabla v_\theta+v_\theta\Delta\beta,
 \qquad
 \dn u_\theta|_{\bd}=\dn L(z;d)+v_\theta\,\dn\beta.
\end{equation}

\subsection{Evolution problems}\label{sec:evolution}
\label{prop:evolution-realisation}

For the evolution problem \eqref{eq:evo-bvp} we give two constructions:
the direct construction, using the space-time Global Relation and the time-stepping construction,
in which each implicit stage is a stationary problem of the type of
Section~\ref{sec:elliptic-linear}.

\subsubsection{Space-time Global Relation construction}
Let $\mathcal G:\mathcal V\to\R^M$ be the linear operator extracting $M$
spatial trace coefficients (e.g., nodal boundary values), with $\pre$ and
$\un$ the index sets of the prescribed and missing coefficients, and let $\widehat{\boldsymbol\gamma}_\pre(t;d)$ denote the
prescribed coefficients, with
$\widehat{\boldsymbol\gamma}_\pre(0;d)=(\mathcal Gu_0)_\pre$. We choose
spatial extension functions $\ell_j$ with $\mathcal G\ell_j=e_j$,
$j=1,\ldots,M$, temporal basis functions
$\boldsymbol\psi=(\psi_1,\ldots,\psi_{M_t})^\top$ with
$\boldsymbol\psi(0)=0$, and gates $\beta(x)$, $\beta_{\rm t}(t)$ with
$\beta_{\rm t}(0)=0$ and $\mathcal G[\beta v]=0$ for every smooth $v$, all
with the regularity required for the approximation to belong to
$\mathcal V$. The missing coefficients are parameterized by:
\begin{equation}\label{eq:evo-trace-basis}
 \widehat{\boldsymbol\gamma}_\un(t;z,d)=(\mathcal Gu_0)_\un+Z\boldsymbol\psi(t),
 \qquad Z\in\R^{|\un|\times M_t},\quad z=\operatorname{\bf vec}(Z)\in\R^m, 
\end{equation}
with $m=|\un|\,M_t$, so that each of the $|\un|$ missing spatial coefficients carries $M_t$
temporal degrees of freedom. The lift is then:
\begin{equation}\label{eq:evo-lift}
 L(z;d)(x,t)=u_0(x)+\sum_{j=1}^M\ell_j(x)\bigl(\widehat\gamma_j(t;z,d)-(\mathcal Gu_0)_j\bigr).
\end{equation}
Since $\boldsymbol\psi(0)=0$ and
$\widehat{\boldsymbol\gamma}_\pre(0;d)=(\mathcal Gu_0)_\pre$, we have
$\widehat{\boldsymbol\gamma}(0;z,d)=\mathcal Gu_0$ and hence,
$L(z;d)(\cdot,0)=u_0$, while $\mathcal G\ell_j=e_j$ gives
$\mathcal GL(z;d)(\cdot,t)=\widehat{\boldsymbol\gamma}(t;z,d)$ by
linearity. The conditions $\beta_{\rm t}(0)=0$ and $\mathcal G[\beta v]=0$
transfer both properties to the lifted approximation:
\begin{equation}\label{eq:evo-field}
 u_\theta(x,t;z,d)=L(z;d)(x,t)+\beta(x)\beta_{\rm t}(t)\,v_\theta(x,t;d),
\end{equation}
which therefore satisfies $u_\theta(\cdot,0;z,d)=u_0$ and
$\mathcal Gu_\theta(\cdot,t;z,d)=\widehat{\boldsymbol\gamma}(t;z,d)$ for every $z$, $\theta$ and $t\in[0,\mathcal T]$. By construction, the
prescribed trace coefficients equal $\widehat{\boldsymbol\gamma}_\pre(t;d)$,
so that the initial condition and the prescribed trace coefficients extracted by $\mathcal G$ are exact.

Fix tests $\varphi_k\in\mathcal H_T$ with upper limits $\tau_k$,
$k=1,\ldots,p$, and split the evolution moment
\eqref{eq:fpinn-unreduced-moments} into a source and a terminal part to obtain:
\begin{equation}\label{eq:evo-source-moment}
 (\mathbf m_{\rm nl}[u;d])_k=\qint_{\Om\times(0,\tau_k)}[\varphi_k(\mathcal N[u;\lambda]-F)],
 \qquad
 (\mathbf m_{\rm end}[u;d])_k=\qint_\Om[\varphi_k(\cdot,\tau_k)\,u(\cdot,\tau_k)].
\end{equation}
The trace collection of $u_\theta$ is affine in $z$, with space-time basis functions given by the boundary traces of $\ell_j(x)\psi_i(t)$, so the assembly of Section~\ref{sec:fpinn-finite-gr} applies to the
quadrature form of \eqref{eq:prelim-evo-GR} and the assembled system \eqref{eq:fpinn-joint-gr} becomes: 
\begin{equation}\label{eq:evo-joint-system}
 \mathsf Bz+\mathsf B'z'
  +\mathbf m_{\rm nl}[u_\theta(\cdot;z,d);d]+\mathbf m_{\rm end}[u_\theta(\cdot;z,d);d]
  =\mathbf f(d),
\end{equation}
where the matrices are those of \eqref{eq:fpinn-assembly-entries} for
the present trace basis and tests, and differ from their elliptic
counterparts in both entries and dimensions. Indexing the columns of
$Z$ by pairs $(j,i)$, $j\in\un$, $i=1,\ldots,M_t$, in the order of
$\operatorname{\bf vec}$, the trace basis function $\boldsymbol\zeta_{(j,i)}$
is the boundary trace collection of $\ell_j(x)\psi_i(t)$, and with
$\boldsymbol\zeta'_l(x,t)$, $l=1,\ldots,m'$, a basis of the auxiliary
unknown traces,
\begin{equation}\label{eq:evo-boundary-vector}
\begin{aligned}
 (\mathsf B)_{k,(j,i)}&=-\qint_{\bd\times(0,\tau_k)}
   [\mathscr B_{\mathcal L}[\boldsymbol\zeta_{(j,i)},\varphi_k]],
 &\qquad \mathsf B&\in\R^{p\times m},\quad m=|\un|M_t,\\
 (\mathsf B')_{kl}&=-\qint_{\bd\times(0,\tau_k)}
   [\mathscr B_{\mathcal L}[\boldsymbol\zeta'_l,\varphi_k]],
 &\qquad \mathsf B'&\in\R^{p\times m'},\\
 (\mathbf f(d))_k&=\qint_\Om[\varphi_k(\cdot,0)u_0]
   +\qint_{\bd\times(0,\tau_k)}
   [\mathscr B_{\mathcal L}[\boldsymbol t^{\rm dat}(d),\varphi_k]],
 &\qquad \mathbf f(d)&\in\R^p,
\end{aligned}
\end{equation}
the first term of $\mathbf f$ being the initial moment $\mathbf c^0(d)$
of \eqref{eq:fpinn-known-initial-moments}. With $\mathsf T\in\R^{n\times p}$
the annihilator of $\mathsf B'$ as in \eqref{eq:fpinn-annihilator},
$n=p-\operatorname{rank}\mathsf B'$, the reduced quantities are
$\mathsf K=\mathsf T\mathsf B\in\R^{n\times m}$,
$\mathbf b(d)=\mathsf T\mathbf f(d)\in\R^n$,
$\mathcal M_{\rm nl}:=\mathsf T\mathbf m_{\rm nl}$ and
$\mathcal M_{\rm end}:=\mathsf T\mathbf m_{\rm end}$, and the finite GR
residual \eqref{eq:fpinn-moment-residual} is:
\begin{equation}\label{eq:evo-residual}
 r(z;\theta,d)=\mathsf Kz-\mathbf b(d)
   +\mathcal M_{\rm nl}[u_\theta(\cdot;z,d);d]
   +\mathcal M_{\rm end}[u_\theta(\cdot;z,d);d].
\end{equation}
At the spatial quadrature nodes $\{x_a\}_{a=1}^{n_A}$, define the matrix
$A\in\R^{n_A\times M}$ with entries $A_{aj}=\ell_j(x_a)$, and the vector
of sampled values of the approximation \eqref{eq:evo-field},
$\mathbf u_A(t)\in\R^{n_A}$, with entries:
\[
 (\mathbf u_A(t))_a=u_\theta(x_a,t;z,d)
  =\bigl(A_{:,\un}Z\boldsymbol\psi(t)\bigr)_a+(\mathbf a_\theta(t))_a,
\]
where, by \eqref{eq:evo-trace-basis} and \eqref{eq:evo-lift}, the vector
$\mathbf a_\theta(t)\in\R^{n_A}$ has entries:
\[
 (\mathbf a_\theta(t))_a
  =u_0(x_a)
   +\Bigl(A_{:,\pre}\bigl(\widehat{\boldsymbol\gamma}_\pre(t;d)-(\mathcal Gu_0)_\pre\bigr)\Bigr)_a
   +\beta_{\rm t}(t)\,\beta(x_a)\,v_\theta(x_a,t;d),
 \qquad a=1,\ldots,n_A,
\]
and collects the $z$-independent contributions (initial datum, prescribed
trace coefficients and gated correction).
Hence, the terminal moment is affine in $z$, with a linear part
$\mathsf J_{\rm end}:=D_z\mathcal M_{\rm end}$ that depends neither on $z$
nor on $\theta$. It follows that, if $z\mapsto\mathcal M_{\rm nl}[u_\theta(\cdot;z,d);d]$
is continuously differentiable, then $r(\cdot;\theta,d)\in C^1(\R^m;\R^n)$ with:
\begin{equation}\label{eq:evo-jacobian}
 J_r(z;\theta,d)=\mathsf K+D_z\mathcal M_{\rm nl}[u_\theta(\cdot;z,d);d]+\mathsf J_{\rm end}.
\end{equation}
Note that, when $\mathcal N[u]=\partial_xG(u)$ and the endpoint values of $G(u)$
involve only prescribed traces, the chain rule equivalently gives $\qint_{\Om\times(0,\tau_k)}[-\partial_x\varphi_k\,G(u)-\varphi_kF]$ for the source moment,
and the endpoint term is moved to the data vector $\mathbf f$ (this approach is used for the Korteweg-De Vries (KdV) equation in Section \ref{sec:numerics}).

\subsubsection{Completion of the space-time relation via the compression map (P)}
\label{sec:evo-completion-P}
In the described space-time construction, the tested equations are obtained by pairing each adjoint test with several upper time limits, so
$r(z;\theta,d)=0$ is an overdetermined nonlinear system. According to
Section~\ref{sec:training-objective}, we must therefore select a square
defining map $\mathcal C$. In the least-squares rule (LS), the Jacobian
\eqref{eq:fpinn-ls-hessian} contains second derivatives of the source
moments with respect to $z$, which are expensive to form over the whole
space-time domain. Hence, we instead apply the rule (P) to compress the system by an appropriate projection rule, as described below.

By \eqref{eq:evo-jacobian}, the terminal moment is affine in $z$,
$\mathcal M_{\rm end}[u_\theta(\cdot;z,d);d]=\mathsf J_{\rm end}z
+\mathcal M_{\rm end}[u_\theta(\cdot;0,d);d]$, so the residual
\eqref{eq:evo-residual} splits into a part that is linear in $z$ and the
nonlinear source moment:
\begin{equation}\label{eq:evo-end-affine}
 r(z;\theta,d)=\mathsf K_\star z-\mathbf b_\star(\theta,d)
   +\mathcal M_{\rm nl}[u_\theta(\cdot;z,d);d],
 \qquad
 \mathsf K_\star:=\mathsf K+\mathsf J_{\rm end}\in\R^{n\times m},
\end{equation}
with $\mathbf b_\star(\theta,d):=\mathbf b(d)-\mathcal M_{\rm end}[u_\theta(\cdot;0,d);d]\in\R^n$.
The matrix $\mathsf K_\star$ collects the complete linear action of the unknown traces on the Global Relations.
If the source moment did not depend on $z$, \eqref{eq:evo-end-affine}
would be an overdetermined linear system for $z$, whose least-squares
solution is $z=\mathsf K_\star^{+}(\mathbf b_\star-\mathcal M_{\rm nl})$,
where:
\[
 \mathsf K_\star^{+}:=(\mathsf K_\star^\top\mathsf K_\star)^{-1}\mathsf K_\star^\top\in\R^{m\times n},
\]
is its Moore-Penrose pseudo-inverse (assuming that $\mathsf K_\star$ has full column rank). However, since the source moment does depend on $z$, we define the completion as the fixed point:
\[
 z=\mathsf K_\star^{+}\bigl(\mathbf b_\star(\theta,d)
   -\mathcal M_{\rm nl}[u_\theta(\cdot;z,d);d]\bigr).
\]
Indeed, since $\mathsf K_\star^{+}\mathsf K_\star=I$, applying $\mathsf K_\star^{+}$ to \eqref{eq:evo-end-affine} gives
$\mathsf K_\star^{+}r(z;\theta,d)=z-\mathsf K_\star^{+}\bigl(\mathbf b_\star(\theta,d)-\mathcal M_{\rm nl}[u_\theta(\cdot;z,d);d]\bigr)$,
so the fixed point is exactly a zero of $\mathsf K_\star^{+}r$, i.e., rule (P) of \eqref{eq:fpinn-defining-map} with $\Pi=\mathsf K_\star^{+}$.

The advantage of this choice is visible in its Jacobian, 
$J_C=D_z\mathcal C=\mathsf K_\star^{+}(\mathsf K_\star+D_z\mathcal M_{\rm nl})
=I+\mathsf K_\star^{+}D_z\mathcal M_{\rm nl}$, which involves only first
derivatives of the source moment and is a perturbation of the identity. Hence, 
if $q:=\|\mathsf K_\star^{+}D_z\mathcal M_{\rm nl}\|_2<1$, then $J_C$ is
nonsingular with $\operatorname{cond}_2(J_C)\le(1+q)/(1-q)$. We can solve $C(z;\theta,d)=\mathsf K_\star^{+}\,r(z;\theta,d)=0$ by Newton's method:
\begin{equation}\label{eq:evo-projected-newton}
 z^{j+1}=z^j-J_C(z^j;\theta,d)^{-1}\mathcal C(z^j;\theta,d),
 \qquad
 J_C=I+\mathsf K_\star^{+}D_z\mathcal M_{\rm nl}[u_\theta(\cdot;z^j,d);d].
\end{equation}

\subsection{Evolution problems: time-stepping construction}
\label{prop:evolution-time-stepping}

An alternative implementation for evolution problems discretizes time first and applies the stationary Green identity of the stage operator $\mathcal L_x+\varsigma_\ell/\Delta t_\ell$ to each implicit stage, which is \eqref{eq:GR} for the heat stages and the third-order identity \eqref{eq:app-kdv-concomitant} for the KdV stages. Specifically, let
$0=t_0<t_1<\cdots<t_{N_T}=\mathcal T$, $\Delta t_\ell=t_\ell-t_{\ell-1}$,
and $u^0=u_0$. Given $u^{\ell-1}$, replacing $\partial_tu$ in
\eqref{eq:evo-bvp} by the backward difference at $t_\ell$ gives the
stage problem:
\begin{equation}\label{eq:time-stepping-stage}
 \mathcal L_xu^\ell+\mathcal N_\ell[u^\ell;\lambda]=F_\ell,
 \qquad
 \mathcal N_\ell[u;\lambda]:=\mathcal N[u;\lambda]+\frac{u}{\Delta t_\ell},
 \qquad
 F_\ell:=F(\cdot,t_\ell)+\frac{u^{\ell-1}}{\Delta t_\ell},
\end{equation}
with the boundary conditions of \eqref{eq:evo-bvp} evaluated at $t_\ell$. With $d_\ell:=(d,t_\ell,\Delta t_\ell,u^{\ell-1})$ the stage datum, fix $\varsigma_\ell\in\{0,1\}$ and write \eqref{eq:time-stepping-stage} as
$(\mathcal L_x+\varsigma_\ell/\Delta t_\ell)u^\ell+q_\ell=0$ with
$q_\ell=\mathcal N[u^\ell;\lambda]+\frac{1-\varsigma_\ell}{\Delta t_\ell}u^\ell-F_\ell$ (for $\varsigma_\ell=0$ the term $u^\ell/\Delta t_\ell$ is part of the
source, for $\varsigma_\ell=1$ it is part of the operator). Since the
added term is of order zero, the identity \eqref{eq:prelim-evo-lagrange}
holds for $\mathcal L_x+\varsigma_\ell/\Delta t_\ell$ with the same
concomitant, and for spatial tests $\varphi_k$ with
$(\mathcal L_x^*+\varsigma_\ell/\Delta t_\ell)\varphi_k=0$, $k=1,\ldots,p$,
the stage global relation is:
\[
 -\int_{\bd}\mathscr B_{\mathcal L}[\boldsymbol\gamma u^\ell,\varphi_k]\,\mathrm ds
 +\int_\Om\varphi_k\,q_\ell\,\mathrm dx=0 ,
\]
i.e., the stationary global relation \eqref{eq:GR} with source $q_\ell$.
Each stage is therefore constructed exactly as in
Section~\ref{sec:elliptic-linear}: the nodal lift \eqref{eq:affine-lift-method},
with the prescribed boundary values evaluated at $t_\ell$, the gate
$\beta$, the matrix $A$ and a stage network
$v_{\theta_\ell}(\cdot;d_\ell)$, $\theta_\ell\in\Theta$. In particular,
the vector of sampled values $\mathbf u_A^\ell(z)\in\R^{n_A}$ has entries:
\[
 (\mathbf u_A^\ell(z))_a=u_{\theta_\ell}(x_a;z,d_\ell)
   =(A_{:,\un}z)_a+(\mathbf a_{\theta_\ell})_a,
\]
with $\mathbf a_{\theta_\ell}$ the vector $\mathbf a_\theta$ of
Section~\ref{sec:elliptic-linear} for the stage data and the stage
network. (For higher-order operators, such as the KdV equation, the nodal
lift is replaced by a Hermite lift of the endpoint traces, with
$A_{:,\un}$ and $\mathbf a_{\theta_\ell}$ defined analogously. See details in 
Appendix~\ref{app:evo-ts-kdv}). Similarly, define the vectors $\mathbf N(z)$, $\mathbf u_A^{\ell-1}$,
$\mathbf F^\ell\in\R^{n_A}$ with entries:
\[
 (\mathbf N(z))_a=\mathcal N[u_{\theta_\ell}(\cdot;z,d_\ell);\lambda](x_a),
 \qquad
 (\mathbf u_A^{\ell-1})_a=u^{\ell-1}(x_a),
 \qquad
 (\mathbf F^\ell)_a=F(x_a,t_\ell).
\]
The sampled stage source is then the vector $\mathbf q_\ell(z)\in\R^{n_A}$, given by:
\[
 \mathbf q_\ell(z)=\mathbf N(z)+\frac{1-\varsigma_\ell}{\Delta t_\ell}\mathbf u_A^\ell(z)
   -\frac{\mathbf u_A^{\ell-1}}{\Delta t_\ell}-\mathbf F^\ell,
\]
and the assembly of Section~\ref{sec:elliptic-linear}, with the stage
tests, gives the finite GR residual of the stage:
\begin{equation}\label{eq:time-stepping-residual}
 r(z;\theta_\ell,d_\ell)=\mathsf Kz-\mathbf b(d_\ell)
   +Q\Bigl(\mathbf N(z)+\frac{1-\varsigma_\ell}{\Delta t_\ell}\mathbf u_A^\ell(z)
          -\frac{\mathbf u_A^{\ell-1}}{\Delta t_\ell}-\mathbf F^\ell\Bigr).
\end{equation}
Since $D_z\mathbf u_A^\ell=A_{:,\un}$ and $\mathbf a_{\theta_\ell}$,
$\mathbf u_A^{\ell-1}$, $\mathbf F^\ell$ do not depend on $z$, the Jacobian
is, for continuously differentiable $\mathbf N$, given by:
\begin{equation}\label{eq:time-stepping-jacobian}
 J_r(z;\theta_\ell,d_\ell)=\mathsf K+Q\Bigl(D_z\mathbf N(z)+\frac{1-\varsigma_\ell}{\Delta t_\ell}A_{:,\un}\Bigr).
\end{equation}
The stage Fokas-PINN is 
$u^\ell:=u_{\theta_\ell}(\cdot;z^*(\theta_\ell,d_\ell),d_\ell)$, where the lift
carries the prescribed coefficients at $t_\ell$ and the gated correction
leaves them unchanged, so the prescribed boundary data are exact at every
step.
\par Naturally, this approach is accompanied with a corresponding accumulation of time-discretization errors. We provide an analysis of these errors in
Appendix~\ref{app:time-stepping-error}.

\section{Numerical experiments}
\label{sec:numerics}

In this section, we implement FokasNets on elliptic problems in two and three dimensions, evolution problems and solution operators, with manufactured solutions, as well as benchmarks from the literature. In this section, $u^*$ denotes the exact solution.

\subsection{Experimental setup}
\label{sec:numerics-setting}

Every FokasNet is compared with two baselines with the same network, collocation, optimizer and budget. In each example, the FokasNet is compared with (a) the standard PINN, which models the solution directly, i.e., $u_\theta=v_\theta$, with all conditions as penalties in
\eqref{eq:prelim-physics-objective}, and (b) a lifted PINN, which uses the lifted approximation \eqref{eq:fpinn-trial-field} of the Fokas-PINN but trains the missing trace coefficients $z$ jointly with $\theta$. For operators, we use the DeepONet architecture with the standard branch-trunk network and a PCA encoding $\mathrm{enc}(d)$ of the datum as branch input, given by:
\begin{equation}
\label{eq:numerics-op-pinn}
 u_\theta(x;d)=\sum_{k=1}^{p_{\mathrm{op}}}b_k(\mathrm{enc}(d))\,t_k(x),
\end{equation}
trained using the Physics-Informed loss \eqref{eq:operator-loss}(PI-DeepONet). For the lifted DeepONet \cite{lu2022comprehensive,sukumar2022exact} and the Fokas Operator Network, we have: 
\begin{equation}
\label{eq:numerics-op-fields}
 u_\theta(x;d)=\ell_0(x;d)+\sum_{j\in\pre}\ell_j(x)\,g_D(x_j^\partial)
   +\sum_{j\in\un}\ell_j(x)\,c_j(d)+\beta(x)\sum_{k=1}^{p_{\mathrm{op}}}b_k(\mathrm{enc}(d))\,t_k(x),
\end{equation}
with:
\begin{equation}
\label{eq:numerics-op-trace}
 c(d)=\begin{cases}
 z_\theta(d) & \text{(lifted PI-DeepONet: $m$ additional branch outputs)},\\[2pt]
 z^*(\theta,d) & \text{(FONet: the completion of each datum)}.
 \end{cases}
\end{equation}

For our experiments, we will consider various different boundary configurations. Specifically, on the unit disk $\mathbb D$ and the unit square $\square=(0,1)^2$ we use the following:
\begin{equation}
\label{eq:numerics-bc}
\begin{aligned}
 &\mathsf M_{\mathbb D}: &&\GD=\{(x,y)\in\partial\mathbb D:\ y\ge0\}, &&\GN=\{(x,y)\in\partial\mathbb D:\ y<0\},\\
 &\mathsf M_{\square}: &&\GD=\{x\in\{0,1\},\ 0\le y\le1\}, &&\GN=\{0<x<1,\ y\in\{0,1\}\},\\
 &\mathsf N_{\mathbb D}: &&\GD=\emptyset, &&\GN=\partial\mathbb D,\\
 &\mathsf A_\alpha: &&\GD=\{(\cos\vartheta,\sin\vartheta):\ 0\le\vartheta<\alpha\}, &&\GN=\partial\mathbb D\setminus\GD,
\end{aligned}
\end{equation}
For the KdV equation, the spatial domain is the interval $\Om=(-1,1)$ and the boundary consists of the two endpoints $x=\pm1$. Since the spatial operator is of third order, each endpoint carries the three traces $u$, $\partial_xu$ and $\partial_x^2u$, and we prescribe $u(-1,t)$, $u(1,t)$ and $\partial_xu(1,t)$, the classical well-posed configuration \cite{himonas2019korteweg}, which we denote by $\mathsf M_{\mathrm{KdV}}$. Each setup therefore changes the amount of information available to the lift, and therefore, the amount of information that the Global Relation solve must complete. 

For all examples, we report the relative $L^2$ errors, given by:
\begin{equation}
\label{eq:numerics-errors}
\begin{aligned}
 E_\Om&=\frac{\|u_\theta-u^*\|_{L^2(\Om)}}{\|u^*\|_{L^2(\Om)}},
 &\qquad
 E_D&=\frac{\|u_\theta-g_D\|_{L^2(\GD)}}{\|g_D\|_{L^2(\GD)}},\\
 E_N&=\frac{\|\dn u_\theta-w_N\|_{L^2(\GN)}}{\|w_N\|_{L^2(\GN)}},
 &\qquad
 E_{\mathrm{miss}}&=\frac{\|u_\theta-u^*\|_{L^2(\GN)}}{\|u^*\|_{L^2(\GN)}}.
\end{aligned}
\end{equation}
These are the interior error, the Dirichlet and Neumann data errors and the completed trace error (the value trace that is not prescribed). Furthermore, we report the pooled error over all test points in the interior and boundary, given by:
\begin{equation}
\label{eq:numerics-closed}
E_{\mathrm{test}}=\biggl(\frac{\sum_i w_i\,|u_\theta(x_i)-u^*(x_i)|^2}{\sum_i w_i\,|u^*(x_i)|^2}\biggr)^{1/2},
\end{equation}
i.e., the relative $\ell^2$ error over the held-out interior and boundary test points, with $w_i$ the quadrature weight of $x_i$ (normalized to mean one, so that interior and boundary points count alike, and $w_i=1$ for the random test points of the
operators). For evolution problems the test points cover $Q_{\mathcal T}$ and $\bd\times(0,\mathcal T)$. Finally, for the solution operator examples, the errors are taken per test datum, followed by the median over the data.

To facilitate the readability of this Section, the problems, solutions, training and testing setups, finite GRs  and lifts are specified in Tables~\ref{tab:app-problems}-\ref{tab:app-lifts} of Appendix~\ref{app:setup}. All results in the tables are medians over the 5 seeds for each method.

\subsection{Elliptic problems}\label{sec:numerics-elliptic}

We consider various PDEs and boundary settings and geometries. The PDEs and results are gathered in Table~\ref{tab:elliptic}. For the manufactured problems, we construct the forcing term such that the PDEs admit a smooth solution. The forms of these solutions (denoted $S$) are given in Table \ref{tab:app-problems} of the Appendix. 

\begin{figure}[!t]
\centering
\includegraphics[width=0.85\textwidth]{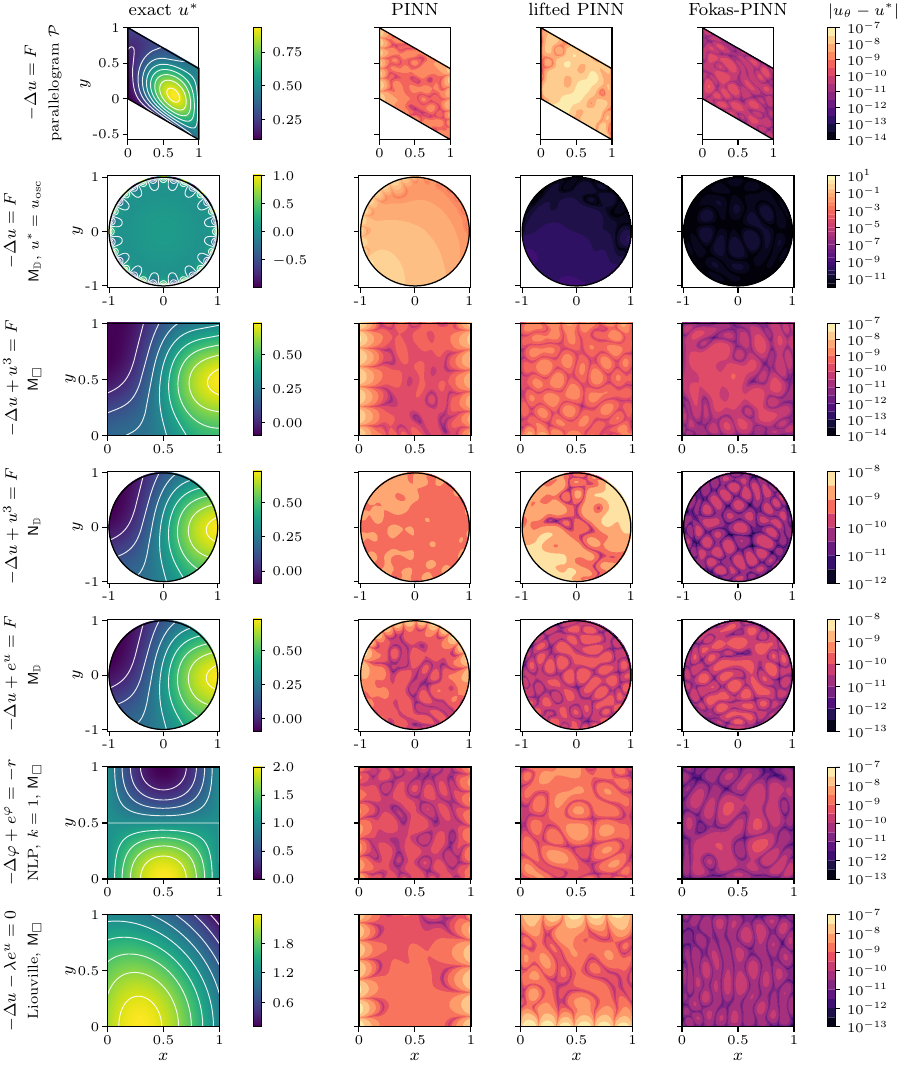}
\caption{True solution and error contours for the 2D elliptic problems in Table \ref{tab:elliptic}.  Each contour shows the median-error seed of each method. }
\label{fig:fields_elliptic}
\end{figure}

The results show that the Fokas-PINN has the smallest total, interior, Dirichlet and completed trace errors in all six problems, achieving up to 2 - 5 orders of magnitude improvement in accuracy. Note that, we include the Poisson problem as a simple verification benchmark, since the source does not depend on the network, so the completion is a fixed linear solve and the trace is exact at the prescribed trace coordinates. Indeed, for this PDE, the Fokas-PINN achieves near machine precision in the completed trace. Furthermore, in the Poisson problem on the disk, the standard PINN is not able to resolve the boundary oscillation of $u_{\rm osc}$, whereas the lifted models are able to achieve high precision. Figure~\ref{fig:fields_elliptic} shows the solution and error contours. 

\begin{figure}[!t]
\centering
\includegraphics[width=0.99\textwidth]{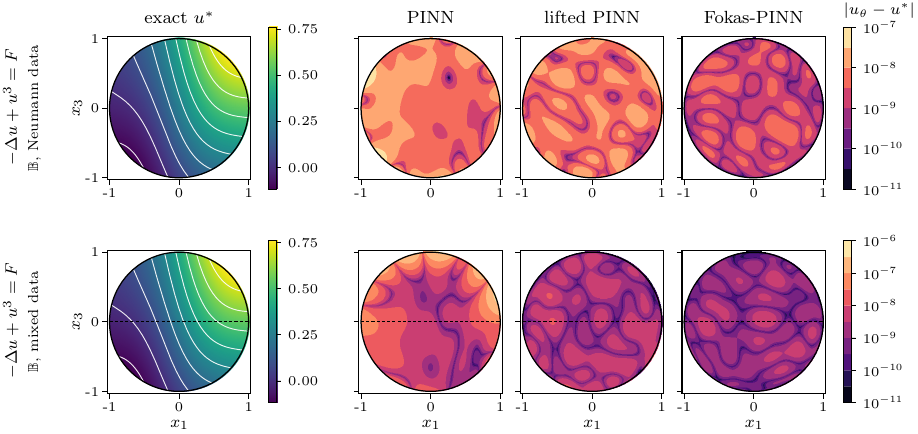}
\caption{True solution and error contours of a slice of the 3D ball example.}
\label{fig:fields_ball}
\end{figure}

Furthermore, we solve $-\Delta u+u^3=F$ on the unit ball $\mathbb B\subset\R^3$, with $F$ and the boundary data taken from the manufactured solution $S_3$ \eqref{eq:app-S3}, for two sets of boundary conditions. In the first, the flux $\partial_n u=w_N$ is
prescribed on the whole sphere (corresponding to the 2D $\mathsf N_{\mathbb D}$ domain used above). In the second, the Dirichlet data $u=g_D$ is prescribed on the upper hemisphere $\{x_3\ge0\}$ and
the flux data on the lower one. The completion of the missing boundary values is described in Appendix~\ref{app:stat-ball} and Tables~\ref{tab:constructions} and~\ref{tab:app-lifts} list the
corresponding lifts and GR constructions. As shown in Table~\ref{tab:elliptic}, the Fokas-PINN is the most accurate method in both cases. The total test error is reduced by $\times 3.7$ and $\times 2.2$  relative to the lifted PINN, and the completed trace by $\times 240$ and $\times 4$ in the Neumann and mixed boundary condition cases, respectively. The solution and error contours are given in Figure~\ref{fig:fields_ball}.

\begin{table}[!t]\centering\footnotesize\setlength{\tabcolsep}{4pt}
\begin{adjustbox}{max width=\linewidth}
\begin{tabular}{@{}llccccc@{}}
\toprule
problem & method & $E_{\mathrm{test}}$ & interior & Dirichlet data & Neumann data & completed trace \\
\midrule
\multicolumn{7}{@{}l}{\textit{Manufactured solutions}}\\[1pt]
$-\Delta u=F$ & PINN & $1.0{\times}10^{-8}$ & $7.1{\times}10^{-9}$ & $4.6{\times}10^{-8}$ & $9.8{\times}10^{-9}$ & $2.1{\times}10^{-8}$ \\
{\scriptsize parallelogram $\mathcal P$} & lifted PINN & $3.8{\times}10^{-8}$ & $3.7{\times}10^{-8}$ & $\mathbf{2.2{\times}10^{-16}}$ & $1.0{\times}10^{-7}$ & $7.3{\times}10^{-8}$ \\
{\scriptsize $u^*=T+\tfrac12\beta W$} & Fokas--PINN & $\mathbf{1.4{\times}10^{-10}}$ & $\mathbf{1.4{\times}10^{-10}}$ & $\mathbf{2.2{\times}10^{-16}}$ & $\mathbf{7.6{\times}10^{-10}}$ & $\mathbf{6.9{\times}10^{-12}}$ \\
\addlinespace[3pt]
$-\Delta u=F$ & PINN & $9.6{\times}10^{-1}$ & $1.2{\times}10^{0}$ & $6.7{\times}10^{-1}$ & $2.2{\times}10^{-1}$ & $6.8{\times}10^{-1}$ \\
{\scriptsize $\mathsf{M}_{\mathbb{D}}$, $u^*=u_{\rm osc}$} & lifted PINN & $6.1{\times}10^{-10}$ & $8.0{\times}10^{-10}$ & $6.4{\times}10^{-12}$ & $1.0{\times}10^{-10}$ & $3.8{\times}10^{-10}$ \\
 & Fokas--PINN & $\mathbf{3.0{\times}10^{-11}}$ & $\mathbf{4.0{\times}10^{-11}}$ & $\mathbf{1.3{\times}10^{-15}}$ & $\mathbf{3.7{\times}10^{-12}}$ & $\mathbf{5.1{\times}10^{-15}}$ \\
\addlinespace[3pt]
$-\Delta u+u^3=F$ & PINN & $1.6{\times}10^{-8}$ & $1.1{\times}10^{-8}$ & $4.2{\times}10^{-8}$ & $2.6{\times}10^{-8}$ & $2.3{\times}10^{-8}$ \\
{\scriptsize $\mathsf{M}_{\square}$, $u^*=S$} & lifted PINN & $2.8{\times}10^{-9}$ & $2.5{\times}10^{-9}$ & $\mathbf{3.1{\times}10^{-16}}$ & $1.1{\times}10^{-7}$ & $7.9{\times}10^{-9}$ \\
 & Fokas--PINN & $\mathbf{2.4{\times}10^{-10}}$ & $\mathbf{2.5{\times}10^{-10}}$ & $\mathbf{3.1{\times}10^{-16}}$ & $\mathbf{3.3{\times}10^{-9}}$ & $\mathbf{2.8{\times}10^{-12}}$ \\
\addlinespace[3pt]
$-\Delta u+u^3=F$ & PINN & $2.6{\times}10^{-9}$ & $2.6{\times}10^{-9}$ & -- & $5.2{\times}10^{-8}$ & $3.1{\times}10^{-9}$ \\
{\scriptsize $\mathsf{N}_{\mathbb{D}}$, $u^*=S$} & lifted PINN & $6.2{\times}10^{-9}$ & $6.0{\times}10^{-9}$ & -- & $9.6{\times}10^{-9}$ & $6.7{\times}10^{-9}$ \\
 & Fokas--PINN & $\mathbf{2.2{\times}10^{-10}}$ & $\mathbf{2.3{\times}10^{-10}}$ & -- & $\mathbf{9.5{\times}10^{-9}}$ & $\mathbf{8.6{\times}10^{-11}}$ \\
\addlinespace[3pt]
$-\Delta u+e^{u}=F$ & PINN & $2.4{\times}10^{-9}$ & $1.9{\times}10^{-9}$ & $1.0{\times}10^{-8}$ & $7.4{\times}10^{-8}$ & $1.5{\times}10^{-9}$ \\
{\scriptsize $\mathsf{M}_{\mathbb{D}}$, $u^*=S$} & lifted PINN & $2.5{\times}10^{-10}$ & $2.6{\times}10^{-10}$ & $5.9{\times}10^{-12}$ & $\mathbf{3.7{\times}10^{-9}}$ & $2.4{\times}10^{-10}$ \\
 & Fokas--PINN & $\mathbf{1.9{\times}10^{-10}}$ & $\mathbf{2.0{\times}10^{-10}}$ & $\mathbf{3.1{\times}10^{-13}}$ & $1.4{\times}10^{-8}$ & $\mathbf{1.9{\times}10^{-11}}$ \\
\midrule
\multicolumn{7}{@{}l}{\textit{Benchmarks from the literature}}\\[1pt]
$-\Delta u+e^u=-r$ & PINN & $1.1{\times}10^{-9}$ & $7.7{\times}10^{-10}$ & $3.7{\times}10^{-9}$ & $7.4{\times}10^{-9}$ & $1.2{\times}10^{-9}$ \\
{\scriptsize nonlinear Poisson, $\mathsf{M}_{\square}$} & lifted PINN & $1.8{\times}10^{-9}$ & $1.8{\times}10^{-9}$ & $\mathbf{2.0{\times}10^{-16}}$ & $4.4{\times}10^{-8}$ & $1.9{\times}10^{-9}$ \\
{\scriptsize \cite{urban2025unveiling}} & Fokas--PINN & $\mathbf{1.1{\times}10^{-10}}$ & $\mathbf{1.2{\times}10^{-10}}$ & $\mathbf{2.0{\times}10^{-16}}$ & $\mathbf{2.7{\times}10^{-9}}$ & $\mathbf{1.2{\times}10^{-11}}$ \\
\addlinespace[3pt]
$-\Delta u-\lambda e^{u}=0$ & PINN & $5.1{\times}10^{-9}$ & $3.7{\times}10^{-9}$ & $1.8{\times}10^{-8}$ & $1.2{\times}10^{-8}$ & $6.2{\times}10^{-9}$ \\
{\scriptsize Liouville, $\mathsf{M}_{\square}$} & lifted PINN & $7.3{\times}10^{-9}$ & $6.8{\times}10^{-9}$ & $\mathbf{3.5{\times}10^{-16}}$ & $2.9{\times}10^{-7}$ & $1.5{\times}10^{-8}$ \\
{\scriptsize \cite{liouville1853equation,gel1963some}} & Fokas--PINN & $\mathbf{3.8{\times}10^{-11}}$ & $\mathbf{4.0{\times}10^{-11}}$ & $\mathbf{3.5{\times}10^{-16}}$ & $\mathbf{1.5{\times}10^{-9}}$ & $\mathbf{1.0{\times}10^{-11}}$ \\
\midrule
\multicolumn{7}{@{}l}{\textit{Three-dimensional problems}}\\[1pt]
$-\Delta u+u^3=F$ & PINN & $3.7{\times}10^{-8}$ & $3.5{\times}10^{-8}$ & -- & $1.9{\times}10^{-6}$ & $7.6{\times}10^{-8}$ \\
{\scriptsize unit ball $\mathbb B$, Neumann data} & lifted PINN & $2.6{\times}10^{-8}$ & $2.6{\times}10^{-8}$ & -- & $\mathbf{1.2{\times}10^{-7}}$ & $2.7{\times}10^{-8}$ \\
{\scriptsize $u^*=S_3$} & Fokas--PINN & $\mathbf{6.9{\times}10^{-9}}$ & $\mathbf{7.1{\times}10^{-9}}$ & -- & $2.5{\times}10^{-7}$ & $\mathbf{1.1{\times}10^{-10}}$ \\
\addlinespace[3pt]
$-\Delta u+u^3=F$ & PINN & $1.6{\times}10^{-7}$ & $1.5{\times}10^{-7}$ & $4.2{\times}10^{-7}$ & $1.6{\times}10^{-6}$ & $9.6{\times}10^{-8}$ \\
{\scriptsize unit ball $\mathbb B$, mixed data} & lifted PINN & $1.1{\times}10^{-8}$ & $1.2{\times}10^{-8}$ & $\mathbf{2.5{\times}10^{-11}}$ & $2.9{\times}10^{-7}$ & $1.7{\times}10^{-8}$ \\
{\scriptsize $u^*=S_3$} & Fokas--PINN & $\mathbf{5.2{\times}10^{-9}}$ & $\mathbf{5.3{\times}10^{-9}}$ & $\mathbf{2.5{\times}10^{-11}}$ & $\mathbf{2.2{\times}10^{-7}}$ & $\mathbf{4.0{\times}10^{-9}}$ \\
\bottomrule
\end{tabular}
\end{adjustbox}
\caption{Error results on the test set for the elliptic problems. Best errors are in {\bf bold}.}
\label{tab:elliptic}
\end{table}

Finally, to show the power of the GR completion, we perform a simple experiment on the disc, whereby we vary $\alpha$ on $A_{\alpha}$, which changes the arc on which the Dirichlet data are known. The results are given in Table~\ref{tab:app-sweep}, which show that the Fokas-PINN is robust to the information provided on the boundary. Specifically, as $\alpha$ shrinks to zero, the Fokas-PINN retains the same total error (reaching an order of magnitude higher accuracy compared to the PINN benchmarks) and retains or improves the boundary value errors, achieving accuracies up to three orders of magnitude higher on the Dirichlet data and two orders on the completed trace.

\begin{table}[!h]\centering\footnotesize\setlength{\tabcolsep}{4pt}
\begin{adjustbox}{max width=\linewidth}
\begin{tabular}{@{}llccccc@{}}
\toprule
problem & method & $E_{\mathrm{test}}$ & interior & Dirichlet data & Neumann data & completed trace \\
\midrule
$\mathsf{A}_0$ & PINN & $2.6{\times}10^{-9}$ & $2.6{\times}10^{-9}$ & -- & $5.2{\times}10^{-8}$ & $3.1{\times}10^{-9}$ \\
{\scriptsize $=\mathsf{N}_{\mathbb{D}}$} & lifted PINN & $6.2{\times}10^{-9}$ & $6.0{\times}10^{-9}$ & -- & $9.6{\times}10^{-9}$ & $6.7{\times}10^{-9}$ \\
 & Fokas--PINN & $\mathbf{2.2{\times}10^{-10}}$ & $\mathbf{2.3{\times}10^{-10}}$ & -- & $\mathbf{9.5{\times}10^{-9}}$ & $\mathbf{8.6{\times}10^{-11}}$ \\
\addlinespace[3pt]
$\mathsf{A}_{\pi/8}$ & PINN & $3.0{\times}10^{-9}$ & $2.9{\times}10^{-9}$ & $3.1{\times}10^{-9}$ & $5.6{\times}10^{-8}$ & $3.8{\times}10^{-9}$ \\
 & lifted PINN & $2.0{\times}10^{-9}$ & $2.0{\times}10^{-9}$ & $1.9{\times}10^{-11}$ & $1.3{\times}10^{-8}$ & $2.2{\times}10^{-9}$ \\
 & Fokas--PINN & $\mathbf{1.9{\times}10^{-10}}$ & $\mathbf{1.9{\times}10^{-10}}$ & $\mathbf{7.2{\times}10^{-14}}$ & $\mathbf{4.6{\times}10^{-9}}$ & $\mathbf{5.3{\times}10^{-12}}$ \\
\addlinespace[3pt]
$\mathsf{A}_{\pi/2}$ & PINN & $3.0{\times}10^{-9}$ & $2.4{\times}10^{-9}$ & $7.6{\times}10^{-9}$ & $5.3{\times}10^{-8}$ & $2.9{\times}10^{-9}$ \\
 & lifted PINN & $5.5{\times}10^{-9}$ & $5.4{\times}10^{-9}$ & $7.6{\times}10^{-11}$ & $4.2{\times}10^{-8}$ & $8.5{\times}10^{-9}$ \\
 & Fokas--PINN & $\mathbf{2.2{\times}10^{-10}}$ & $\mathbf{2.3{\times}10^{-10}}$ & $\mathbf{1.5{\times}10^{-13}}$ & $\mathbf{7.8{\times}10^{-9}}$ & $\mathbf{2.1{\times}10^{-11}}$ \\
\addlinespace[3pt]
$\mathsf{A}_{\pi}$ & PINN & $3.6{\times}10^{-9}$ & $2.8{\times}10^{-9}$ & $1.5{\times}10^{-8}$ & $1.1{\times}10^{-7}$ & $2.3{\times}10^{-9}$ \\
{\scriptsize $=\mathsf{M}_{\mathbb{D}}$} & lifted PINN & $3.4{\times}10^{-10}$ & $3.5{\times}10^{-10}$ & $1.6{\times}10^{-11}$ & $\mathbf{9.3{\times}10^{-9}}$ & $4.3{\times}10^{-10}$ \\
 & Fokas--PINN & $\mathbf{2.2{\times}10^{-10}}$ & $\mathbf{2.2{\times}10^{-10}}$ & $\mathbf{2.7{\times}10^{-13}}$ & $2.0{\times}10^{-8}$ & $\mathbf{2.4{\times}10^{-11}}$ \\
\bottomrule
\end{tabular}
\end{adjustbox}
\caption{Error results for $-\Delta u+u^3=F$ in the unit disk with $u^*=S$, $u=g_D$ on the arc $\{\vartheta\in[0,\alpha)\}$ and $\partial_nu=w_N$ on the rest of the circle. Best errors are in {\bf bold}.}
\label{tab:app-sweep}
\end{table}

\subsection{Evolution problems} \label{sec:numerics-evolution}

\begin{figure}[!t]
\centering
\includegraphics[width=\textwidth]{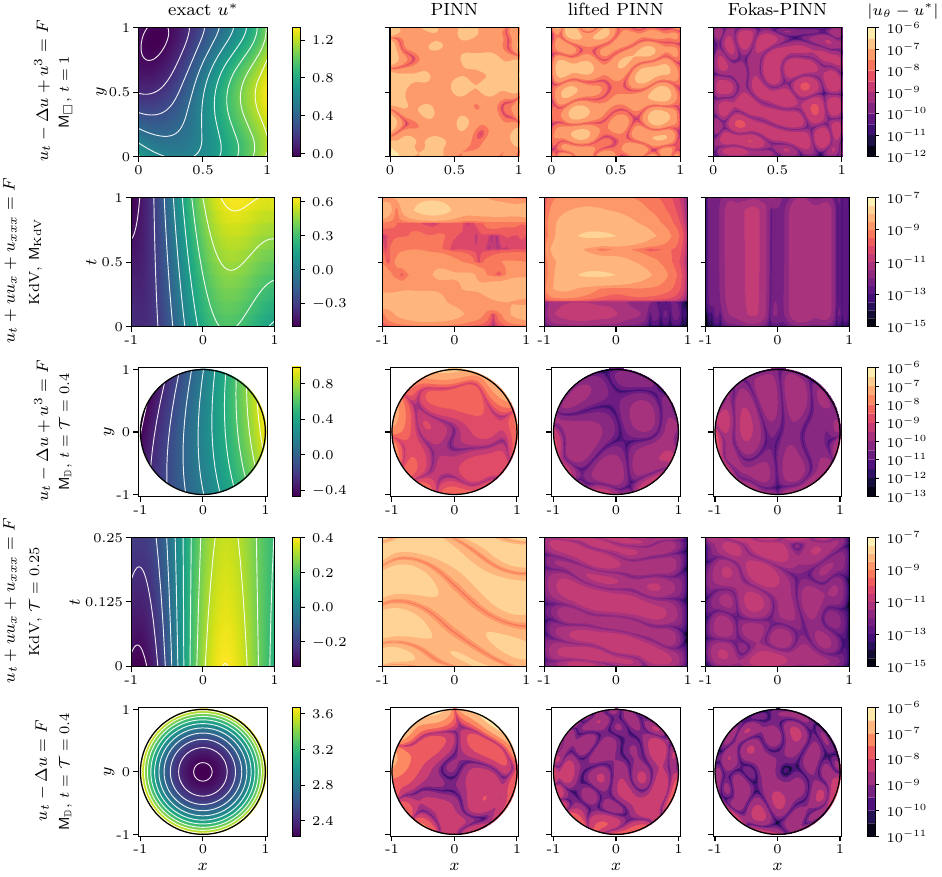}
\caption{True solution and error contours for the evolution problems in Table~\ref{tab:evolution}. Rows 1-2: time stepping, heat equation on the square at $t=\mathcal T=1$ and
KdV equation in the $(x,t)$ plane. Rows 3-4: space-time global relation, heat
equation on the disk at $t=\mathcal T=0.4$ and KdV equation in the $(x,t)$ plane. }
\label{fig:fields_evolution}
\end{figure}

Next, we consider evolution PDEs. As for the elliptic examples, we use problems with manufactured solutions as well as from the literature. The form of the manufactured solutions we use are again given in Table \ref{tab:app-problems} of the Appendix. Table~\ref{tab:evolution} contains the results, using both the time-stepping construction of Section~\ref{prop:evolution-time-stepping}, as well as the space-time relation of Section~\ref{sec:evo-completion-P}. With time stepping, where backward Euler is exact in time for the chosen solutions, the Fokas-PINN reduces the total error by one to two orders of magnitude, achieving an error of $3.3\times 10^{-12}$ on the forced KdV \cite{himonas2019korteweg}. Furthermore, the time-stepping Fokas-PINN achieves two to three orders of magnitude reduction in the completed trace.  With the space-time GR, the Fokas-PINN improves the total accuracy in two of the three examples (tying with the lifted PINN in one case). In all three examples the completed trace error improves, reaching nearly one order of magnitude more accurate results. The solution and error contours are given in Figure \ref{fig:fields_evolution}.

\begin{table}[!t]\centering\footnotesize\setlength{\tabcolsep}{4pt}
\begin{adjustbox}{max width=\linewidth}
\begin{tabular}{@{}llccccc@{}}
\toprule
problem & method & $E_{\mathrm{test}}$ & interior & Dirichlet data & Neumann data & completed trace \\
\midrule
\multicolumn{7}{@{}l}{\textit{Time stepping (10 backward-Euler levels, $\Delta t=0.1$): errors at $t=\mathcal T=1$}}\\[1pt]
$u_t-\Delta u+u^3=F$ & PINN & $1.5{\times}10^{-7}$ & $1.4{\times}10^{-7}$ & $3.1{\times}10^{-7}$ & $8.1{\times}10^{-7}$ & $2.2{\times}10^{-7}$ \\
{\scriptsize $\mathsf{M}_{\square}$} & lifted PINN & $9.3{\times}10^{-8}$ & $9.5{\times}10^{-8}$ & $\mathbf{2.3{\times}10^{-16}}$ & $2.8{\times}10^{-6}$ & $1.7{\times}10^{-7}$ \\
 & Fokas--PINN & $\mathbf{1.7{\times}10^{-9}}$ & $\mathbf{1.9{\times}10^{-9}}$ & $\mathbf{2.3{\times}10^{-16}}$ & $\mathbf{3.5{\times}10^{-8}}$ & $\mathbf{1.6{\times}10^{-10}}$ \\
\addlinespace[3pt]
$u_t+uu_x+u_{xxx}=F$ & PINN & $1.0{\times}10^{-8}$ & $1.0{\times}10^{-8}$ & $2.4{\times}10^{-9}$ & $1.4{\times}10^{-8}$ & $1.9{\times}10^{-8}$ \\
{\scriptsize KdV, $\mathsf M_{\mathrm{KdV}}$} & lifted PINN & $2.5{\times}10^{-9}$ & $2.5{\times}10^{-9}$ & $\mathbf{2.8{\times}10^{-16}}$ & $\mathbf{5.3{\times}10^{-16}}$ & $6.5{\times}10^{-9}$ \\
 & Fokas--PINN & $\mathbf{3.3{\times}10^{-12}}$ & $\mathbf{3.4{\times}10^{-12}}$ & $\mathbf{2.8{\times}10^{-16}}$ & $\mathbf{5.3{\times}10^{-16}}$ & $\mathbf{1.3{\times}10^{-11}}$ \\
\midrule
\multicolumn{7}{@{}l}{\textit{Space--time Global Relation: errors over $Q_{\mathcal T}$}}\\[1pt]
$u_t-\Delta u+u^3=F$ & PINN & $1.6{\times}10^{-8}$ & $1.5{\times}10^{-8}$ & $3.2{\times}10^{-8}$ & $5.1{\times}10^{-8}$ & $8.4{\times}10^{-9}$ \\
{\scriptsize $\mathsf{M}_{\mathbb{D}}$, $\mathcal T=0.4$} & lifted PINN & $\mathbf{3.7{\times}10^{-10}}$ & $\mathbf{3.6{\times}10^{-10}}$ & $4.8{\times}10^{-11}$ & $3.7{\times}10^{-9}$ & $6.7{\times}10^{-10}$ \\
 & Fokas--PINN & $\mathbf{3.7{\times}10^{-10}}$ & $3.8{\times}10^{-10}$ & $\mathbf{3.5{\times}10^{-11}}$ & $\mathbf{3.2{\times}10^{-9}}$ & $\mathbf{3.5{\times}10^{-10}}$ \\
\addlinespace[3pt]
$u_t+uu_x+u_{xxx}=F$ & PINN & $4.3{\times}10^{-8}$ & $4.2{\times}10^{-8}$ & $8.2{\times}10^{-8}$ & $6.9{\times}10^{-8}$ & $1.7{\times}10^{-7}$ \\
{\scriptsize KdV, $\mathsf M_{\mathrm{KdV}}$, $\mathcal T=0.25$} & lifted PINN & $3.9{\times}10^{-11}$ & $4.0{\times}10^{-11}$ & $\mathbf{5.0{\times}10^{-17}}$ & $\mathbf{2.1{\times}10^{-16}}$ & $1.1{\times}10^{-10}$ \\
 & Fokas--PINN & $\mathbf{2.6{\times}10^{-11}}$ & $\mathbf{2.7{\times}10^{-11}}$ & $\mathbf{5.0{\times}10^{-17}}$ & $\mathbf{2.1{\times}10^{-16}}$ & $\mathbf{6.2{\times}10^{-11}}$ \\
\addlinespace[3pt]
$u_t-\Delta u=F$ & PINN & $1.5{\times}10^{-8}$ & $1.4{\times}10^{-8}$ & $3.4{\times}10^{-8}$ & $4.0{\times}10^{-8}$ & $7.8{\times}10^{-9}$ \\
{\scriptsize PINNacle heat, $\mathsf{M}_{\mathbb{D}}$, $\mathcal T=0.4$} & lifted PINN & $1.8{\times}10^{-9}$ & $1.7{\times}10^{-9}$ & $2.2{\times}10^{-10}$ & $2.3{\times}10^{-8}$ & $3.7{\times}10^{-9}$ \\
{\scriptsize \cite{hao2024pinnacle}} & Fokas--PINN & $\mathbf{6.3{\times}10^{-10}}$ & $\mathbf{6.4{\times}10^{-10}}$ & $\mathbf{3.2{\times}10^{-11}}$ & $\mathbf{6.2{\times}10^{-9}}$ & $\mathbf{6.2{\times}10^{-10}}$ \\
\bottomrule
\end{tabular}
\end{adjustbox}
\caption{Error results for the evolution problems. For the space-time implementation, the errors are computed over the space-time domain $Q_{\mathcal T}$ and its boundary $\partial\Omega\times(0,\mathcal T)$. For the time-stepping implementation, they are computed at the final time $t=\mathcal T$. Best errors are in {\bf bold}.}
\label{tab:evolution}
\end{table}

\subsection{Solution operators}\label{sec:numerics-operators}

We now consider solution operator examples for both elliptic and evolution PDEs. As we will see, in the operator setting, the improvements provided by FONets are particularly prevalent. The results are given in Table~\ref{tab:operators}, showing that the FONet is the most accurate method in every case, across PDE families and GR implementations (elliptic, time-stepping and space-time). More specifically, it reduces the total error by $\times 10 - 40$ on the manufactured stationary families, $\times 3 - 17$ on the benchmarks and $\times 9$ on the
space-time operator, with corresponding trace error improvements between 2-4 orders of magnitude in most examples. 
\par We highlight the cubic reaction example with the Neumann data, $\mathsf N_{\mathbb D}$, where the improvements are four orders of magnitude both in the total and the completed trace accuracy. Here, the baselines are not able to capture the underlying dynamics due to the well-known ill-conditioning of pure-Neumann semilinear problems for residual-based methods, whereas the   GR component of the FONet is precisely the mechanism that provides the necessary boundary data. Furthermore, we note that, for the linear Poisson operator, the accuracy in the completed trace is near machine-precision, since the GR completion is a fixed linear map of each datum. An interesting observation is also that the FONets produce lower errors on the boundary value data than in the interior for every example.    

\begin{figure}[!h]
\centering
\includegraphics[width=\textwidth]{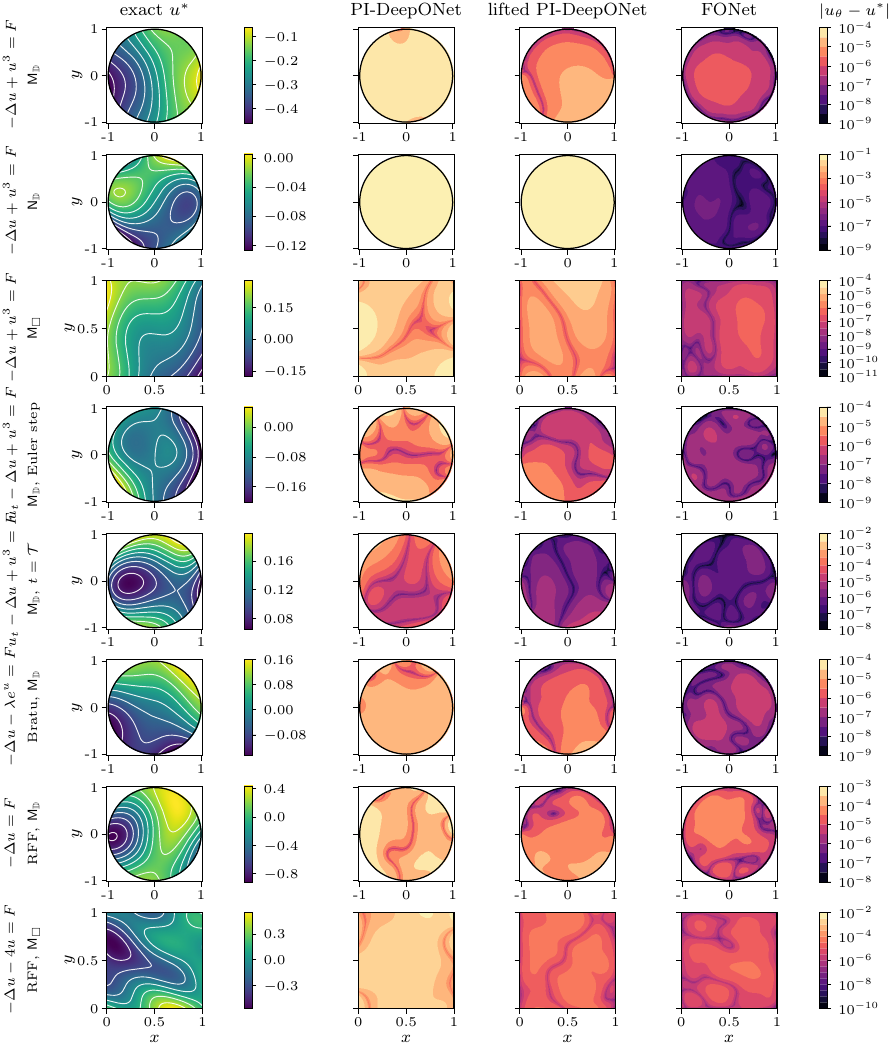}
\caption{Solution operators on one unseen test datum (the median-error datum of the FONet at its median seed; all methods at that seed and datum): exact solution $u^*$ and pointwise error. }
\label{fig:fields_operator}
\end{figure}

In Figure~\ref{fig:fields_operator} we plot the contours for the solution operators applied to an unseen test datum.

\begin{table}[!h]\centering\footnotesize\setlength{\tabcolsep}{4pt}
\begin{adjustbox}{max width=\linewidth}
\begin{tabular}{@{}llccccc@{}}
\toprule
problem & method & $E_{\mathrm{test}}$ & interior & Dirichlet data & Neumann data & completed trace \\
\midrule
\multicolumn{7}{@{}l}{\textit{Manufactured data}}\\[1pt]
$-\Delta u+u^3=F$ & PI-DeepONet & $1.1{\times}10^{-4}$ & $1.1{\times}10^{-4}$ & $1.4{\times}10^{-4}$ & $2.2{\times}10^{-4}$ & $1.2{\times}10^{-4}$ \\
{\scriptsize $\mathsf{M}_{\mathbb{D}}$} & lifted PI-DeepONet & $3.3{\times}10^{-5}$ & $3.4{\times}10^{-5}$ & $1.2{\times}10^{-6}$ & $1.7{\times}10^{-4}$ & $5.6{\times}10^{-5}$ \\
 & FONet & $\mathbf{3.5{\times}10^{-6}}$ & $\mathbf{3.6{\times}10^{-6}}$ & $\mathbf{1.0{\times}10^{-9}}$ & $\mathbf{1.8{\times}10^{-5}}$ & $\mathbf{8.8{\times}10^{-8}}$ \\
\addlinespace[3pt]
$-\Delta u+\lambda u^3=F$ & PI-DeepONet & $1.5{\times}10^{-4}$ & $1.5{\times}10^{-4}$ & $2.0{\times}10^{-4}$ & $2.7{\times}10^{-4}$ & $1.6{\times}10^{-4}$ \\
{\scriptsize $\mathsf{M}_{\mathbb{D}}$, $\lambda\in[0.2,2]$} & lifted PI-DeepONet & $6.4{\times}10^{-5}$ & $6.3{\times}10^{-5}$ & $4.0{\times}10^{-6}$ & $3.2{\times}10^{-4}$ & $1.0{\times}10^{-4}$ \\
 & FONet & $\mathbf{1.6{\times}10^{-6}}$ & $\mathbf{1.6{\times}10^{-6}}$ & $\mathbf{3.7{\times}10^{-10}}$ & $\mathbf{9.9{\times}10^{-6}}$ & $\mathbf{2.8{\times}10^{-8}}$ \\
\addlinespace[3pt]
$-\Delta u+u^3=F$ & PI-DeepONet & $1.4{\times}10^{-1}$ & $1.4{\times}10^{-1}$ & -- & $7.6{\times}10^{-3}$ & $1.3{\times}10^{-1}$ \\
{\scriptsize $\mathsf{N}_{\mathbb{D}}$} & lifted PI-DeepONet & $2.6{\times}10^{-2}$ & $2.7{\times}10^{-2}$ & -- & $1.1{\times}10^{-3}$ & $2.4{\times}10^{-2}$ \\
 & FONet & $\mathbf{2.7{\times}10^{-6}}$ & $\mathbf{2.7{\times}10^{-6}}$ & -- & $\mathbf{2.0{\times}10^{-5}}$ & $\mathbf{2.1{\times}10^{-6}}$ \\
\addlinespace[3pt]
$-\Delta u+u^3=F$ & PI-DeepONet & $2.2{\times}10^{-4}$ & $2.1{\times}10^{-4}$ & $3.0{\times}10^{-4}$ & $3.1{\times}10^{-4}$ & $2.5{\times}10^{-4}$ \\
{\scriptsize $\mathsf{M}_{\square}$} & lifted PI-DeepONet & $4.4{\times}10^{-5}$ & $4.4{\times}10^{-5}$ & $\mathbf{4.2{\times}10^{-16}}$ & $1.5{\times}10^{-4}$ & $6.3{\times}10^{-5}$ \\
 & FONet & $\mathbf{1.5{\times}10^{-6}}$ & $\mathbf{1.6{\times}10^{-6}}$ & $\mathbf{4.2{\times}10^{-16}}$ & $\mathbf{6.2{\times}10^{-6}}$ & $\mathbf{4.8{\times}10^{-9}}$ \\
\midrule
\multicolumn{7}{@{}l}{\textit{Benchmark operators}}\\[1pt]
$-\Delta u-\lambda e^{u}=F$ & PI-DeepONet & $1.5{\times}10^{-4}$ & $1.5{\times}10^{-4}$ & $1.8{\times}10^{-4}$ & $2.4{\times}10^{-4}$ & $1.5{\times}10^{-4}$ \\
{\scriptsize Bratu, $\mathsf{M}_{\mathbb{D}}$, $\lambda\in[0.25,0.75]$} & lifted PI-DeepONet & $6.4{\times}10^{-5}$ & $6.4{\times}10^{-5}$ & $3.5{\times}10^{-6}$ & $2.9{\times}10^{-4}$ & $1.0{\times}10^{-4}$ \\
 & FONet & $\mathbf{3.7{\times}10^{-6}}$ & $\mathbf{3.9{\times}10^{-6}}$ & $\mathbf{1.1{\times}10^{-8}}$ & $\mathbf{1.6{\times}10^{-5}}$ & $\mathbf{1.0{\times}10^{-6}}$ \\
\addlinespace[3pt]
$-\Delta u=F$ & PI-DeepONet & $1.4{\times}10^{-3}$ & $1.3{\times}10^{-3}$ & $2.6{\times}10^{-3}$ & $1.5{\times}10^{-3}$ & $1.6{\times}10^{-3}$ \\
{\scriptsize RFF, $\mathsf{M}_{\mathbb{D}}$} & lifted PI-DeepONet & $2.2{\times}10^{-4}$ & $2.2{\times}10^{-4}$ & $7.2{\times}10^{-6}$ & $2.9{\times}10^{-4}$ & $2.9{\times}10^{-4}$ \\
 & FONet & $\mathbf{8.2{\times}10^{-5}}$ & $\mathbf{8.5{\times}10^{-5}}$ & $\mathbf{1.2{\times}10^{-15}}$ & $\mathbf{1.9{\times}10^{-4}}$ & $\mathbf{3.1{\times}10^{-15}}$ \\
\addlinespace[3pt]
$-\Delta u-4u=F$ & PI-DeepONet & $3.3{\times}10^{-3}$ & $3.2{\times}10^{-3}$ & $5.5{\times}10^{-3}$ & $2.2{\times}10^{-3}$ & $4.3{\times}10^{-3}$ \\
{\scriptsize RFF, $\mathsf{M}_{\square}$} & lifted PI-DeepONet & $2.2{\times}10^{-4}$ & $2.2{\times}10^{-4}$ & $\mathbf{7.5{\times}10^{-16}}$ & $1.8{\times}10^{-4}$ & $2.6{\times}10^{-4}$ \\
 & FONet & $\mathbf{5.8{\times}10^{-5}}$ & $\mathbf{5.9{\times}10^{-5}}$ & $\mathbf{7.5{\times}10^{-16}}$ & $\mathbf{7.4{\times}10^{-5}}$ & $\mathbf{1.6{\times}10^{-5}}$ \\
\midrule
\multicolumn{7}{@{}l}{\textit{Evolution operators}}\\[1pt]
$-\Delta u+u/\Delta t+u^3=\tilde F$ & PI-DeepONet & $8.6{\times}10^{-5}$ & $7.9{\times}10^{-5}$ & $1.8{\times}10^{-4}$ & $3.1{\times}10^{-4}$ & $8.1{\times}10^{-5}$ \\
{\scriptsize $\mathsf{M}_{\mathbb{D}}$, one level} & lifted PI-DeepONet & $2.1{\times}10^{-5}$ & $2.1{\times}10^{-5}$ & $1.3{\times}10^{-6}$ & $1.8{\times}10^{-4}$ & $3.8{\times}10^{-5}$ \\
 & FONet & $\mathbf{1.5{\times}10^{-6}}$ & $\mathbf{1.5{\times}10^{-6}}$ & $\mathbf{2.9{\times}10^{-8}}$ & $\mathbf{1.4{\times}10^{-5}}$ & $\mathbf{1.4{\times}10^{-6}}$ \\
\addlinespace[3pt]
$-\Delta u+u/\Delta t+u^3=\tilde F$ & PI-DeepONet & $1.1{\times}10^{-4}$ & $1.0{\times}10^{-4}$ & $2.0{\times}10^{-4}$ & $3.5{\times}10^{-4}$ & $1.1{\times}10^{-4}$ \\
{\scriptsize rollout over 10 levels, $t=1$} & lifted PI-DeepONet & $3.0{\times}10^{-5}$ & $3.0{\times}10^{-5}$ & $1.6{\times}10^{-6}$ & $2.1{\times}10^{-4}$ & $5.1{\times}10^{-5}$ \\
 & FONet & $\mathbf{2.2{\times}10^{-6}}$ & $\mathbf{2.3{\times}10^{-6}}$ & $\mathbf{1.0{\times}10^{-8}}$ & $\mathbf{1.5{\times}10^{-5}}$ & $\mathbf{8.3{\times}10^{-7}}$ \\
\addlinespace[3pt]
$u_t-\Delta u+u^3=F$ & PI-DeepONet & $1.5{\times}10^{-3}$ & $1.3{\times}10^{-3}$ & $2.1{\times}10^{-3}$ & $3.8{\times}10^{-3}$ & $1.1{\times}10^{-3}$ \\
{\scriptsize $\mathsf{M}_{\mathbb{D}}$, space--time, $\mathcal T=0.4$} & lifted PI-DeepONet & $2.6{\times}10^{-5}$ & $2.5{\times}10^{-5}$ & $2.5{\times}10^{-6}$ & $1.4{\times}10^{-4}$ & $3.9{\times}10^{-5}$ \\
 & FONet & $\mathbf{2.9{\times}10^{-6}}$ & $\mathbf{3.3{\times}10^{-6}}$ & $\mathbf{2.9{\times}10^{-7}}$ & $\mathbf{3.0{\times}10^{-5}}$ & $\mathbf{3.8{\times}10^{-6}}$ \\
\bottomrule
\end{tabular}
\end{adjustbox}
\caption{Error metrics for the solution operators. For the
level operator, ``one level'' is a single time step from the exact previous state, and ``rollout'' is ten consecutive steps from $u_0$, each output used as the next input, with the errors measured at the final time $t=1$. The RFF rows use the random-Fourier-feature data \eqref{eq:app-grf}. Best errors are in {\bf bold}.}
\label{tab:operators}
\end{table}

\subsection{Computational cost} \label{sec:numerics-cost}

\begin{figure}[!h]
\centering
\includegraphics[width=\textwidth]{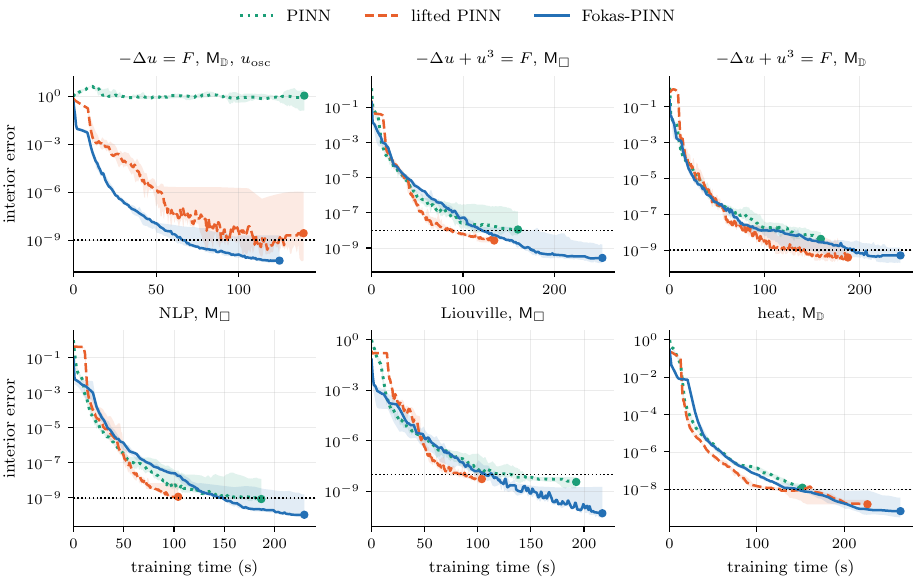}
\caption{Cost comparison for single PDE problems. The plots show the interior error on the monitor set against training time, with the bands corresponding to the range and the lines the medians. The dotted line marks the error target for the given problem as specified in Table \ref{tab:cost_operator}. }
\label{fig:cost_single}
\end{figure}

\begin{figure}[!h]
\centering
\includegraphics[width=0.85\textwidth]{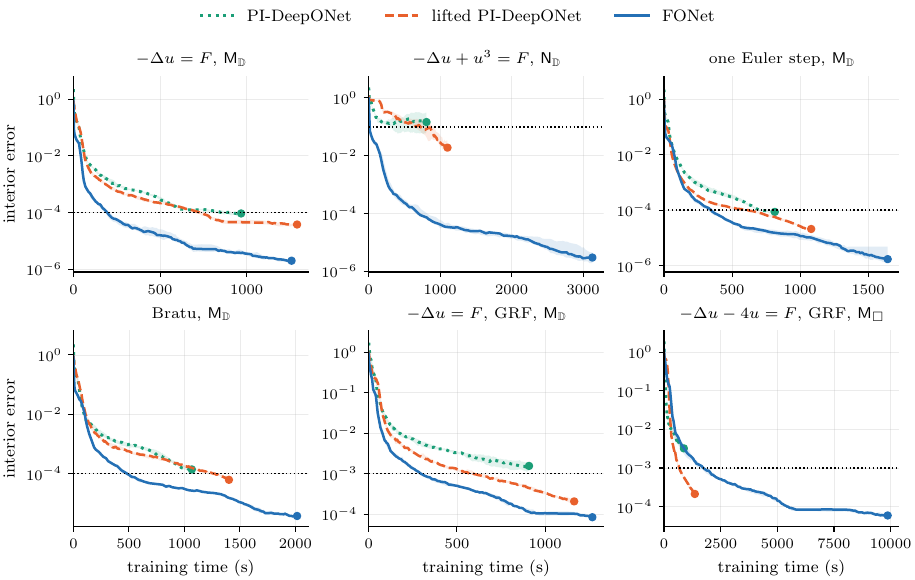}
\caption{Cost comparison for solution operators. The dotted line marks the error target for the given problem as specified in Table \ref{tab:cost_operator}.}
\label{fig:cost_operator}
\end{figure}

Finally, we discuss the computational cost of FokasNets compared to the established benchmarks. Figures~\ref{fig:cost_single} and~\ref{fig:cost_operator} show the interior error during training. The cost per evaluation, the training time and the time to a problem-specific target are also given in Tables~\ref{tab:cost_single} and~\ref{tab:cost_operator}. All times are medians over the five seeds, and the cost per evaluation is the training time divided by the number of loss evaluations, including those of the line search. The target of a problem is the smallest power of ten reached by the more accurate baseline, and the time to target is the training time at which the interior error first falls below it. In what follows, we state cost ratios relative to the lifted baselines, which share the lifted construction, so that the ratio isolates the cost of the completion. 

For single problems, one evaluation of the Fokas-PINN costs as much as one of the lifted PINN when the completion does not depend on the network (Poisson) and $\times1.3$-$2.0$ as much otherwise. It reaches the target first in two of the six problems and within $\times1.5$ of the faster baseline in the others, and attains the smallest final error in all six. In the 3D ball case one evaluation of the Fokas-PINN costs $\times 4.3$ and $\times 5.3$ as much as one of the lifted PINN. Hence, in these cases and with the given training configurations, the GR completion improves the attainable accuracy, but does not in general shorten the time to a given baseline accuracy.

This changes significantly for the solution operator problems, where one evaluation of the FONet costs approximately $\times 0.7$-$0.9$ as much as one of the lifted PI-DeepONet when the completion is a fixed linear map (Helmholtz, Poisson), $\times 1.3$-$1.7$ times as much for the nonlinear operators and $2.8$ times for $\mathsf N_{\mathbb D}$. Most importantly, the FONet reaches the target first for all six operators, $\times 1.7$ - $3.3$ earlier than the more accurate baseline, and $54$ times earlier for $\mathsf N_{\mathbb D}$. Many techniques can be considered to further speed-up the training and implementation of FokasNets. We consider one such approach in Section \ref{sec:numerics-gn-sec}, with more extensive analysis and experiments deferred  to future work.

\subsection{Gauss-Newton training with the implicit Jacobian}\label{sec:numerics-gn-sec}

The runs of Section~\ref{sec:numerics-cost} are trained from the loss and its gradient, for which the completion contributes
one adjoint solve \eqref{eq:fpinn-adjoint-gradient} per evaluation, and the curvature used by SSBroyden is built from successive
gradients. Proposition~\ref{prop:fpinn-implicit-differentiation}, however, provides more than the gradient of the loss. Specifically, the sensitivity \eqref{eq:fpinn-implicit-sensitivity} is the derivative of the completion itself, and through \eqref{eq:fpinn-field-sensitivity} it gives the derivative of every residual of the lifted approximation with respect to $\theta$. In this subsection, we use this to train a Fokas-PINN with a Gauss-Newton method, and we show what this changes in the
cost and in the attainable accuracy, for the single example $-\Delta u+u^3=F$ on $\mathsf M_{\mathbb D}$, with solution $S$ of Table~\ref{tab:app-problems}. The completion is the rule (R), solved by damped Newton from the linear solution and emitted through \eqref{eq:fpinn-ift-output}.

To see why the Jacobian is affordable, we write the loss \eqref{eq:loss-general} of the completed approximation as a sum of squares, $\mathcal J(\theta)=\|\mathbf r(\theta)\|_2^2$, where the vector $\mathbf r(\theta)\in\R^{N}$ collects the weighted PDE residuals at the $N_\Om$ collocation points and the weighted Neumann residuals at the $N_N$ nodes of $\GN$. A Gauss-Newton step then solves:
\begin{equation}\label{eq:numerics-gn-step}
 (\mathbf J^\top\mathbf J+\lambda I)\,\delta=\mathbf J^\top\mathbf r,\qquad \mathbf J:=D_\theta\mathbf r\in\R^{N\times N_\theta},
\end{equation}
and the question is the cost of $\mathbf J$. Automatic differentiation produces a Jacobian either one column at a time, by a forward sweep per parameter, or one row at a time, by a reverse sweep per residual, and with $N=2112$ and $N_\theta=1185$
neither is acceptable at every step. For the Fokas-PINN, the structure of $\mathbf r$ removes this difficulty. Consider the PDE residual at a collocation point $x_i$, $r_i=c_i\,[\Delta u_\theta-R(u_\theta)+F](x_i)$, with $c_i$ the square root of the
quadrature weight. Inserting $u_\theta=L(z)+\beta v_\theta$ and holding $z$ fixed, the derivative with respect to $\theta$ gives:
\begin{equation}\label{eq:numerics-gn-row}
 \partial_\theta r_i=c_i\Big[\beta\,\partial_\theta\Delta v_\theta+2\nabla\beta\cdot\partial_\theta\nabla v_\theta
   +\big(\Delta\beta-R'(u_\theta)\,\beta\big)\,\partial_\theta v_\theta\Big](x_i),
\end{equation}
i.e., a combination, with known coefficients, of the gradients with respect to $\theta$ of the four numbers $(v_\theta,\partial_xv_\theta,\partial_yv_\theta,\Delta v_\theta)(x_i)$, and the Neumann rows involve $\partial_\theta v_\theta$ at
the boundary nodes in the same way. Since a residual depends on the network only at its own point, these four gradients at all collocation points are obtained by four reverse sweeps in total, whatever the number of parameters, and they assemble the whole
matrix $\partial_\theta\mathbf r$.

The dependence of $\mathbf r$ on $\theta$ through the completion is what characterizes the Fokas-PINN, this allows us to use Proposition~\ref{prop:fpinn-implicit-differentiation}. By the chain rule and
\eqref{eq:fpinn-implicit-sensitivity}, we obtain:
\begin{equation}\label{eq:numerics-gn-jacobian}
 \mathbf J=\partial_\theta\mathbf r+\partial_z\mathbf r\;D_\theta z^*,
 \qquad D_\theta z^*=-(J_C^\star)^{-1}\,\partial_\theta\mathcal C\in\R^{m\times N_\theta},
\end{equation}
where $\partial_z\mathbf r\in\R^{N\times m}$ is the derivative of the residuals through the lift, i.e.,
$-c_iR'(u_\theta)\,D_zL(x_i)$ for a PDE row and $c_j\,\dn D_zL(x_j)$ for a Neumann row. The second term of
\eqref{eq:numerics-gn-jacobian} is therefore a matrix of rank at most $m$, and none of its factors requires new computation. Specifically, $\partial_\theta\mathcal C$ is, for the rule (R), the moment matrix applied to $R'(u_\theta)\beta\,\partial_\theta v_\theta$, i.e., a weighted sum of the gradients already computed for
\eqref{eq:numerics-gn-row}, $J_C^\star$ is the $m\times m$ Jacobian factored in the last Newton iteration of the completion and $D_zL$ is the lift basis at the collocation points. Hence, inside a second-order method the completion adds to the cost of a step an amount set by $m$, the size of the trace, and not by the number of collocation points or of parameters. This is the opposite of the first-order situation of Table~\ref{tab:cost_single}, where every evaluation of the loss pays for the completion in full.

The same Jacobian also provides the initialization. If the first layer of a one-hidden-layer network $v_\theta$ is drawn at random and frozen, the four numbers in \eqref{eq:numerics-gn-row} are linear in the output weights, and $\mathbf r$ depends on them nonlinearly only through $R(u_\theta)$ and the completion. A few Gauss--Newton steps on the output weights alone, with the corresponding columns of \eqref{eq:numerics-gn-jacobian}, therefore fit the network before training, as in extreme learning machines \cite{huang2006extreme,dwivedi2020physics,schiassi2021extreme}, with the difference that the completion is part of the fit and the trace is solved consistently with the random features. The first layer is then released and all weights are trained, by Gauss-Newton until its steps become short and then by the SSBroyden iteration of Section~\ref{sec:numerics-setting}, initialized with the last Gauss--Newton approximation of the inverse Hessian as in hybrid
methods for nonlinear least squares \cite{fletcher1987hybrid}.

\begin{figure}[!h]
\centering
\includegraphics[width=\textwidth]{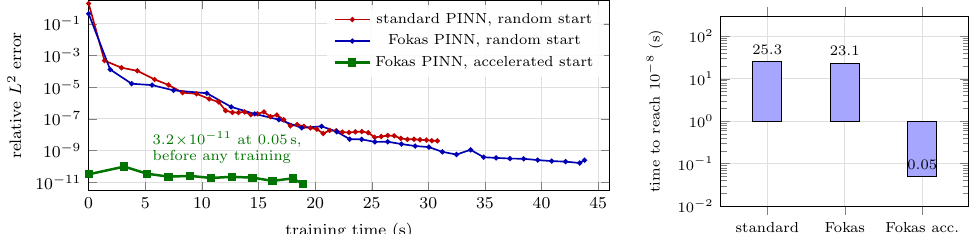}
\caption{Left: error versus training time (measured runs, CPU). Right: first logged time at which the error is $\le10^{-8}$ (random-start runs: $25$\,s and $23$\,s, upper bounds since the error is logged every 250 iterations; accelerated: already true at the start). Only the accelerated run ever reaches $10^{-10}$.}
\label{fig:gn}
\end{figure}

Table~\ref{tab:gn} and Figure~\ref{fig:gn} report the resulting runs, here only with one seed per case, on a CPU in double precision. From a random initialization, Gauss-Newton with the Jacobian \eqref{eq:numerics-gn-jacobian} brings the Fokas-PINN to $2.5\times10^{-10}$ in $44$\,s, an accuracy that the first-order runs of Table~\ref{tab:cost_single} reach after $240$\,s on a GPU. With the random-feature start and $296$ hidden units, the Fokas-PINN is at $3.2\times10^{-11}$ before training and at $7.4\times10^{-12}$ after $19$\,s. Note that the random-feature start accelerates the standard PINN as well, so the gain in wall-clock time is not to be attributed to the completion; what the completion changes is the accuracy obtained per parameter. At the same width the standard PINN starts at $5\times10^{-8}$ and training does not improve it, and only with $600$ hidden units, i.e., twice the parameters, does it reach the accuracy of the Fokas-PINN.

\begin{table}[!h]
\centering\small
\begin{tabular}{@{}llrccr@{}}\toprule
method & start & $N_\theta$ & $E_\Om$ after start & $E_\Om$ final & time (s)\\\midrule
PINN        & random                                   & 1185 & $1.9$               & $4.1\times10^{-9}$  & 30.8\\
Fokas-PINN  & random                                   & 1185 & $0.43$              & $2.5\times10^{-10}$ & 43.8\\
PINN        & accelerated, $m_{\rm h}=296$$^\dagger$   & 1185 & $5.0\times10^{-8}$  & $1.9\times10^{-8}$  & 15.9\\
PINN        & accelerated, $m_{\rm h}=600$$^\dagger$   & 2401 & $5.5\times10^{-11}$ & $5.5\times10^{-11}$ & 16.5\\
Fokas-PINN  & accelerated, $m_{\rm h}=296$             & 1185 & {\bf $3.2\times10^{-11}$} & {\bf $7.4\times10^{-12}$} & 18.9\\\bottomrule
\end{tabular}
\caption{Second-order training on $\mathsf M_{\mathbb D}$, $-\Delta u+u^3=F$. Relative interior error $E_\Om$ after the start
(the random initialization, or the random-feature fit with $m_{\rm h}$ hidden units) and at the end of training, and the
wall-clock time on a CPU (Apple M5, double precision, one seed per case). $^\dagger$Polish with an Armijo backtracking search
instead of the strong-Wolfe search of the other rows.}
\label{tab:gn}
\end{table}

With these results, we show that the implicit Jacobian makes second-order training of FokasNets extremely efficient, able to achieve significant speed-ups compared to standard PINNs, as well as the higher accuracies. A comprehensive study that includes the solution operator problems and the lifted neural network benchmark across many examples is left for future work.

\section{Conclusion}\label{sec:conclusion}

In this work we developed FokasNets, which are Physics-Informed Neural Networks and Operators that add an additional layer to the training process which solves the Fokas method's Global Relation. The Global Relation provides a deterministic integral equation that couples the boundary values and is used to complete the unknown boundary traces during the neural approximation training. Hence, at every evaluation of the loss a finite Global Relation is solved for the unknown trace, which is then used to create a lifted solution approximation that satisfies the value traces identically. 

FokasNets address a known source of difficulty in Physics-Informed training, namely, that the imposition of boundary conditions as a soft penalty in the loss function can often give rise to large errors both on the boundary and in the interior. The Global Relation step of FokasNets allows the unknown trace to be computed from the known data and the current network, rather than learned, leading to more efficient solution representations. We provide explicit constructions of FokasNets for elliptic and evolution PDEs and their solution operators, using a DeepONet representation. Our numerical results show that FokasNets are able to produce significantly more accurate approximations across a suite of PDEs and solution operators. Importantly, in almost all the solution operator examples, the accuracy improvement is multiple orders of magnitude. We also show that, by taking advantage of the structure of FokasNets, appropriate second-order training schemes using the implicit Jacobian can be used to achieve significant speedups compared to both standard PINNs and the ``naive'' training implementation of the FokasNet.

Many exciting future research directions can be explored. Specifically, the proposed methods and construction of FokasNets are generalizable to many other families of PDEs and their corresponding Global Relation formulations (e.g., integrable nonlinear equations for which the Fokas method was first developed). Furthermore, FokasNets can be combined with various neural network solution representations, using different lifting mechanisms (e.g., from methods from the field of TFC), or even different approaches in incorporating the Global Relation solver into the training process. Finally, extending FokasNets by further optimizing computational efficiency and considering more general domain geometries is an important next step.

\section*{Acknowledgments}
K. C. Georgiou was supported by the Defense Advanced Research Projects Agency (DARPA) under Agreement No. HR00112590073, the Office of Naval Research (ONR) MURI Grant No. N00014-24-1-2545 through Purdue University Subaward F9000086602080, and the U.S. Department of Energy, Office of Science, Office of Advanced Scientific Computing Research (ASCR) under Award Number DE-SC0023191. G. Em Karniadakis and K. Shukla were supported by the Pacific Northwest National Laboratory (operated by Battelle Memorial Institute for the U.S. Department of Energy under Prime Contract DE-AC05-76RL01830) under Joint Appointment Agreement No. 801034 (PNNL JAA) and Subcontract No. 831126 for the SciDAC Learning-Accelerated Domain Science project (SciDAC LEADS).

Approved for public release; distribution is unlimited. This report was prepared as an account of work sponsored by an agency of the United States Government. Neither the United States Government nor any agency thereof, nor any of their employees, makes any warranty, express or implied, or assumes any legal liability or responsibility for the accuracy, completeness, or usefulness of any information, apparatus, product, or process disclosed, or represents that its use would not infringe privately owned rights. Reference herein to any specific commercial product, process, or service by trade name, trademark, manufacturer, or otherwise does not necessarily constitute or imply its endorsement, recommendation, or favoring by the United States Government or any agency thereof. The views and opinions of authors expressed herein do not necessarily state or reflect those of the United States Government or any agency thereof.

\appendix

\section{Experimental setup and additional results}
\label{app:setup}

\paragraph{Manufactured solutions.} With $(\tilde x,\tilde y)=(x,y)$ on the disk and $(\tilde x,\tilde y)=(2x-1,2y-1)$ on the square, we consider the following manufactured solutions:
\begin{equation}
\label{eq:numerics-S}
\begin{gathered}
 S=\sum_{k=1}^{8}c_ks_k(\tilde x,\tilde y),
 \qquad
 c=(0.25,\ 0.2,\ -0.1,\ 0.12,\ 0.08,\ 0.07,\ 0.025,\ 0.015),\\
 \mathbf s=\bigl(1,\ \tilde x,\ \tilde y,\ \sin 2\tilde x\cos \tilde y\,e^{0.3\tilde y},\ \cos 3\tilde y\sin \tilde x,\ e^{0.5\tilde x}\cos 2\tilde y,\ \sin(4\tilde x+\tilde y),\ \cos(\tilde x-3\tilde y)\bigr),
\end{gathered}
\end{equation}
\begin{equation}
\label{eq:numerics-osc}
 u_{\rm osc}=\operatorname{Re}(x+\mathrm{i}y)^{16}+0.1(1-r^2)^2 ,
\end{equation}
so that $-\Delta u_{\rm osc}=0.8-1.6r^2$ is smooth. For time stepping, $u^*=S+tP$ (heat) and $u^*=S_1+tP_1$ (KdV), with:
\begin{equation}
\label{eq:numerics-P}
\begin{gathered}
 P=\sum_{k=1}^{8}c_kp_k(\tilde x,\tilde y),\qquad
 \mathbf p=\bigl(1,\ \tilde x,\ \tilde y,\ \tilde x\tilde y,\ \tilde x^2-\tilde y^2,\ (\tilde x^2+\tilde y^2)^2,\ \tilde x^4\tilde y^2,\ \tilde y^8\bigr),\\
 S_1=0.3\sin2x+0.2\cos3x\,e^{0.2x}+0.1x^2,\qquad P_1=0.25+0.2x-0.1x^2+0.12x^3 .
\end{gathered}
\end{equation}

For the 3D example, with $x=(x_1,x_2,x_3)$, we construct the solution:
\begin{equation}
\label{eq:app-S3}
\begin{aligned}
 S_3={}&0.25+0.2x_1-0.1x_2+0.15x_3+0.12\sin(2x_1)\cos(x_2)e^{0.3x_3}+0.08\cos(3x_2)\sin(x_1)\,x_3\\
 &+0.07e^{0.5x_1}\cos(2x_2)\sin(x_3)+0.025\sin(4x_1+x_2-x_3)+0.015\cos(x_1-3x_2+2x_3).
\end{aligned}
\end{equation}

On the parallelogram $\mathcal P$ of Appendix~\ref{app:stat-par}, in its sheared coordinates $(\xi,\eta)=(x,y-sx)$,
\begin{equation}
\label{eq:app-par}
\begin{gathered}
 u^*=T(\xi,\eta)+\tfrac12\beta(\xi,\eta)W(\xi,\eta),\qquad
 T=0.3\xi^2\eta+0.2\xi\cos2\eta+0.1e^{\eta}+0.15\xi^3(1-\eta)^2,\\
 W=\sin(2\xi)\cos(\eta)e^{0.3\eta}+0.5\cos(\xi-3\eta),\qquad \beta=16\xi(1-\xi)\eta(1-\eta).
\end{gathered}
\end{equation}
Since $\beta W$ vanishes on the boundary, the traces of $u^*$ on the slanted sides are those of $T$, polynomials of degree
three in $\xi$. The solution is smooth: it does not contain the corner mode $r^{3/4}\sin(3\vartheta/4)$ of the $120^\circ$ corners
(Appendix~\ref{app:stat-par}).

For the manufactured evolution problems, we consider the following:
\begin{equation}
\label{eq:app-st-solutions}
\begin{aligned}
 u^*_{\rm heat}&=(0.3+x+0.2xy)\,e^{-t}+0.15\sin t\,(x^2-y^2)+0.1\,e^{-2t}x^3,\\
 u^*_{\rm KdV}&=e^{-t}\bigl(0.3\sin2x+0.2\cos3x\,e^{0.2x}\bigr)+0.15\sin t\,(x^2-0.1x^3)+0.1 .
\end{aligned}
\end{equation}

\paragraph{Operator data.} A datum of the manufactured families is generated by:
\begin{equation}
\label{eq:numerics-op-family}
 u^*=\sum_{k=1}^{18}\xi_ks_k,
 \qquad
 \xi_k\sim\mathcal U(-\sigma_k,\sigma_k),
\end{equation}
with $s_1,\ldots,s_8$ of \eqref{eq:numerics-S}, $\sigma_k=|c_k|$, and $s_9,\ldots,s_{18}$ the harmonic polynomials $\operatorname{Re}(\tilde x+\mathrm{i}\tilde y)^n$, $\operatorname{Im}(\tilde x+\mathrm{i}\tilde y)^n$, $n=1,\ldots,5$, with amplitudes
$0.05/n$. The training, validation and test data use the modes with $n\le4$. 

\par The time-stepping operator is trained on
$u^*=\sum_k(\xi_k+t\zeta_k)s_k$, $\zeta_k\sim\mathcal U(-\sigma_k/2,\sigma_k/2)$, at random levels $t_\ell\in\{0.1,\ldots,1\}$ with
the exact previous state, and the space-time operator on $u^*=\sum_k\xi_ks_k(x)\pi_k(t)$ with the time profiles
$\pi_k\in\{e^{-t},\cos t,e^{-2t},1+\tfrac12\sin t\}$ assigned cyclically. The random-Fourier-feature
(RFF) families of the two linear operators are:
\begin{equation}
\label{eq:app-grf}
 u^*=\sum_{j=1}^{24}\xi_j\cos(w_j\cdot x+b_j),\qquad \xi_j\sim\mathcal U(-\sigma,\sigma),\quad \sigma=\tfrac12\sqrt{2/24},
\end{equation}
with frequencies $w_j\sim\mathcal N(0,\ell^{-2}I)$, $\ell=0.5$, and phases $b_j\sim\mathcal U(0,2\pi)$ fixed for the family. Table \ref{tab:app-problems} gathers all the problems considered in Section \ref{sec:numerics}, their domains, settings and solutions.

Regarding training, the key configurations are given in Table \ref{tab:app-training}. An NVIDIA RTX A6000 GPU was used for the single PDE problems and NVIDIA Quadro RTX 6000 for the solution operator examples. For all networks we use Xavier-normal weights and zero biases. Training is performed by Adam with learning rate $10^{-3}$ followed by the self-scaled Broyden method (SSBroyden) with a strong-Wolfe line search \cite{kiyani2025optimizing,urban2025unveiling}, using double precision. The validation loss is evaluated every 50 iterations for elliptic problems and every 100 iterations for evolution problems and operators, on held-out points distinct from the nodes (512 interior and 128 boundary points for single problems and 1024 and 256 points on each validation datum for operators), and the checkpoint with the smallest validation loss is reported. Test points are disjoint from the training and validation points. Finally, the operator branch input is a PCA encoding of the datum with at most 16 components per data block. The level operator is rolled out on 128 unseen trajectories.

\begin{table}[!h]\centering\footnotesize\setlength{\tabcolsep}{3pt}
\begin{adjustbox}{max width=\linewidth}
\begin{tabular}{@{}>{\raggedright\arraybackslash}p{2.6cm}>{\raggedright\arraybackslash}p{3.3cm}>{\raggedright\arraybackslash}p{3.6cm}>{\raggedright\arraybackslash}p{3.7cm}>{\raggedright\arraybackslash}p{2.2cm}@{}}
\toprule
problem & equation & domain, data & solution / data family & source\\
\midrule
\multicolumn{5}{@{}l}{\textit{Elliptic (Section~\ref{sec:numerics-elliptic})}}\\[1pt]
Poisson, oscillatory & $-\Delta u=F$ & $\mathbb D$, $\mathsf M_{\mathbb D}$ & $u_{\rm osc}$ \eqref{eq:numerics-osc} & manufactured\\
parallelogram & $-\Delta u=F$ & $\mathcal P$, Dirichlet on the vertical, Neumann on the slanted sides & $T+\tfrac12\beta W$ \eqref{eq:app-par} & manufactured\\
cubic & $-\Delta u+u^3=F$ & $\square$, $\mathsf M_\square$; $\mathbb D$, $\mathsf N_{\mathbb D}$; $\mathbb D$, $\mathsf A_\alpha$ & $S$ \eqref{eq:numerics-S} & manufactured\\
exponential & $-\Delta u+e^u=F$ & $\mathbb D$, $\mathsf M_{\mathbb D}$ & $S$ & manufactured\\
nonlinear Poisson & $-\Delta u+e^u=-r$ & $\square$, $\mathsf M_\square$ & $1+\sin(\pi x)\cos(\pi y)$ & \cite{urban2025unveiling}\\
Liouville & $-\Delta u-\lambda e^u=0$, $\lambda=1$ & $\square$, $\mathsf M_\square$ & $\log\frac{8|f'|^2}{\lambda(1+|f|^2)^2}$, $f(z)=z+\tfrac12z^2$ & \cite{liouville1853equation,gel1963some}\\
ball & $-\Delta u+u^3=F$ & $\mathbb B\subset\R^3$, Neumann data & $S_3$ \eqref{eq:app-S3} & manufactured\\
ball, mixed & $-\Delta u+u^3=F$ & $\mathbb B$, Dirichlet on $\{x_3\ge0\}$, Neumann on $\{x_3<0\}$ & $S_3$ & manufactured\\
\midrule
\multicolumn{5}{@{}l}{\textit{Evolution (Section~\ref{sec:numerics-evolution})}}\\[1pt]
heat, time stepping & $u_t-\Delta u+u^3=F$ & $\square\times(0,1]$, $\mathsf M_\square$ & $S+tP$ \eqref{eq:numerics-P} & manufactured\\
KdV, time stepping & $u_t+uu_x+u_{xxx}=F$ & $(-1,1)\times(0,1]$, $\mathsf M_{\mathrm{KdV}}$ & $S_1+tP_1$ & manufactured\\
heat, space--time & $u_t-\Delta u+u^3=F$ & $\mathbb D\times(0,0.4]$, $\mathsf M_{\mathbb D}$ & $u^*_{\rm heat}$ \eqref{eq:app-st-solutions} & manufactured\\
KdV, space--time & $u_t+uu_x+u_{xxx}=F$ & $(-1,1)\times(0,0.25]$, $\mathsf M_{\mathrm{KdV}}$ & $u^*_{\rm KdV}$ & manufactured\\
PINNacle heat & $u_t-\Delta u=F$ & $\mathbb D\times(0,0.4]$, $\mathsf M_{\mathbb D}$ & $e^{|x|^2/2+2t}$ & \cite{hao2024pinnacle} (HNd, $d=2$)\\
\midrule
\multicolumn{5}{@{}l}{\textit{Solution operators (Section~\ref{sec:numerics-operators})}}\\[1pt]
cubic & $-\Delta u+u^3=F$, $(F,g_D,w_N)\mapsto u$ & $\mathbb D$, $\mathsf M_{\mathbb D}$; $\mathbb D$, $\mathsf N_{\mathbb D}$; $\square$, $\mathsf M_\square$ & \eqref{eq:numerics-op-family} & manufactured\\
cubic, $\lambda$ & $-\Delta u+\lambda u^3=F$, $\lambda\sim\mathcal U(0.2,2)$ & $\mathbb D$, $\mathsf M_{\mathbb D}$ & \eqref{eq:numerics-op-family} & manufactured\\
Bratu & $-\Delta u-\lambda e^u=F$, $\lambda\sim\mathcal U(0.25,0.75)$ & $\mathbb D$, $\mathsf M_{\mathbb D}$ & \eqref{eq:numerics-op-family} & \cite{gel1963some}\\
RFF Poisson & $-\Delta u=F$ & $\mathbb D$, $\mathsf M_{\mathbb D}$ & \eqref{eq:app-grf} & \cite{lu2021learning,wang2021learning}\\
RFF Helmholtz & $-\Delta u-4u=F$ & $\square$, $\mathsf M_\square$ & \eqref{eq:app-grf} & \cite{lu2021learning,wang2021learning}\\
level operator & $-\Delta u+u/\Delta t+u^3=\tilde F$, $\Delta t=0.1$ & $\mathbb D$, $\mathsf M_{\mathbb D}$ & levels of $\sum_k(\xi_k+t\zeta_k)s_k$ & manufactured\\
space--time & $u_t-\Delta u+u^3=F$, $(u_0,F,g_D,w_N)\mapsto u$ & $\mathbb D\times(0,0.4]$, $\mathsf M_{\mathbb D}$ & $\sum_k\xi_ks_k(x)\pi_k(t)$ & manufactured\\
\bottomrule
\end{tabular}
\end{adjustbox}
\caption{Problems of Section~\ref{sec:numerics}. Manufactured problems take $F$, $g_D$, $w_N$ (and $u_0$) from the stated solution and the literature benchmarks are posed with the mixed partitions of \eqref{eq:numerics-bc}.}
\label{tab:app-problems}
\end{table}

\begin{table}[!h]\centering\footnotesize\setlength{\tabcolsep}{3pt}
\begin{adjustbox}{max width=\linewidth}
\begin{tabular}{@{}>{\raggedright\arraybackslash}p{2.5cm}>{\raggedright\arraybackslash}p{2.9cm}>{\raggedright\arraybackslash}p{2.6cm}>{\raggedright\arraybackslash}p{4.0cm}>{\raggedright\arraybackslash}p{3.0cm}>{\raggedright\arraybackslash}p{1.6cm}@{}}
\toprule
problems & network & iterations (Adam + SSBroyden) & collocation nodes & test points & data (train / val / test)\\
\midrule
elliptic, $\mathbb D$ & $2\times32$, $\tanh$ & $1000+10\,000$ & $32$ Gauss radii $\times$ $65$ angles; $65$ boundary nodes & $9264+512$ & --\\
elliptic, $\square$ & $2\times32$, $\tanh$ & $1000+10\,000$ & $32\times32$ Gauss grid, $32$ Gauss nodes per side & $2304+256$ & --\\
parallelogram & $2\times32$, $\tanh$ & $1000+10\,000$ & $40\times40$ Gauss grid (sheared), $24$ nodes per side & $2304+256$ & --\\
ball & $2\times32$, $\tanh$ & $1000+10\,000$ & $16$ Gauss radii $\times$ $882$ directions, $882$ boundary nodes & $28\,224+1176$ & --\\
time stepping, heat & $2\times32$ per level & $300+1500$ per level & as elliptic, $\square$ & $2304+256$ at $t=1$ & --\\
time stepping, KdV & $2\times32$ per level & $300+1500$ per level & $64$ Gauss nodes & $64+2$ at $t=1$ & --\\
space--time, heat & $2\times32$, input $(x,t)$ & $1000+5000$ & $14\times29$ disk nodes $\times$ $10$ Gauss times & $12$ Gauss times $\times(2328+128)$ & --\\
space--time, KdV & $2\times32$, input $(x,t)$ & $1000+5000$ & $40$ Gauss nodes $\times$ $12$ Gauss times & $12$ Gauss times $\times(64+2)$ & --\\
stationary and level operators & branch--trunk, $p_{\rm op}=32$, trunk $2\times24$, branch $2\times32$ & $1000+10\,000$ & as the single problems & $4096+256$ per datum & $512/128/256$\\
space-time operator & as above, trunk input $(x,t)$ & $1000+10\,000$ & $14\times29$ disk nodes $\times$ $10$ Gauss times & $4096+256$ per datum & $256/64/128$\\
\bottomrule
\end{tabular}
\end{adjustbox}
\caption{Training and testing configurations for the problems of Section \ref{sec:numerics}.}
\label{tab:app-training}
\end{table}

\begin{table}[!t]\centering\footnotesize\setlength{\tabcolsep}{4pt}
\begin{tabular}{@{}llcccc@{}}
\toprule
problem & method & ms/eval & training time (s) & time to target (s) & final error \\
\midrule
\multicolumn{6}{@{}l}{\textit{Manufactured examples}}\\[1pt]
$-\Delta u=F$ & PINN & 7.4 & 140 & -- & $1.2{\times}10^{0}$ \\
{\scriptsize $\mathsf{M}_{\mathbb{D}}$, $u^*=u_{\rm osc}$} & lifted PINN & 8.5 & 139 & 89\,(4/5) & $8.0{\times}10^{-10}$ \\
{\scriptsize target $10^{-9}$} & Fokas-PINN & 8.5 & 125 & \textbf{65}\,(4/5) & $\mathbf{4.0{\times}10^{-11}}$ \\
\addlinespace[3pt]
$-\Delta u+u^3=F$ & PINN & 7.4 & 160 & 158\,(3/5) & $1.1{\times}10^{-8}$ \\
{\scriptsize $\mathsf{M}_{\square}$} & lifted PINN & 8.7 & 134 & \textbf{84} & $2.5{\times}10^{-9}$ \\
{\scriptsize target $10^{-8}$} & Fokas-PINN & 16.9 & 252 & 121 & $\mathbf{2.5{\times}10^{-10}}$ \\
\addlinespace[3pt]
$-\Delta u+u^3=F$ & PINN & 7.4 & 159 & -- & $2.8{\times}10^{-9}$ \\
{\scriptsize $\mathsf{M}_{\mathbb{D}}$} & lifted PINN & 8.5 & 187 & \textbf{131} & $3.5{\times}10^{-10}$ \\
{\scriptsize target $10^{-9}$} & Fokas-PINN & 11.3 & 243 & 174\,(4/5) & $\mathbf{2.2{\times}10^{-10}}$ \\
\midrule
\multicolumn{6}{@{}l}{\textit{Benchmarks from the literature}}\\[1pt]
$-\Delta\varphi+e^{\varphi}=-r$ & PINN & 8.7 & 187 & 163\,(4/5) & $7.7{\times}10^{-10}$ \\
{\scriptsize NLP, $\mathsf{M}_{\square}$} & lifted PINN & 10.2 & 104 & 97\,(2/5) & $1.8{\times}10^{-9}$ \\
{\scriptsize target $10^{-9}$} & Fokas-PINN & 18.0 & 230 & \textbf{146} & $\mathbf{1.2{\times}10^{-10}}$ \\
\addlinespace[3pt]
$-\Delta u-\lambda e^{u}=0$ & PINN & 8.9 & 193 & 128 & $3.7{\times}10^{-9}$ \\
{\scriptsize Liouville, $\mathsf{M}_{\square}$} & lifted PINN & 10.1 & 104 & \textbf{83} & $6.8{\times}10^{-9}$ \\
{\scriptsize target $10^{-8}$} & Fokas-PINN & 19.7 & 217 & 109 & $\mathbf{4.0{\times}10^{-11}}$ \\
\addlinespace[3pt]
$u_t-\Delta u=F$ & PINN & 12.8 & 152 & 150\,(1/5) & $1.4{\times}10^{-8}$ \\
{\scriptsize $\mathsf{M}_{\mathbb{D}}$} & lifted PINN & 12.3 & 227 & \textbf{127} & $1.7{\times}10^{-9}$ \\
{\scriptsize target $10^{-8}$} & Fokas-PINN & 21.2 & 264 & 145 & $\mathbf{6.4{\times}10^{-10}}$ \\
\midrule
\multicolumn{6}{@{}l}{\textit{Three-dimensional problems}}\\[1pt]
$-\Delta u+u^3=F$ & PINN & 6.5 & 134 & \textbf{99} & $3.5{\times}10^{-8}$ \\
{\scriptsize $\mathbb{B}$, Neumann data} & lifted PINN & 7.4 & 169 & 129 & $2.6{\times}10^{-8}$ \\
{\scriptsize target $10^{-7}$} & Fokas-PINN & 31.6 & 655 & 339 & $\mathbf{7.1{\times}10^{-9}}$ \\
\addlinespace[3pt]
$-\Delta u+u^3=F$ & PINN & 7.0 & 145 & -- & $1.5{\times}10^{-7}$ \\
{\scriptsize $\mathbb{B}$, mixed data} & lifted PINN & 7.8 & 165 & \textbf{94} & $1.2{\times}10^{-8}$ \\
{\scriptsize target $10^{-7}$} & Fokas-PINN & 41.6 & 855 & 432 & $\mathbf{5.3{\times}10^{-9}}$ \\
\bottomrule
\end{tabular}
\caption{Computational cost comparison for the single problems of Section~\ref{sec:numerics-cost}. In parentheses we give the number of seeds that reach the target (omitted when all the seeds reach the threshold). }
\label{tab:cost_single}
\end{table}

\begin{table}[!t]\centering\footnotesize\setlength{\tabcolsep}{4pt}
\begin{tabular}{@{}llcccc@{}}
\toprule
problem & method & ms/eval & training time (s) & time to target (s) & final error \\
\midrule
\multicolumn{6}{@{}l}{\textit{Manufactured examples}}\\[1pt]
$-\Delta u+u^3=F$ & PI-DeepONet & 34.0 & 892 & 776\,(1/5) & $1.1{\times}10^{-4}$ \\
{\scriptsize $\mathsf{M}_{\mathbb{D}}$} & lifted PI-DeepONet & 38.1 & 1218 & 861 & $3.4{\times}10^{-5}$ \\
{\scriptsize target $10^{-4}$} & FONet & 59.9 & 1696 & \textbf{402} & $\mathbf{3.6{\times}10^{-6}}$ \\
\addlinespace[3pt]
$-\Delta u+u^3=F$ & PI-DeepONet & 33.3 & 808 & 376\,(2/5) & $1.4{\times}10^{-1}$ \\
{\scriptsize $\mathsf{N}_{\mathbb{D}}$} & lifted PI-DeepONet & 38.9 & 1102 & 698 & $2.7{\times}10^{-2}$ \\
{\scriptsize target $10^{-1}$} & FONet & 109.3 & 3130 & \textbf{13} & $\mathbf{2.7{\times}10^{-6}}$ \\
\addlinespace[3pt]
$-\Delta u+u/\Delta t+u^3=\tilde F$ & PI-DeepONet & 31.7 & 814 & 706 & $7.9{\times}10^{-5}$ \\
{\scriptsize $\mathsf{M}_{\mathbb{D}}$, one Euler step} & lifted PI-DeepONet & 35.3 & 1082 & 588 & $2.1{\times}10^{-5}$ \\
{\scriptsize target $10^{-4}$} & FONet & 59.2 & 1644 & \textbf{351} & $\mathbf{1.5{\times}10^{-6}}$ \\
\midrule
\multicolumn{6}{@{}l}{\textit{Benchmarks from the literature}}\\[1pt]
$-\Delta u-\lambda e^{u}=F$ & PI-DeepONet & 40.8 & 1066 & -- & $1.5{\times}10^{-4}$ \\
{\scriptsize $\mathsf{M}_{\mathbb{D}}$, $\lambda\in[0.25,0.75]$} & lifted PI-DeepONet & 44.6 & 1401 & 1280 & $6.4{\times}10^{-5}$ \\
{\scriptsize target $10^{-4}$} & FONet & 56.3 & 1610 & \textbf{459} & $\mathbf{3.9{\times}10^{-6}}$ \\
\addlinespace[3pt]
$-\Delta u=F$ & PI-DeepONet & 39.8 & 908 & -- & $1.3{\times}10^{-3}$ \\
{\scriptsize RFF, $\mathsf{M}_{\mathbb{D}}$} & lifted PI-DeepONet & 44.0 & 1164 & 572 & $2.2{\times}10^{-4}$ \\
{\scriptsize target $10^{-3}$} & FONet & 40.9 & 1267 & \textbf{301} & $\mathbf{8.5{\times}10^{-5}}$ \\
\addlinespace[3pt]
$-\Delta u-4u=F$ & PI-DeepONet & 41.4 & 875 & -- & $3.2{\times}10^{-3}$ \\
{\scriptsize RFF, $\mathsf{M}_{\square}$} & lifted PI-DeepONet & 44.4 & 1360 & 673 & $2.2{\times}10^{-4}$ \\
{\scriptsize target $10^{-3}$} & FONet & 29.2 & 992 & \textbf{203} & $\mathbf{5.9{\times}10^{-5}}$ \\
\bottomrule
\end{tabular}
\caption{Computational comparison for the solution operators of Section~\ref{sec:numerics-cost}.}
\label{tab:cost_operator}
\end{table}

\section{Stationary constructions}
\label{app:stationary}

This appendix gives, for each spatial domain, the trace representation, the lift and gate, the tests and the resulting finite
Global Relation \eqref{eq:general-joint-system}-\eqref{eq:residual}. The sizes, solvers and the checks of each assembly with the
exact solution are summarized in Table~\ref{tab:constructions}, and the lifts in Table~\ref{tab:app-lifts}. In Table \ref{tab:constructions}, we also report the trace recovered by the GR completion when the exact solution is inserted in the moments (this corresponds to the error $\delta$ of
Proposition~\ref{prop:ca-completion-error}, induced by the numerical approximation of the GR), as well as the condition number that enters the stability constant $C_S$.

\begin{table}[!h]\centering\footnotesize\setlength{\tabcolsep}{3pt}
\begin{adjustbox}{max width=\linewidth}
\begin{tabular}{@{}>{\raggedright\arraybackslash}p{2.2cm}>{\raggedright\arraybackslash}p{3.3cm}>{\raggedright\arraybackslash}p{3.6cm}cl>{\raggedright\arraybackslash}p{2.0cm}c>{\raggedright\arraybackslash}p{1.1cm}@{}}
\toprule
domain & used in & tests $\varphi_k$ & $n\times m$ & rule, solver & trace error, exact data & cond & details\\
\midrule
\multicolumn{8}{@{}l}{\textit{Spatial relations}}\\[1pt]
$\mathbb D$ & elliptic ($\mathsf M_{\mathbb D}$, $\mathsf N_{\mathbb D}$, $\mathsf A_\alpha$); disk operators, level operator & $\operatorname{Re},\operatorname{Im}(x+\mathrm iy)^k$, $k\le32$ & $32\times32$ to $65\times65$ & (R), Newton & $\le6\times10^{-14}$ & $28$--$170$ & \S\ref{app:stat-disk}\\
$\square$ & elliptic $\mathsf M_\square$, benchmarks on $\square$, square operators & $\sin(n\pi x)\{\cosh,\sinh\}(n\pi(y-\tfrac12))$, $n\le64$ & $128\times64$ & (LS), Gauss--Newton & $3.1\times10^{-12}$ & $6.0\times10^{3}$ & \S\ref{app:stat-square}\\
$\mathcal P$ & parallelogram & $\sin(n\pi x)\{\cosh,\sinh\}(n\pi(y-y_0))$, $n\le12$ & $24\times18$ & (LS), one solve & $1.2\times10^{-11}$ & $9.5\times10^{4}$ & \S\ref{app:stat-par}\\
$\mathbb B$ & ball & $r^lY_{lm}$, $l\le16$ & $289\times289$ & (R), Newton & $6.8\times10^{-11}$ & $12$ & \S\ref{app:stat-ball}\\
$\square$, level & heat, time stepping & as $\square$ & $128\times64$ & (LS), Gauss--Newton & $2.7\times10^{-11}$\,$^a$ & $38$ & \S\ref{app:evo-ts-heat}\\
$(-1,1)$, level & KdV, time stepping & $e^{\mu_jx}$, $\mu_j^3=1/\Delta t$ & $3\times3$ & (R), Newton & $4.9\times10^{-15}$\,$^a$ & $2.7$ & \S\ref{app:evo-ts-kdv}\\
\midrule
\multicolumn{8}{@{}l}{\textit{Space--time relations}}\\[1pt]
$\mathbb D\times(0,\mathcal T)$ & heat, PINNacle heat, space--time operator & $e^{-\kappa^2t}I_n(\kappa r)e^{\mathrm in\vartheta}$, 12 $\kappa$, 5 windows & $840\times112$ & (P), projected Newton & $3.8\times10^{-12}$ & $3.4\times10^{3}$ & \S\ref{app:evo-st-heat}\\
$(-1,1)\times(0,\mathcal T)$ & KdV & $e^{-\mathrm ikx-\mathrm ik^3t}$, 10 $k$, 5 windows & $100\times24$ & (P), projected Newton & $1.8\times10^{-12}$ & $5.4\times10^{2}$ & \S\ref{app:evo-st-kdv}\\
\bottomrule
\end{tabular}
\end{adjustbox}
\caption{The finite Global Relations of all problems. The size $n\times m$ counts the equations and the unknown trace coefficients, and the rule is the completion rule of \eqref{eq:fpinn-defining-map} with the solver used for it. The last two columns are computed with the exact solution in place of the network: the trace error is the relative error of the trace that the completion recovers from exact data, which is the accuracy floor of the construction, and cond is the condition number of the Jacobian $J_r=D_zr$ of the relation with respect to the unknown coefficients at the exact solution. For the time-stepping levels and the space-time relations it is the condition number of the linear part $\mathsf K$ of the relation. Where a row covers several problems, the range over them is given. $^a$: for the time-stepping relations the residual of the relation at the exact solutions reported, since a backward-Euler level has no exact trace of its own to compare with.}
\label{tab:constructions}
\end{table}

\begin{table}[!h]\centering\footnotesize\setlength{\tabcolsep}{3pt}
\begin{adjustbox}{max width=\linewidth}
\begin{tabular}{@{}>{\raggedright\arraybackslash}p{2.2cm}>{\raggedright\arraybackslash}p{3.0cm}>{\raggedright\arraybackslash}p{4.4cm}>{\raggedright\arraybackslash}p{4.0cm}>{\raggedright\arraybackslash}p{2.6cm}@{}}
\toprule
construction & completed trace & trace representation & lift $L$ & gate $\beta$\\
\midrule
$\mathbb D$ & $u$ on $\GN$ & nodal values at 65 equispaced boundary nodes (trigonometric interpolation) & harmonic extension (harmonic cardinal functions) & $1-r^2$\\
$\square$ & $u$ on $y\in\{0,1\}$ & nodal values at 32 Gauss nodes per side; prescribed corner values; corner-pinned interpolation & Coons transfinite interpolation of the four side traces & $16xy(1-x)(1-y)$\\
$\mathcal P$ & $u$ on the slanted sides & Legendre polynomials of degree $\le8$ per side & Coons interpolation in $(\xi,\eta)=(x,y-sx)$ & $16\xi(1-\xi)\eta(1-\eta)$\\
$\mathbb B$ & $u$ on $\partial\mathbb B$ & real spherical harmonics $Y_{lm}$, $l\le16$ & harmonic extension $\sum z_{lm}r^lY_{lm}$ & $1-|x|^2$\\
$(-1,1)$, level & $\partial_xu(-1)$, $\partial_x^2u(\pm1)$ & three numbers per level & quintic Hermite interpolant of the six endpoint traces & $(1-x^2)^3$\\
$\mathbb D\times(0,\mathcal T)$ & $u$ on $\GN\times(0,\mathcal T)$ & 14 nodes $\times$ $\psi_p(t)=\frac t{\mathcal T}P_{p-1}(\frac{2t}{\mathcal T}-1)$, $p\le8$ & $u_0$ $+$ harmonic extension of the trace increments & $(1-r^2)\,t/\mathcal T$\\
$(-1,1)\times(0,\mathcal T)$ & $\partial_xu(-1,t)$, $\partial_x^2u(\pm1,t)$ & $3\times\psi_p(t)$, $p\le8$ & $u_0$ $+$ quintic Hermite interpolant of the increments & $(1-x^2)^3\,t/\mathcal T$\\
\bottomrule
\end{tabular}
\end{adjustbox}
\caption{The trace representation, lift and gate of each construction. }
\label{tab:app-lifts}
\end{table}

Below we detail the constructions of the lifts and GR systems that were used in our numerical experiments. Throughout,
$Q_0=H^\top\operatorname{diag}(\omega_a)$ with $H_{ak}=\varphi_k(x_a)$ as in \eqref{eq:elliptic-matrices}.

\subsection{Disk}
\label{app:stat-disk}

On the disk the harmonic test functions and the interpolation basis of the boundary trace are the same Fourier family. This has
two consequences that make the disk the simplest construction: the tested GR \eqref{eq:GR} becomes the
Dirichlet-to-Neumann map of the Laplacian, written node by node, and the unknown flux on $\GD$ is removed by simply discarding the
equations of the Dirichlet nodes.

\paragraph{Trace.} The boundary nodes are the $M=65$ equispaced points $x_j^\partial=(\cos\vartheta_j,\sin\vartheta_j)$,
$\vartheta_j=2\pi(j-1)/M$, and the value and flux traces are stored by their nodal values $\mathbf g,\mathbf w\in\R^{65}$. As in
\eqref{eq:affine-lift-method}, the nodes on $\GD$ (index set $\pre$) carry the data, $\mathbf g_\pre=(g_D(x_j^\partial))_{j\in\pre}$,
and the nodes on $\GN$ (index set $\un$) carry the unknowns, $\mathbf g_\un=z$. For $\mathsf M_{\mathbb D}$ there are $32$ unknown
nodes, for $\mathsf A_{\pi/2}$ and $\mathsf A_{\pi/8}$ $48$ and $60$, and for $\mathsf N_{\mathbb D}$ all $65$.

\paragraph{Tests and lift.} Consider the $65$ harmonic polynomials:
\begin{equation}
\label{eq:app-disk-harmonics}
 \phi_1=1,\qquad \phi_{2\nu}=\operatorname{Re}(x+\mathrm{i}y)^\nu,\qquad \phi_{2\nu+1}=\operatorname{Im}(x+\mathrm{i}y)^\nu,\qquad \nu=1,\ldots,32,
\end{equation}
which restrict to the Fourier modes $1,\cos\nu\vartheta,\sin\nu\vartheta$ on the circle. First, we note that,
since $\partial_r(r^\nu e^{\mathrm i\nu\vartheta})=\nu\,r^{\nu-1}e^{\mathrm i\nu\vartheta}$, each $\phi_k$ is an eigenfunction of the
normal derivative on the unit circle, with eigenvalue equal to its degree:
\begin{equation}
\label{eq:app-disk-eig}
 \dn\phi_k=\varrho_k\,\phi_k\ \text{ on }\partial\mathbb D,\qquad \varrho_1=0,\quad \varrho_{2\nu}=\varrho_{2\nu+1}=\nu .
\end{equation}
Second, the trapezoidal rule on the $65$ nodes integrates every product $\phi_k\phi_{k'}$ exactly (a trigonometric polynomial of
degree at most $64$). Let $\mathcal B\in\R^{65\times65}$ be the matrix with entries $\mathcal B_{jk}=\phi_k(x_j^\partial)$, whose $k$-th column is the mode
$\phi_k$ evaluated at the $65$ boundary nodes. It is the trigonometric interpolation matrix, invertible, so that $\mathcal B\mathbf a$
are the nodal values of the trace with Fourier coefficients $\mathbf a$ and $\mathcal B^{-1}\mathbf g$ are the Fourier coefficients of
the trace with nodal values $\mathbf g$. The lift is the
harmonic extension of the trigonometric interpolant of the nodal values, and the gate vanishes on the circle:
\begin{equation}
\label{eq:app-disk-lift}
 L[\mathbf g](x)=\sum_{k=1}^{M}\phi_k(x)\,(\mathcal B^{-1}\mathbf g)_k=\sum_{j=1}^{M}\ell_j(x)\,g_j,\qquad
 \ell_j(x)=\sum_{k=1}^{M}\phi_k(x)\,(\mathcal B^{-1})_{kj},\qquad \beta(x)=1-|x|^2 ,
\end{equation}
so that in \eqref{eq:affine-lift-method} the $\ell_j$ are these harmonic cardinal functions and $\ell_0=0$. The lift takes the
value $g_j$ at node $j$ and is harmonic. The neural approximation \eqref{eq:fpinn-trial-field} therefore has exactly the nodal
boundary values $\mathbf g$, as required. The tests of the relation are the same polynomials, $\varphi_k=\phi_k$, $p=65$.

\paragraph{The Global Relation as a Dirichlet-to-Neumann map.} Insert $\varphi_k=\phi_k$ in \eqref{eq:GR} and evaluate the two boundary
integrals with the trapezoidal rule. By the two facts above, with $c_1=2\pi$ and $c_k=\pi$ for $k\ge2$, we obtain:
\[
 \qint_{\partial\mathbb D}\phi_k\,\dn u\,\mathrm ds=c_k\,(\mathcal B^{-1}\mathbf w)_k,\qquad
 \qint_{\partial\mathbb D}u\,\dn\phi_k\,\mathrm ds=\varrho_k\qint_{\partial\mathbb D}\phi_k\,u\,\mathrm ds=\varrho_k\,c_k\,(\mathcal B^{-1}\mathbf g)_k ,
\]
i.e., the boundary integrals pick out the $k$-th Fourier coefficient of the flux and $\varrho_k$ times the $k$-th Fourier coefficient
of the trace. The quadrature form of \eqref{eq:GR} for the test $\phi_k$ is therefore:
\[
 c_k(\mathcal B^{-1}\mathbf w)_k-\varrho_kc_k(\mathcal B^{-1}\mathbf g)_k=(Q_0\mathbf q)_k,\qquad k=1,\ldots,65 ,
\]
with $Q_0$ the volume quadrature of the tests as in \eqref{eq:elliptic-matrices} ($32$ Gauss--Legendre radii times the $65$ boundary
directions). These are $65$ equations in Fourier coefficients. Multiplying them by $\mathcal BC^{-1}$, $C=\operatorname{diag}(c)$,
returns to nodal quantities:
\begin{equation}
\label{eq:app-disk-nodal}
 \mathbf w=\mathsf D\,\mathbf g+\mathcal BC^{-1}Q_0\,\mathbf q,\qquad
 \mathsf D=\mathcal B\operatorname{diag}(\varrho)\mathcal B^{-1}.
\end{equation}
$\mathsf D$ maps the nodal values of a harmonic function to its nodal normal derivatives (mode $\nu$ is multiplied by $\nu$. Notice that this is the discrete Dirichlet-to-Neumann map of the Laplacian). Relation \eqref{eq:app-disk-nodal} states that the flux is the Dirichlet-to-Neumann map of the trace plus the contribution of the source, and it is the joint system \eqref{eq:general-joint-system} in nodal form.

Row $j$ of \eqref{eq:app-disk-nodal} contains the flux only through the nodal value $w_j$, which is the
datum $w_N(x_j^\partial)$ at the nodes of $\GN$ and the unknown flux $z'$ at the nodes of $\GD$. The annihilator $\mathsf T$ of
\eqref{eq:fpinn-annihilator} is therefore the selection of the rows $j\in\un$, and the reduced system
\eqref{eq:general-reduced-system} is:
\begin{equation}
\label{eq:app-disk-reduced}
 \mathsf Kz-\mathbf b+Q\,\mathbf q(z)=0,\qquad
 \mathsf K=\mathsf D_{\un,\un},\qquad Q=\bigl(\mathcal BC^{-1}Q_0\bigr)_{\un,:},\qquad \mathbf b=\mathbf w_N-\mathsf D_{\un,\pre}\,\mathbf g_\pre ,
\end{equation}
with the subscripts selecting rows and columns of the nodal matrices. It has as many equations as unknowns ($32\times32$ for
$\mathsf M_{\mathbb D}$, $48\times48$ for $\mathsf A_{\pi/2}$, $60\times60$ for $\mathsf A_{\pi/8}$) and is solved by Newton's method
(rule (R) of \eqref{eq:fpinn-defining-map}) with the Jacobian \eqref{eq:jacobian}, warm-started from the previous evaluation of the
loss. Note that, the level operator of Section~\ref{sec:numerics-operators} uses the same relation with the stage source
$q=u^3+u/\Delta t-\tilde F$, i.e.\ $\varsigma_\ell=0$ in \eqref{eq:time-stepping-stage}.

Finally, for the special case $\mathsf N_{\mathbb D}$ all $65$ rows are kept ($\mathsf T=I$, $\mathsf K=\mathsf D$). With Neumann data alone the linear part of the relation cannot determine the constant part of the trace, since a constant has zero normal derivative. The nonlinear term in the GR fixes this, since the $l=0$ row of the reduced system is the compatibility identity $\int_\Omega u^3=\int_\Omega F+\oint_{\partial\Omega}w_N$, which changes when a constant is added to $u$.

\subsection{Square}
\label{app:stat-square}

On a polygon there is no family of harmonic functions that turns the GR into a node-by-node map as on the disk, so the unknown
flux on $\GD$ cannot be removed by discarding equations. Instead, the tests are chosen to vanish on $\GD$: the flux there is then
multiplied by zero in \eqref{eq:GR} and never enters the relation ($\mathsf B'=0$ and $\mathsf T=I$ in
\eqref{eq:fpinn-annihilator}). Such tests span fewer directions than a free harmonic family, so more of them are used than there
are unknowns, and the completion is a least-squares solve.

\paragraph{Trace.} Each side of $\square$ carries the $32$ Gauss--Legendre nodes $t_j$ of $(0,1)$, $M=128$, and the four corner values
$c_{00},c_{10},c_{01},c_{11}$ are separate trace coordinates. For $\mathsf M_\square$ (Dirichlet data on $x\in\{0,1\}$, Neumann data on
$y\in\{0,1\}$) the corners and the nodes of the two vertical sides are prescribed by $g_D$, and the unknowns are the $m=64$ nodal
values on the two horizontal sides. On a side with parameter $t\in(0,1)$ and corner values $c_0$, $c_1$, the value trace is the
interpolant of the nodal values pinned to the corners:
\begin{equation}
\label{eq:app-square-side}
 g_s(t)=\sum_{j=1}^{32}\Lambda_j(t)\,g_{s,j}+\Lambda_0(t)\,c_0+\Lambda_{33}(t)\,c_1,\qquad
 \Lambda_j(t)=\frac{t(1-t)}{t_j(1-t_j)}\,\mathsf L_j(t),
\end{equation}
where $\mathsf L_j$ are the Lagrange polynomials of the Gauss nodes, $\Lambda_0(t)=(1-t)-\sum_j(1-t_j)\Lambda_j(t)$ and
$\Lambda_{33}(t)=t-\sum_jt_j\Lambda_j(t)$. The factor $t(1-t)$ makes every $\Lambda_j$ vanish at the corners, and $\Lambda_0$,
$\Lambda_{33}$ restore the linear function through the corner values, so that $g_s(t_j)=g_{s,j}$, $g_s(0)=c_0$ and $g_s(1)=c_1$.
The nodal values of a side therefore determine its trace polynomial of degree $33$ without disturbing the corners, which the
adjacent sides share.

\paragraph{Lift and gate.} The lift is the Coons transfinite interpolant of the four side traces $g_{\rm W},g_{\rm E},g_{\rm S},g_{\rm N}$
on $x=0$, $x=1$, $y=0$, $y=1$:
\begin{equation}
\label{eq:app-square-coons}
\begin{aligned}
 L[\hat g](x,y)={}&(1-x)\,g_{\rm W}(y)+x\,g_{\rm E}(y)+(1-y)\,g_{\rm S}(x)+y\,g_{\rm N}(x)\\
 &-\bigl[(1-x)(1-y)\,c_{00}+x(1-y)\,c_{10}+(1-x)y\,c_{01}+xy\,c_{11}\bigr].
\end{aligned}
\end{equation}
The first line blends the side traces linearly across the square (for example, on the side $x=0$, it gives
$g_{\rm W}(y)+(1-y)c_{00}+y\,c_{01}$, and the second line removes exactly the bilinear corner term $(1-y)c_{00}+y\,c_{01}$, so that
$L[\hat g]=g_{\rm W}$ there). Since every side trace is linear in the nodal and corner values, so is $L$, and the $\ell_j$ and
$\ell_0$ of \eqref{eq:affine-lift-method} are its coefficients of the unknown nodal values and of the prescribed values,
respectively. Note that, unlike the harmonic lift of the disk, $L$ is not harmonic and, hence, $\Delta L$ and $\dn L$ are evaluated analytically and enter
the residual \eqref{eq:residual}. The gate is $\beta=16\,x(1-x)\,y(1-y)$.

\paragraph{The Global Relation with tests vanishing on $\GD$.} The tests are, for $\nu=1,\ldots,64$, given by:
\begin{equation}
\label{eq:app-stat-square-tests}
 \varphi_\nu^{c}=\sin(\nu\pi x)\,\frac{\cosh\bigl(\nu\pi(y-\tfrac12)\bigr)}{\cosh(\nu\pi/2)},\qquad
 \varphi_\nu^{s}=\sin(\nu\pi x)\,\frac{\sinh\bigl(\nu\pi(y-\tfrac12)\bigr)}{\sinh(\nu\pi/2)},
\end{equation}
$p=128$. They are harmonic (obtained by separation of variables), they vanish on $x\in\{0,1\}=\GD$, and the denominators bound them by one on
$\overline\square$. The symmetric and antisymmetric pairs see, respectively, the sum and the difference of the traces on the two
horizontal sides. For such a test the GR \eqref{eq:GR} reads:
\begin{equation}
\label{eq:app-stat-square-gr}
 \sum_{s\in\{{\rm W},{\rm E},{\rm S},{\rm N}\}}\qint_0^1 g_s(t)\,\dn\varphi\bigl(x_s(t)\bigr)\,\mathrm dt+\qint_\square\varphi\,q\,\mathrm dx
 =\sum_{s\in\{{\rm S},{\rm N}\}}\qint_0^1\varphi\bigl(x_s(t)\bigr)\,w_N\bigl(x_s(t)\bigr)\,\mathrm dt ,
\end{equation}
with $x_s(t)$ the parameterization of side $s$. The flux appears only on the horizontal sides, where it is the datum $w_N$, whereas the trace appears on all four sides, as data on W and E and through the unknown nodal values on S and N. The side integrals are evaluated by
a $400$-point Gauss rule per side (the tests oscillate up to $\nu=64$ along a side). Inserting \eqref{eq:app-square-side}, the
coefficients of the unknown nodal values form the columns of $\mathsf K$, and the data (Dirichlet traces, corner values and
prescribed flux) form $\mathbf b$.

\par For the volume term, $q$ is known only through its values $\mathbf q$ on the $32\times32$ Gauss grid, which are also the interior
collocation points of the loss. The tests oscillate too fast ($\nu$ up to $64$) for the $32$-point Gauss rule to integrate
$\varphi_kq$ accurately, so instead $q$ is replaced by its tensor Lagrange interpolant on the grid, $q(x,y)\approx\sum_{i,j}q(x_i,y_j)\mathsf L_i(x)\mathsf L_j(y)$, which is accurate for a smooth $q$, and the tests are integrated
against this polynomial exactly. Since each test is a product of a function of $x$ and a function of $y$, the integrals factorize:
\begin{equation}
\label{eq:app-square-Q}
 Q_{k,(i,j)}=\int_0^1\sin(\nu_k\pi x)\,\mathsf L_i(x)\,\mathrm dx\int_0^1 C_{\nu_k}(y)\,\mathsf L_j(y)\,\mathrm dy ,
\end{equation}
with $C_{\nu}$ the hyperbolic factor of the test. These one-dimensional integrals are computed once with a fine rule, and the volume
term at every evaluation of the loss is the product $Q\mathbf q$.

The result is the $128\times64$ system \eqref{eq:residual}. Since it is rectangular, the least-squares completion depends on the
relative weight of the equations, and the normal derivatives of the tests grow like $\nu\pi$. Therefore, every row of $\mathsf K$, $Q$ and
$\mathbf b$ is divided by the norm of the corresponding row of $\mathsf K$, which amounts to normalizing each test. The row-scaled quantities are the ones denoted $\mathsf K$, $Q$, $\mathbf b$ and $r$ in \eqref{eq:fpinn-ls-objective}. The completion is the least-squares solution (rule (LS) of \eqref{eq:fpinn-defining-map}),
computed by Gauss-Newton iterations with reduced QR steps, warm-started from the previous evaluation of the loss, followed by one Newton step on the stationarity condition $J_r^\top r=0$, which is the step that produces $z_{\rm out}$ in \eqref{eq:fpinn-ift-output}. The same construction is used for the square benchmarks and the square operators.

\subsection{Parallelogram}
\label{app:stat-par} 

The parallelogram is the square construction on a sheared domain. It shows what happens to the tests vanishing on $\GD$ when the
Neumann sides are no longer aligned with the coordinates: the family is the same, but it determines fewer trace modes, and the
trace space has to be chosen accordingly.

\paragraph{Geometry and data.} $\mathcal P$ has the vertices $A=(0,0)$, $B=(1,s)$, $C=(1,1+s)$, $D=(0,1)$ with $s=-1/\sqrt3$, so
that the corners at $A$ and $C$ have angle $120^\circ$. The sheared coordinates $(\xi,\eta)=(x,\,y-sx)$ map $\mathcal P$ onto the
unit square, and a function of $(\xi,\eta)$ has the physical Laplacian $\partial_\xi^2-2s\,\partial_\xi\partial_\eta+(1+s^2)\partial_\eta^2$.
The Dirichlet data are prescribed on the vertical sides $DA$ ($\xi=0$) and $BC$ ($\xi=1$), the Neumann data on the slanted sides
$AB$ ($\eta=0$) and $DC$ ($\eta=1$), both taken from the manufactured solution \eqref{eq:app-par}, whose traces on the slanted
sides are cubic in $\xi$.

\paragraph{Trace, lift and gate.} On each slanted side the unknown value trace is expanded in Legendre polynomials,
$g_s(\xi)=\sum_{j=0}^{8}z_{s,j}P_j(2\xi-1)$, so $m=18$. The lift is the Coons interpolant \eqref{eq:app-square-coons} written in the
coordinates $(\xi,\eta)$, with the Dirichlet data on $DA$, $BC$, the Legendre expansions on $AB$, $DC$, and the corner values taken
from the data. Its Laplacian and normal derivatives are evaluated analytically through the shear. The gate is
$\beta=16\xi(1-\xi)\eta(1-\eta)$.

\paragraph{The Global Relation.} The tests are the family of the square, obtained by separation of variables and shifted to the
centre of $\mathcal P$:
\begin{equation}
\label{eq:app-par-tests}
 \varphi_\nu^{c,s}=\frac{\sin(\nu\pi x)\,\{\cosh,\sinh\}\bigl(\nu\pi(y-y_0)\bigr)}{\cosh(0.8\nu\pi)},\qquad y_0=\frac{1+s}{2},\qquad \nu=1,\ldots,12,
\end{equation}
$p=24$. They are harmonic and vanish on $x\in\{0,1\}=\GD$, so again $\mathsf B'=0$ and $\mathsf T=I$, and the GR \eqref{eq:GR}
reads:
\begin{equation}
\label{eq:app-par-gr}
 \qint_{AB\cup DC}u\,\dn\varphi\,\mathrm ds=\qint_{AB\cup DC}\varphi\,w_N\,\mathrm ds-\qint_{DA\cup BC}g_D\,\dn\varphi\,\mathrm ds
 -\qint_{\mathcal P}\varphi\,q\,\mathrm dx ,
\end{equation}
with the unknown trace on the left and the data and the volume term on the right. The side integrals use a $160$-point Gauss rule
per side and the volume integral the $40\times40$ Gauss grid in $(\xi,\eta)$, which also carries the interior collocation points of
the loss. With $\nu\le12$ the tests are mild enough for the Gauss rule to integrate $\varphi q$ directly. Inserting the Legendre
expansions gives the $24\times18$ system \eqref{eq:residual}, with the rows scaled as on the square. Since $R\equiv0$, the volume
term does not depend on $u$, and the completion is a single linear least-squares solve (rule (LS) of \eqref{eq:fpinn-defining-map}).

\subsection{Ball}
\label{app:stat-ball}

The ball with Neumann data on the whole sphere is the three-dimensional copy of the disk construction. The spherical harmonics
play the role of the Fourier modes, the solid harmonics $|x|^lY_{lm}$ the role of the harmonic polynomials, and the GR becomes
diagonal.

\paragraph{Trace, lift and gate.} The trace is expanded in the real spherical harmonics $Y_{lm}$ of degree $l\le16$, normalized to
unit root-mean-square on the sphere, so that $\int_{\partial\mathbb B}Y_{lm}Y_{l'm'}\,\mathrm ds=4\pi\,\delta_{ll'}\delta_{mm'}$. The
unknown is the vector $z\in\R^m$ of the $m=17^2=289$ coefficients $z_{lm}$, one for each degree $l\le16$ and order $|m'|\le l$, and
since no value is prescribed anywhere, there is
no prescribed part of the trace and no auxiliary flux. The lift is the harmonic extension of the trace, and the gate vanishes on
the sphere:
\begin{equation}
\label{eq:app-ball-lift}
 L[z](x)=\sum_{l\le16}\sum_{m}z_{lm}\,|x|^lY_{lm}\bigl(x/|x|\bigr),\qquad \beta(x)=1-|x|^2 .
\end{equation}
Each term $|x|^lY_{lm}(x/|x|)$ is a solid harmonic, a homogeneous harmonic polynomial of degree $l$ in $x$, so $L[z]$ is harmonic
and takes the boundary values $\sum z_{lm}Y_{lm}$.

\paragraph{The Global Relation.} The tests are the solid harmonics themselves, $\varphi_{lm}=|x|^lY_{lm}$, $p=289$. As on the disk,
each test is an eigenfunction of the normal derivative on the sphere, with eigenvalue equal to its degree:
\begin{equation}
\label{eq:app-ball-eig}
 \dn\varphi_{lm}=\partial_r\bigl(r^lY_{lm}\bigr)\big|_{r=1}=l\,Y_{lm} .
\end{equation}
Insert $\varphi_{lm}$ in \eqref{eq:GR} and expand the trace. By the orthogonality of the $Y_{lm}$, the trace integral picks out one
coefficient, and we obtain:
\begin{equation}
\label{eq:app-ball-gr}
 4\pi\,l\,z_{lm}-\qint_{\partial\mathbb B}Y_{lm}\,w_N\,\mathrm ds+\qint_{\mathbb B}\varphi_{lm}\,q\,\mathrm dx=0,\qquad l\le16 .
\end{equation}
This is the system \eqref{eq:residual} with $\mathsf K=4\pi\operatorname{diag}(l)$, diagonal, $\mathbf b$ the flux moments and $Q$ the volume quadrature of the solid harmonics. Every equation with $l\ge1$ determines one coefficient from the prescribed flux and the volume term. The equation with $l=0$ contains $z_{00}$ only through the reaction in $q$, exactly as for $\mathsf N_{\mathbb D}$ on the disk, and the system is solved in the same way, i.e., Newton's method (rule (R) of \eqref{eq:fpinn-defining-map}), started from the minimum-norm
solution of the linear part, with the rows $l\ge1$ divided by $4\pi l$.

The sphere integrals use Gauss-Legendre nodes in $\cos\vartheta$ times equispaced azimuths ($18\times49=882$ nodes), which
integrate products of harmonics of degree $\le16$ exactly, and the volume integrals use $16$ Gauss--Legendre radii times the same
$882$ directions, $14\,112$ nodes, which are also the interior collocation points of the loss.

\section{Evolution constructions}
\label{app:evolution-gr}

This appendix gives the constructions of the evolution problems: the space-time construction of Section~\ref{sec:evolution} for the heat and KdV equations, and the time-stepping construction of Section~\ref{prop:evolution-time-stepping} for the same two
equations. Below, we detail three changes in the space-time construction relative to the stationary constructions of Appendix~\ref{app:stationary}.

First, the relation is the space--time GR \eqref{eq:prelim-evo-GR}. Its tests solve the backward adjoint equation in
$\Omega\times(0,\mathcal T)$ rather than Laplace's equation, each test is paired with an upper time limit $\tau$, and besides the
boundary terms and the source moment the relation contains the terminal moment $\int_\Omega\varphi(x,\tau)u(x,\tau)\,\mathrm dx$,
which depends on the approximation. In the residual \eqref{eq:evo-residual} these are the source moment $\mathcal M_{\rm nl}$ and
the terminal moment $\mathcal M_{\rm end}$ of \eqref{eq:evo-source-moment}, while the initial moment
$\int_\Omega\varphi(x,0)u_0\,\mathrm dx$ is data and goes to $\mathbf b$.

Second, the traces are functions of time. Each unknown spatial trace coefficient carries $M_t$ temporal coefficients in the basis
$\boldsymbol\psi$ of \eqref{eq:evo-trace-basis}, for which we take $M_t=8$ and:
\begin{equation}
\label{eq:app-psi}
 \psi_p(t)=\frac{t}{\mathcal T}\,P_{p-1}\Bigl(\frac{2t}{\mathcal T}-1\Bigr),\qquad p=1,\ldots,8,
\end{equation}
Legendre polynomials multiplied by $t/\mathcal T$ so that $\psi_p(0)=0$ and the lift \eqref{eq:evo-lift} equals $u_0$ at $t=0$.
The temporal gate is $\beta_{\rm t}(t)=t/\mathcal T$, so that the network term $\beta(x)\beta_{\rm t}(t)v_\theta$ vanishes both on
the lateral boundary and at $t=0$.

Third, the system is overdetermined by design. Every test is used with the five nested upper limits $\tau_l=l\mathcal T/5$,
$l=1,\ldots,5$, so that it yields five equations, one per window $(0,\tau_l)$, and the equations outnumber the unknowns several
times. The completion is rule (P) of \eqref{eq:fpinn-defining-map} with the compression map $\Pi=\mathsf K_\star^{+}$, where $\mathsf K_\star=\mathsf K+\mathsf J_{\rm end}$ collects the complete linear action of the unknown traces, through the boundary terms and through the terminal moment (Section~\ref{sec:evo-completion-P}), and the projected equation is solved by Newton's method
\eqref{eq:evo-projected-newton}.

The time-stepping construction changes nothing of the above because it discretizes time first. Each backward-Euler level is the
stationary stage problem \eqref{eq:time-stepping-stage}, and the only choice is whether the term $u^\ell/\Delta t_\ell$ is kept in
the source ($\varsigma_\ell=0$, the stationary construction is reused unchanged) or in the operator ($\varsigma_\ell=1$, the
tests solve the modified equation). Spatial integrals are evaluated on the collocation nodes of Table~\ref{tab:app-training}, and
time integrals over a window $(0,\tau_l)$ by a Gauss--Legendre rule on that window.

\subsection{Space--time construction}
\label{app:evo-st}

\subsubsection{Heat equation on $\mathbb D\times(0,0.4]$}
\label{app:evo-st-heat}

The space-time heat construction follows that of the disk construction of Appendix~\ref{app:stat-disk} with one change. Specifically, the tests separate into a decaying exponential in time and a function of $x$ that solves the modified Helmholtz equation instead of Laplace's equation. The nodal Dirichlet-to-Neumann map of the disk is replaced by that of $\Delta-\kappa^2$, and everything else carries over.

\paragraph{Trace and lift.} Here $\mathcal L_x=-\Delta$ and the boundary form of \eqref{eq:prelim-evo-GR} is the one of the
elliptic GR, $\mathscr B_{\mathcal L}(u,\varphi)=\varphi\dn u-u\dn\varphi$. The value trace is stored at $M=29$ equispaced boundary nodes,
prescribed at the $15$ nodes of $\GD$ and unknown at the $14$ nodes of $\GN$. Each unknown nodal value is a function of time with
$M_t=8$ coefficients, so $m=14\times 8=112$, and the unknown nodal traces are
$\mathbf g_N(t)=\mathbf g_N(0)+Z\boldsymbol\psi(t)$ with $Z\in\R^{14\times8}$ and $z=\operatorname{\bf vec}(Z)$, as in
\eqref{eq:evo-trace-basis}. The spatial extension functions $\ell_j$ of the lift \eqref{eq:evo-lift} are the harmonic cardinal
functions \eqref{eq:app-disk-lift} of the disk for $29$ nodes, so at every $t$ the spatial part of the lift is the harmonic
extension of the nodal trace, and the gate is $\beta(x)\beta_{\rm t}(t)=(1-|x|^2)\,t/\mathcal T$.

\paragraph{Tests.} For twelve decay rates $\kappa_i^2=0.4+35.6\,(i-1)/11$, $i=1,\ldots,12$, the tests are:
\begin{equation}
\label{eq:app-heat-tests}
 \varphi_{i,k}(x,t)=e^{-\kappa_i^2t}\,\phi_{i,k}(x),\qquad
 \phi_{i,1}=\frac{I_0(\kappa_ir)}{I_0(\kappa_i)},\quad
 \phi_{i,2\nu}=\frac{I_\nu(\kappa_ir)}{I_\nu(\kappa_i)}\cos\nu\vartheta,\quad
 \phi_{i,2\nu+1}=\frac{I_\nu(\kappa_ir)}{I_\nu(\kappa_i)}\sin\nu\vartheta,
\end{equation}
$\nu=1,\ldots,14$, with $I_\nu$ the modified Bessel functions, so that there are $29$ spatial modes for each $\kappa_i$. Since
$\Delta\phi_{i,k}=\kappa_i^2\phi_{i,k}$, these tests satisfy $-\partial_t\varphi-\Delta\varphi=0$, the backward adjoint equation of
the heat operator, so $\varphi_{i,k}\in\mathcal H_T$. The range of $\kappa_i$ makes the time weights $e^{-\kappa_i^2t}$ vary from
slowly to rapidly decaying over $(0,\mathcal T)$, so that the windowed time integrals see the trace on different time scales. On
the circle $\phi_{i,k}$ are the Fourier modes \eqref{eq:app-disk-harmonics} of the disk, normalized to one at $r=1$, and, as on the disk, each is an eigenfunction of the normal derivative, now with the eigenvalue:
\begin{equation}
\label{eq:app-heat-eig}
 \dn\phi_{i,k}=\varrho_\nu(\kappa_i)\,\phi_{i,k}\ \text{ on }\partial\mathbb D,\qquad \varrho_\nu(\kappa)=\frac{\kappa I_\nu'(\kappa)}{I_\nu(\kappa)} ,
\end{equation}
which replaces the degree $\nu$ of \eqref{eq:app-disk-eig}.

\paragraph{The Global Relation.} Insert $\varphi_{i,k}$ with the upper limit $\tau_l$ in \eqref{eq:prelim-evo-GR}. The boundary
integrals are evaluated with the trapezoidal rule on the $29$ nodes, as on the disk, and give $e^{-\kappa_i^2t}$ times the $k$-th
Fourier coefficients of the flux and of $\varrho_\nu(\kappa_i)$ times the trace. Multiplying the $29$ equations of the family $i$ by
$\mathcal BC^{-1}$, exactly as in \eqref{eq:app-disk-nodal}, turns them into nodal equations:
\begin{equation}
\label{eq:app-heat-nodal}
 \int_0^{\tau_l}\!e^{-\kappa_i^2t}\bigl(\mathsf D_i\mathbf g(t)-\mathbf w(t)\bigr)\,\mathrm dt
 +\mathsf Q_i\Bigl[\int_0^{\tau_l}\!e^{-\kappa_i^2t}\,\mathbf q(t)\,\mathrm dt
   +e^{-\kappa_i^2\tau_l}\,\mathbf u_A(\tau_l)-\mathbf u_A(0)\Bigr]=0 ,
\end{equation}
where $\mathsf D_i=\mathcal B\operatorname{diag}(\varrho_\nu(\kappa_i))\mathcal B^{-1}$ is the nodal Dirichlet-to-Neumann map of
$\Delta-\kappa_i^2$, $\mathsf Q_i=\mathcal BC^{-1}\Phi_i^\top\operatorname{diag}(\omega_a)$ with $(\Phi_i)_{ak}=\phi_{i,k}(x_a)$ is the
volume quadrature of the modes $\phi_{i,k}$ brought to nodal form, and $\mathbf g(t)$, $\mathbf w(t)$, $\mathbf u_A(t)$, $\mathbf q(t)$
are the nodal traces and the values of $u$ and of $q=u^3-F$ at the area nodes. The three terms in the bracket are the source
moment, the terminal moment at $\tau_l$ and the initial moment, with $\mathbf u_A(0)=u_0$ at the area nodes. The time integrals use
$24$ Gauss nodes per window.

Row $j$ of \eqref{eq:app-heat-nodal} contains the flux only through $w_j(t)$, so, as on the disk, the rows of the nodes of $\GN$ are kept and those of $\GD$ discarded, and $n=12\times5\times14=840$. Inserting $\mathbf g_N(t)=\mathbf g_N(0)+Z\boldsymbol\psi(t)$, the quantities of the residual \eqref{eq:evo-residual} are, for the row $(i,l,j)$ and the column $(n',p)$ with $j,n'\in\GN$:
\begin{equation}
\label{eq:app-heat-blocks}
\begin{aligned}
 \mathsf K_{(i,l,j),(n',p)}&=(\mathsf D_i)_{jn'}\int_0^{\tau_l}e^{-\kappa_i^2t}\psi_p(t)\,\mathrm dt,
 \qquad
 (\mathcal M_{\rm nl})_{(i,l,j)}=\Bigl(\mathsf Q_i\int_0^{\tau_l}e^{-\kappa_i^2t}\,\mathbf q(t)\,\mathrm dt\Bigr)_j,\\
 (\mathcal M_{\rm end})_{(i,l,j)}&=e^{-\kappa_i^2\tau_l}\bigl(\mathsf Q_i\mathbf u_A(\tau_l)\bigr)_j,\\
 \mathbf b_{(i,l,j)}&=\int_0^{\tau_l}\!e^{-\kappa_i^2t}\Bigl(w_j(t)-\sum_{n'\in\GD}(\mathsf D_i)_{jn'}g_{n'}(t)
   -\sum_{n'\in\GN}(\mathsf D_i)_{jn'}g_{n'}(0)\Bigr)\mathrm dt+\bigl(\mathsf Q_i\mathbf u_A(0)\bigr)_j .
\end{aligned}
\end{equation}
The $840\times112$ system is solved as described in the introduction. The area nodes are $14$ Gauss radii times the $29$ boundary directions.

The same relation is used for the PINNacle \cite{hao2024pinnacle} heat benchmark and for the space-time operator, with the same nodes and windows. For the benchmark $\mathcal N=0$, so $\mathcal M_{\rm nl}$ depends neither on $z$ nor on $\theta$, $J_C=I$ in
\eqref{eq:evo-projected-newton}, and the completion is a single linear solve per evaluation.

\subsubsection{KdV equation on $(-1,1)\times(0,0.25]$}
\label{app:evo-st-kdv}

The KdV construction differs from the heat construction in the spatial operator, which is of third order, producing three traces per endpoint that enter the boundary term. The tests are the exponentials of the Fourier transform in $x$, and the nonlinearity is handled in the conservative form $u\,\partial_xu=\partial_x\!\left(\tfrac12u^2\right)$.

\paragraph{Boundary term and traces.} Here $\mathcal L_x=\partial_x^3$, $\mathcal N[u]=u\partial_xu$ and $\mathcal L_x^*=-\partial_x^3$.
Integrating $\varphi\,\partial_x^3u$ by parts three times on $(-1,1)$ gives the boundary form of \eqref{eq:prelim-evo-GR}:
\begin{equation}
\label{eq:app-kdv-concomitant}
 \int_{\bd}\mathscr B_{\mathcal L}(u,\varphi)\,\mathrm ds
   =-\Bigl[\varphi\,\partial_x^2u-\partial_x\varphi\,\partial_xu+\partial_x^2\varphi\,u\Bigr]_{x=-1}^{x=1},
\end{equation}
so the traces are $(u,u_x,u_{xx})$ at each endpoint, six in all. $\mathsf M_{\rm KdV}$ prescribes $u(\pm1,t)$ and $u_x(1,t)$, and
the unknowns are $u_x(-1,t)$, $u_{xx}(-1,t)$ and $u_{xx}(1,t)$, each expanded in the $M_t=8$ functions \eqref{eq:app-psi} as in
\eqref{eq:evo-trace-basis}, so $m=24$. All six traces are encoded by the lift, the quintic Hermite interpolant of the six endpoint
values in Table~\ref{tab:app-lifts}, so there is no auxiliary trace and $\mathsf T=I$.

\paragraph{Tests.} The tests are the plane waves of the Fourier transform in $x$ with the dispersion relation of the KdV
equation:
\begin{equation}
\label{eq:app-kdv-tests}
 \varphi_k(x,t)=e^{-\mathrm{i}kx-\mathrm{i}k^3t},\qquad k\in\{0.5,1,\ldots,5\},
\end{equation}
which satisfy $-\partial_t\varphi_k-\partial_x^3\varphi_k=0$, the backward adjoint equation, so $\varphi_k\in\mathcal H_T$. With
them the volume integrals become Fourier transforms, $\hat u(k,t)=\int_{-1}^1e^{-\mathrm{i}kx}u(x,t)\,\mathrm dx$, and the ten values
of $k$ sample the transform over the range of wavenumbers present in the solution.

\paragraph{The Global Relation.} Insert $\varphi_k$ with the upper limit $\tau_l$ in \eqref{eq:prelim-evo-GR}, using
\eqref{eq:app-kdv-concomitant} for the boundary term:
\begin{multline}
\label{eq:app-kdv-gr}
 \int_0^{\tau_l}e^{-\mathrm{i}k^3t}\Bigl[e^{-\mathrm{i}kx}\bigl(\partial_x^2u+\mathrm{i}k\,\partial_xu-k^2u\bigr)\Bigr]_{x=-1}^{x=1}\mathrm dt
 +\int_0^{\tau_l}\!\!\int_{-1}^{1}\varphi_k\,\bigl(u\partial_xu-F\bigr)\,\mathrm dx\,\mathrm dt\\
 +e^{-\mathrm{i}k^3\tau_l}\hat u(k,\tau_l)-\hat u_0(k)=0 .
\end{multline}
The bracket contains the six traces with explicit coefficients. Those of the unknown traces, $-\mathrm{i}k\,e^{\mathrm{i}k}$ for
$\partial_xu(-1,t)$, $-e^{\mathrm{i}k}$ for $\partial_x^2u(-1,t)$ and $e^{-\mathrm{i}k}$ for $\partial_x^2u(1,t)$, multiplied by the time
integrals of $e^{-\mathrm{i}k^3t}\psi_p(t)$, give $\mathsf K$. The prescribed traces and the initial moment $\hat u_0(k)$ give
$\mathbf b$, and the term $e^{-\mathrm{i}k^3\tau_l}\hat u(k,\tau_l)$ is the terminal moment $\mathcal M_{\rm end}$.

For the source moment the nonlinearity is written in conservative form, $u\partial_xu=\partial_x(u^2/2)$, and integrated by parts:
\begin{equation}
\label{eq:app-kdv-conservative}
 \int_{-1}^{1}\varphi_k\,u\partial_xu\,\mathrm dx
   =\Bigl[\tfrac12\varphi_k\,u^2\Bigr]_{x=-1}^{x=1}+\frac{\mathrm{i}k}{2}\int_{-1}^{1}\varphi_k\,u^2\,\mathrm dx .
\end{equation}
The endpoint term contains only the prescribed values $u(\pm1,t)$ and is moved to $\mathbf b$, and the source moment becomes
$(\mathcal M_{\rm nl})_{k,l}=\qint_{\Om\times(0,\tau_l)}\varphi_k\bigl(\tfrac{\mathrm{i}k}{2}u^2-F\bigr)\,\mathrm dx\,\mathrm dt$, which
involves $u$ but not $\partial_xu$ (the remark at the end of Section~\ref{sec:evolution}). Real and imaginary parts of the ten
tests over the five windows give $n=p=2\times 10\times 5=100$ equations for the $m=24$ unknowns, completed by rule (P) as described in the introduction. The quadrature uses $40$ Gauss nodes in $x$ and $48$ Gauss nodes per window in $t$.

\subsection{Time-stepping construction}
\label{app:evo-ts}

In the time-stepping construction each backward-Euler level is the stationary stage problem \eqref{eq:time-stepping-stage}, with
the stage source $F_\ell=F(\cdot,t_\ell)+u^{\ell-1}/\Delta t_\ell$ built from the previous level. The only new element is the term
$u^\ell/\Delta t_\ell$, and the two examples treat it in the two possible ways.

\subsubsection{Heat equation on $\square$ ($\varsigma_\ell=0$)}
\label{app:evo-ts-heat}

The term $u^\ell/\Delta t_\ell$ is kept in the source, $q^\ell=(u^\ell)^3+u^\ell/\Delta t_\ell-F_\ell$. The stage operator is then
$-\Delta$, the stage relation is the elliptic GR \eqref{eq:GR}, and the harmonic tests, trace, lift, quadrature and matrices of
the square construction of Appendix~\ref{app:stat-square} are used unchanged at every level ($m=64$, $n=128$, $\mathsf B'=0$). This
is the stage residual \eqref{eq:time-stepping-residual} with $\mathbf N(z)=(\mathbf u_A^\ell(z))^3$ and $\varsigma_\ell=0$, and the
only difference from a stationary problem is that $\mathbf F^\ell$ changes from level to level.

\subsubsection{KdV equation ($\varsigma_\ell=1$)}
\label{app:evo-ts-kdv}

The term $u^\ell/\Delta t_\ell$ is kept in the operator, so the stage operator is $\partial_x^3+1/\Delta t_\ell$ and the tests must
solve its adjoint equation $-\varphi'''+\varphi/\Delta t_\ell=0$. Its solutions are the exponentials:
\begin{equation}
\label{eq:app-kdv-level-tests}
 \varphi(x)=e^{\mu x},\qquad \mu^3=\frac1{\Delta t_\ell},\qquad \mu_j=\Delta t_\ell^{-1/3}e^{2\pi\mathrm{i}j/3},\quad j=0,1,2 ,
\end{equation}
one real root and a complex pair, used in the real form $\varphi_1=e^{\mu_0(x-1)}$, $\varphi_2=\operatorname{Re}e^{\mu_1(x+1)}$,
$\varphi_3=\operatorname{Im}e^{\mu_1(x+1)}$, each shifted to the endpoint where it is largest so that all three are bounded by one on $[-1,1]$. By \eqref{eq:app-kdv-concomitant} the stage relation for the test $\varphi_r$ is:
\begin{equation}\label{eq:app-kdv-level-gr}
 \Bigl[\varphi_r\,\partial_x^2u^\ell-\varphi_r'\,\partial_xu^\ell+\varphi_r''\,u^\ell\Bigr]_{x=-1}^{x=1}
 +\qint_{-1}^{1}\varphi_r\bigl(u^\ell\partial_xu^\ell-F_\ell\bigr)\,\mathrm dx=0,
 \qquad r=1,2,3.
\end{equation}
The traces are the six endpoint values of the space-time construction, now at the single time $t_\ell$, with the same three unknowns $z=(\partial_xu(-1),\partial_x^2u(-1),\partial_x^2u(1))$ at $t_\ell$ and the same quintic Hermite lift, so again $\mathsf T=I$. Three tests for three unknowns give a square system, $p=n=m=3$, with:
\begin{align*}
 \mathsf K&=\begin{pmatrix}\varphi_r'(-1)&-\varphi_r(-1)&\varphi_r(1)\end{pmatrix}_{r=1,2,3},\\
 \mathbf b_r&=\varphi_r''(-1)\,u(-1,t_\ell)-\varphi_r''(1)\,u(1,t_\ell)+\varphi_r'(1)\,\partial_xu(1,t_\ell),
\end{align*}
the rows of $\mathsf K$ being the coefficients of the three unknown traces in the bracket and $\mathbf b$ the contribution of the three prescribed ones. This is the stage residual \eqref{eq:time-stepping-residual} with $\varsigma_\ell=1$ and
$\mathbf N(z)=\mathbf u_A^\ell\,\partial_x\mathbf u_A^\ell$, solved by Newton's method (rule (R)). The source moment uses $64$ Gauss nodes,
and the gate is $(1-x^2)^3$, which vanishes at the endpoints together with its first two derivatives and therefore leaves all six traces to the lift.

\section{Error accumulation in the time-stepping construction}
\label{app:time-stepping-error}

We analyse how the stage errors of the time-stepping construction of Section~\ref{prop:evolution-time-stepping} accumulate over the time steps, in the notation of Section~\ref{sec:evolution}. In what follows, we consider a time grid $0=t_0<\cdots<t_{N_T}=\mathcal T$, steps $\Delta t_\ell=t_\ell-t_{\ell-1}$ and stage approximations $u^\ell$, with $u^0=u_0$. The constants and norms are those of Section~\ref{sec:consistency-error}.

\begin{proposition}[Accumulation of time-stepping errors]
\label{prop:ca-time-error}
Let $u^\dagger(\cdot,t)$ be the selected exact evolution
solution, and use the backward-Euler time levels of the
time-stepping construction.
For a previous stage approximation $w$, let $\mathcal S_\ell(w)$ be the
selected solution of:
\[ \mathcal L_xv+\mathcal N[v;\lambda]+\frac{v}{\Delta t_\ell}
= F(\cdot,t_\ell)+\frac{w}{\Delta t_\ell}, \]
with the prescribed boundary data at $t_\ell$.
Assume these maps exist on a set containing the relevant
exact and approximate inputs, and satisfy:
\begin{align}
 \|\mathcal S_\ell(w)-\mathcal S_\ell(\widetilde w)\|_\star
   &\le(1+C_{\rm st}\Delta t_\ell)
                         \|w-\widetilde w\|_\star,
       \label{eq:ca-time-stability}\\
 \|u^\dagger(\cdot,t_\ell)
       -\mathcal S_\ell(u^\dagger(\cdot,t_{\ell-1}))\|_\star
   &\le C_{\rm t}(\Delta t_\ell)^2,
       \label{eq:ca-time-consistency}
\end{align}
with $C_{\rm st},C_{\rm t}\ge0$ independent of the time grid.
Suppose the stage Fokas-PINN approximations satisfy:
\[
 \|u^\ell-\mathcal S_\ell(u^{\ell-1})\|_\star
       \le\varepsilon_\ell^{\rm stage}.
\]
Then, with
$E^\ell=\|u^\ell-u^\dagger(\cdot,t_\ell)\|_\star$, we have:
\begin{equation}
\label{eq:ca-time-recurrence}
 E^\ell
   \le(1+C_{\rm st}\Delta t_\ell)E^{\ell-1}
       +\varepsilon_\ell^{\rm stage}
       +C_{\rm t}(\Delta t_\ell)^2,
\end{equation}
and:
\begin{equation}
\label{eq:ca-time-global}
 E^\ell\le e^{C_{\rm st}t_\ell}
 \left[
  E^0+\sum_{j=1}^{\ell}\varepsilon_j^{\rm stage}
      +C_{\rm t}t_\ell\max_{1\le j\le\ell}\Delta t_j
 \right].
\end{equation}
When a stage satisfies the hypotheses of
Proposition~\ref{prop:ca-completion-error}, its error can
be bounded by \eqref{eq:ca-total-error}, applied to the
exact stage solution $\mathcal S_\ell(u^{\ell-1})$
with the same previous input. Consequently, if $E^0\to0$,
$\max_\ell\Delta t_\ell\to0$ and
$\sum_\ell\varepsilon_\ell^{\rm stage}\to0$,
then $\max_{0\le\ell\le N_T}E^\ell\to0$.
\end{proposition}

\begin{proof}
Add and subtract
$\mathcal S_\ell(u^{\ell-1})$ and
$\mathcal S_\ell(u^\dagger(\cdot,t_{\ell-1}))$.
Using the triangle inequality, stage error estimate,
\eqref{eq:ca-time-stability} and
\eqref{eq:ca-time-consistency} give:
\[
\begin{aligned}
 E^\ell
   &\le\|u^\ell-\mathcal S_\ell(u^{\ell-1})\|_\star
   +\|\mathcal S_\ell(u^{\ell-1})
             -\mathcal S_\ell(u^\dagger(\cdot,t_{\ell-1}))\|_\star
+\|\mathcal S_\ell(u^\dagger(\cdot,t_{\ell-1}))
             -u^\dagger(\cdot,t_\ell)\|_\star\\
   &\le\varepsilon_\ell^{\rm stage}
       +(1+C_{\rm st}\Delta t_\ell)E^{\ell-1}
       +C_{\rm t}(\Delta t_\ell)^2.
\end{aligned}
\]
Iterating this recurrence and using
$1+C_{\rm st}\Delta t_j\le e^{C_{\rm st}\Delta t_j}$
yields:
\[
 E^\ell\le e^{C_{\rm st}t_\ell}
 \left[
 E^0+\sum_{j=1}^{\ell}
       \bigl(\varepsilon_j^{\rm stage}
                 +C_{\rm t}(\Delta t_j)^2\bigr)
 \right].
\]
Since
$\sum_{j=1}^{\ell}(\Delta t_j)^2
 \le t_\ell\max_{1\le j\le\ell}\Delta t_j$,
this proves \eqref{eq:ca-time-global}.
The convergence assertion follows directly.
\end{proof}

The stage stability and temporal consistency assumptions
are properties of the time-discretized PDE, separate from
regularity of the finite completion.
For example, a twice continuously time-differentiable solution
has a backward-difference residual of order $\Delta t_\ell$.
If the stage response to a forcing perturbation is bounded
by $C\Delta t_\ell$ times that perturbation, this gives
\eqref{eq:ca-time-consistency}.
The stage error bounds must retain their dependence on
both the spatial discretization and $\Delta t_\ell$.
Taking $u^0=u_0$ gives $E^0=0$; stage errors of order
$(\Delta t_\ell)^2$ are sufficient to preserve a first-order
global temporal bound.

\section{Convergence of the stop-gradient output}
\label{app:implicit}

We prove the two statements following Proposition \ref{prop:fpinn-implicit-differentiation} and \eqref{eq:fpinn-ift-output}. Specifically, at an exact completion, $z_{\mathrm{out}}$ and its derivative $D^{\mathrm{AD}}_\theta z_{\mathrm{out}}=-\bar J_C^{-1}\partial_\theta\mathcal C(\bar z;\theta,d)$ coincide with $z^*$ and $D_\theta z^*$; and for numerical completions $\widetilde z^{(j)}\to z^*$ both converge to these exact values.

\begin{proof}
Throughout, $\theta$ and $d$ are fixed and suppressed from the
notation: $\mathcal C(z):=\mathcal C(z;\theta,d)$,
$J_C(z):=D_z\mathcal C(z;\theta,d)$,
$\partial_\theta\mathcal C(z):=\partial_\theta\mathcal C(z;\theta,d)$,
and $J_C^\star=J_C(z^*)$. We use the spectral norm $\|\cdot\|_2$ for
square matrices and the Frobenius norm $\|\cdot\|_F$ for
$m\times N_\theta$ matrices, so that $\|AB\|_F\le\|A\|_2\|B\|_F$.

\smallskip\noindent\textit{(i) Exact completion.}
If $\widetilde z=z^*$, then $\bar z=z^*$, $\bar J_C=J_C^\star$ and
$\mathcal C(\bar z)=0$, so the Newton correction in
\eqref{eq:fpinn-ift-output} vanishes and $z_{\mathrm{out}}=z^*$. By
the derivative formula after \eqref{eq:fpinn-ift-output},
$D^{\mathrm{AD}}_\theta z_{\mathrm{out}}
 =-(J_C^\star)^{-1}\partial_\theta\mathcal C(z^*)$,
which is $D_\theta z^*$ by \eqref{eq:fpinn-implicit-sensitivity}. Since
$z_{\mathrm{out}}$ and its attached derivative coincide with $z^*$ and
$D_\theta z^*$, the chain rule through
$u_\theta=L(z_{\mathrm{out}};d)+\beta v_\theta$ gives
\eqref{eq:fpinn-field-sensitivity}, and the chain rule through
$\mathcal J(\theta;z_{\mathrm{out}},d)$ gives
\eqref{eq:fpinn-chain-rule} at $z=z^*$, which equals
\eqref{eq:fpinn-adjoint-gradient} by the proof of
Proposition~\ref{prop:fpinn-implicit-differentiation}.

\smallskip\noindent\textit{(ii) Numerical completions.}
We proceed in four steps.

\smallskip\noindent\textit{Step 1: eventual nonsingularity and bounded
inverses.} Since $\mathcal C$ is continuously differentiable near $z^*$
and $\widetilde z^{(j)}\to z^*$,
\[
 \widetilde J_C^{(j)}=J_C(\widetilde z^{(j)})\to J_C^\star,\qquad
 \mathcal C(\widetilde z^{(j)})\to\mathcal C(z^*)=0,\qquad
 \partial_\theta\mathcal C(\widetilde z^{(j)})
   \to\partial_\theta\mathcal C(z^*).
\]
The set of nonsingular $m\times m$ matrices is open and matrix inversion
is continuous on it. Since $J_C^\star$ is nonsingular, there is $j_0$
such that $\widetilde J_C^{(j)}$ is nonsingular for all $j\ge j_0$ and
$(\widetilde J_C^{(j)})^{-1}\to(J_C^\star)^{-1}$. In particular,
\[
 M:=\sup_{j\ge j_0}\|(\widetilde J_C^{(j)})^{-1}\|_2<\infty .
\]
(Quantitatively, once
$\|\widetilde J_C^{(j)}-J_C^\star\|_2\le\tfrac12\|(J_C^\star)^{-1}\|_2^{-1}$,
a Neumann-series argument gives
$\|(\widetilde J_C^{(j)})^{-1}\|_2\le2\|(J_C^\star)^{-1}\|_2$.)

\smallskip\noindent\textit{Step 2: forward outputs.}
For $j\ge j_0$, subtracting $z^*$ from
\eqref{eq:fpinn-ift-output} and using Step~1,
\[
 \|z^{(j)}_{\mathrm{out}}-z^*\|_2
 \le\|\widetilde z^{(j)}-z^*\|_2
   +\|(\widetilde J_C^{(j)})^{-1}\|_2\,\|\mathcal C(\widetilde z^{(j)})\|_2
 \le\|\widetilde z^{(j)}-z^*\|_2
   +M\,\|\mathcal C(\widetilde z^{(j)})\|_2
 \;\longrightarrow\;0 .
\]

\smallskip\noindent\textit{Step 3: attached sensitivities.}
Write $A_j:=(\widetilde J_C^{(j)})^{-1}$, $A:=(J_C^\star)^{-1}$,
$B_j:=\partial_\theta\mathcal C(\widetilde z^{(j)})$ and
$B:=\partial_\theta\mathcal C(z^*)$, so that
$D^{\mathrm{AD}}_\theta z^{(j)}_{\mathrm{out}}=-A_jB_j$ by
the derivative formula after \eqref{eq:fpinn-ift-output} and $D_\theta z^*=-AB$ by
\eqref{eq:fpinn-implicit-sensitivity}. Adding and subtracting $A_jB$,
\[
 A_jB_j-AB=A_j(B_j-B)+(A_j-A)B,
\]
hence, for $j\ge j_0$,
\[
 \|D^{\mathrm{AD}}_\theta z^{(j)}_{\mathrm{out}}-D_\theta z^*\|_F
 \le M\,\|B_j-B\|_F+\|A_j-A\|_2\,\|B\|_F
 \;\longrightarrow\;0,
\]
since $B_j\to B$ and $A_j\to A$ by Step~1.

\smallskip\noindent\textit{Step 4: approximation and loss.}
The AD derivative of
$u^{(j)}_\theta:=L(z^{(j)}_{\mathrm{out}};d)+\beta v_\theta$ is
$\beta D_\theta v_\theta+D_zL\,D^{\mathrm{AD}}_\theta z^{(j)}_{\mathrm{out}}$,
with $D_zL$ constant because $L$ is affine in $z$; by Step~3 it
converges to \eqref{eq:fpinn-field-sensitivity}. The gradient returned
by backpropagation through $\mathcal J(\theta;z^{(j)}_{\mathrm{out}},d)$
is
\[
 \partial_\theta\mathcal J(\theta;z^{(j)}_{\mathrm{out}},d)
 +\bigl(D^{\mathrm{AD}}_\theta z^{(j)}_{\mathrm{out}}\bigr)^{\top}
  \nabla_z\mathcal J(\theta;z^{(j)}_{\mathrm{out}},d).
\]
Since $\mathcal J$ is continuously differentiable in $(z,\theta)$, both
$\partial_\theta\mathcal J$ and $\nabla_z\mathcal J$ are continuous in
$z$, so by Step~2 they converge to their values at $z^*$; combined with
Step~3, the expression converges to \eqref{eq:fpinn-chain-rule} at
$z=z^*$, i.e.\ to \eqref{eq:fpinn-adjoint-gradient}.
\end{proof}

\bibliographystyle{ieeetr}
\bibliography{sc_fokas}
\end{document}